\documentclass[aop,MSNbibl,seceqn,citesort,dvips]{arximspdf}
\usepackage{mathbh}

\doi{10.1214/11-AOP680}
\volume{40}
\issue{6}
\pubyear{2012}
\firstpage{2299}
\lastpage{2361}

\makeatletter

\newtheorem{Lemma}{Lemma}[section]
\newtheorem{Proposition}[Lemma]{Proposition}
\newtheorem{Theorem}[Lemma]{Theorem}

\newproclaim{Remark}[Lemma]{Remark}

\newtheorem{Corollary}[Lemma]{Corollary}

\newcommand{\prob}{\mathbb{P}}
\newcommand{\bfw}{\mathbf{w}}
\newcommand{\bfc}{\mathbf{c}}
\newcommand{\bfd}{\mathbf{d}}
\newcommand{\E}{\mathbb{E}}
\newcommand{\CC}{\mathcal{C}}
\newcommand{\II}{\mathcal{I}}
\newcommand{\GG}{\mathcal{G}}
\newcommand{\SSS}{\mathcal{S}}
\newcommand{\ZZ}{\mathcal{Z}}

\newcommand{\Rbold}{{\mathbb{R}}}
\newcommand{\bfitx}{{\mathbf x}}
\newcommand{\RR}{\mathcal{R}}
\newcommand{\expec}{\mathbb{E}}

\newcommand{\Var}{\operatorname{Var}}
\newcommand{\Op}{O_{\prob}}
\newcommand{\cluster}{\mathcal{C}}
\newcommand{\weight}{\mathcal{W}}
\newcommand{\Ccal}{\mathcal{C}}
\newcommand{\Cmax}{\cluster_{\max}}
\newcommand{\Poi}{\operatorname{Poi}}
\newcommand{\Exp}{\operatorname{Exp}}
\newcommand{\vep}{\varepsilon}
\newcommand{\indic}[1]{\mathbh{1}_{\{#1\}}}
\newcommand{\indicwo}[1]{\mathbh{1}_{#1}}
\newcommand{\convd}{\stackrel{d}{\longrightarrow}}
\newcommand{\convp}{\stackrel{{\mathbb P}}{\longrightarrow}}
\newcommand{\conn}{\longleftrightarrow}
\newcommand{\nc}{\conn{\hspace{-2.8ex} /}\hspace{1.8ex}}
\newcommand{\ZBP}{Z^{(\mathrm{BP})}}
\newcommand{\XBP}{X^{(\mathrm{BP})}}
\newcommand{\varthetas}{\vartheta}
\renewcommand{\i}{{\mathrm i}}

\newcommand{\bfwit}{\mathbf{w}}
\newcommand{\Ver}{V_n}

\makeatother

\begin{document}
\begin{frontmatter}

\title{Novel scaling limits for critical inhomogeneous random graphs}
\runtitle{Novel scaling limits}

\begin{aug}
\author[A]{\fnms{Shankar} \snm{Bhamidi}\thanksref{t1}\ead[label=e1]{bhamidi@email.unc.edu}},
\author[B]{\fnms{Remco} \snm{van der Hofstad}\corref{}\thanksref{t2}\ead[label=e2]{rhofstad@win.tue.nl}}\\ and
\author[B]{\fnms{Johan S. H.} \snm{van Leeuwaarden}\thanksref{t2}\ead[label=e3]{j.s.h.v.leeuwaarden@tue.nl}}
\runauthor{S. Bhamidi, R. van der Hofstad and J. S. H. van Leeuwaarden}
\affiliation{University of North Carolina, Eindhoven University of
Technology and~Eindhoven University of Technology}
\address[A]{S. Bhamidi\\
Department of Statistics\\
\quad and Operations Research\\
University of North Carolina\\
304 Hanes Hall\\
Chapel Hill, North Carolina 27599\\
USA\\
\printead{e1}} 
\address[B]{R. van der Hofstad\\
J. S. H. van Leeuwaarden\\
Department of Mathematics\\
\quad and Computer Science\\
Eindhoven University of Technology\\
P.O. Box 513\\
5600 MB Eindhoven\\
The Netherlands\\
\printead{e2}\\
\hphantom{E-mail: }\printead*{e3}}
\end{aug}

\thankstext{t1}{Supported in part by PIMS and
NSERC, Canada and NSF Grant DMS-11-05581.}

\thankstext{t2}{Supported in part by the Netherlands
Organisation for Scientific Research (NWO).}

\received{\smonth{10} \syear{2010}}
\revised{\smonth{5} \syear{2011}}

%
\begin{abstract}
We find scaling limits for the sizes of the largest components at
criticality for rank-1 inhomogeneous random graphs with power-law
degrees with power-law exponent $\tau$. We investigate the case where
$\tau\in(3,4)$, so that the degrees have finite variance but infinite
third moment. The sizes of the largest clusters, rescaled by
$n^{-(\tau-2)/(\tau-1)}$, converge to hitting times of a ``thinned''
L\'evy process, a special case of the general multiplicative
coalescents studied by Aldous [\textit{Ann. Probab.} \textbf{25}
(1997) 812--854] and Aldous and Limic [\textit{Electron. J. Probab.}
\textbf{3} (1998) 1--59].

Our results should be contrasted to the case $\tau>4$, so that the
third moment is finite. There, instead, the sizes of the components
rescaled by $n^{-2/3}$ converge to the excursion lengths of an
inhomogeneous Brownian motion, as proved in Aldous [\textit{Ann.
Probab.} \textbf{25} (1997) 812--854] for the
Erd\H{o}s--R\'enyi random graph
and extended to the present setting in Bhamidi, van der Hofstad and
van Leeuwaarden [\textit{Electron. J. Probab.} \textbf{15} (2010)
1682--1703] and Turova [(2009) Preprint].
\end{abstract}

%
\begin{keyword}[class=AMS]
\kwd{60C05}
\kwd{05C80}
\kwd{90B15}.
\end{keyword}
\begin{keyword}
\kwd{Critical random graphs}
\kwd{phase transitions}
\kwd{inhomogeneous networks}
\kwd{thinned L\'evy processes}
\kwd{multiplicative coalescent}.
\end{keyword}

\pdfkeywords{60C05, 05C80, 90B15, Critical random graphs,
phase transitions, inhomogeneous networks, thinned Levy processes,
multiplicative coalescent}

\end{frontmatter}

\section{Introduction}
\label{sec-int}

The critical behavior of random graphs has received tremendous
attention in the past decades. The simplest example of a random graph
is the Erd\H{o}s--R\'enyi random graph, whose critical regime has been
intensely explored (see, e.g.,
\cite{Aldo97,AloSpe00,Boll01,Durr06,JanLucRuc00} and the references
therein). In the past few years, many examples of real-world networks have
been found where the degrees are highly variable and heavy tailed,
unlike the degrees in the Erd\H{o}s--R\'enyi random graph, which
instead are extremely light tailed. As a result, there has been a
concerted effort to define and analyze models for such real-world
networks. See, for example,~\cite{AlbBar02,DorMen02,Newm03a} for major
reviews of real-world networks and models for them.

In this paper, we study how \textit{inhomogeneity} in the random graph model
changes the critical regime of the random graph. In our model,
the vertices have a \textit{weight} associated to them,
and the weight of a vertex moderates its degree. Therefore, by choosing
these weights appropriately, we can generate random graphs with highly
variable degrees. For our class of random graphs,
it is shown in~\cite{Hofs09a}, Theorem 1.1, that when the weights do
not vary
too much, the critical behavior is similar to the one in the
Erd\H{o}s--R\'enyi random graph. See in particular the recent
works~\cite{BhaHofLee09a,Turo09}, where it was shown that if the
degrees have
finite \textit{third} moment, then the scaling limit for the largest critical
components in the critical window are essentially the same as for the
Erd\H{o}s--R\'enyi random graph, as identified by Aldous in~\cite{Aldo97}.

Interestingly, in~\cite{Hofs09a}, Theorem 1.2, it was shown that when
the degrees
have \textit{infinite} third moment, then the sizes of the largest
critical clusters are quite \textit{different}. See also~\cite{HatMol09}
for a related result for the configuration model, another
random graph model having flexible degrees. In this paper, we
identify the \textit{scaling limits} of the
largest critical clusters in the critical window in the regime where
the degrees have infinite third moments. As we shall see, this scaling limit
is rather different compared to that for the
Erd\H{o}s--R\'enyi random graph.
Let us first introduce the model that shall  be the focus of our
investigations for the rest of this article.

\subsection{Model}
\label{sec-mod}

In our random graph model, vertices have \textit{weights},
and the edges are independent with the edge probability
being approximately equal to the rescaled product of the weights
of the two end vertices of the edge. While there are many
different versions of such random graphs (see below),
it will be convenient for us to work
with the so-called \textit{Poissonian random graph} or Norros--Reittu model
\cite{NorRei06}. To define the model, we consider the vertex set
$[n]:=\{1,2,\ldots, n\}$ and suppose each vertex is assigned a \textit
{weight}, vertex
$i$ having weight $w_i$. Now, attach an edge with probability
$p_{ij}$ between vertices $i$ and $j$, where
%
%
\begin{equation} \label{pij-NR} p_{ij} = 1-\exp(-w_i
w_j/\ell_n)
\end{equation}
with $\ell_n$ denoting the \textit{total weight}
%
%
\begin{equation}
\label{ln-def} \ell_n = \sum_{i\in[n]} w_i.
\end{equation}
Different edges are independent. In this model, the average degree
of vertex $i$ is close to $w_i$, which brings \textit{inhomogeneity}
into the model.

There are many adaptations of this model,
for which equivalent results hold. Indeed,
the model considered here is a special case of the
so-called \textit{rank-1 inhomogeneous random graph}
introduced in great generality by Bollob\'as,
Janson and Riordan~\cite{BolJanRio07}. It is
asymptotically equivalent with many related models,
such as the \textit{random graph with given prescribed degrees}
or the Chung--Lu model, where instead
%
%
\begin{equation}
p_{ij}=\max(w_iw_j/\ell_n, 1),
\end{equation}
and which has been studied intensively by Chung and Lu (see
\cite{ChuLu02a,ChuLu02b,ChuLu03,ChuLu06c,ChuLu06}).
A~further adaptation is the \textit{generalized random graph}
introduced by Britton, Deijfen and Martin-L\"of in
\cite{BriDeiMar-Lof05}, for which
%
%
\begin{equation} \label{pij-GRG} p_{ij} = \frac{w_i w_j}{\ell
_n+w_iw_j}.
\end{equation}
See Janson~\cite{Jans08a} for conditions under which these random graphs
are \textit{asymptotically equivalent}, meaning that all events
have asymptotically equal probabilities. As discussed in more detail
in~\cite{Hofs09a}, Section 1.3, these conditions apply
in the setting to be studied in this paper. Therefore, all results
proved here also hold for these related rank-1 models.

Let us now specify how the weights are chosen.
We let the weight sequence $\bfw= (w_i)_{i\in[n]}$ be defined by
%
%
\begin{equation} \label{choicewi} w_i = [1-F]^{-1}(i/n),
\end{equation}
where $F$ is a distribution function on $[0,\infty)$
for which we assume that there exists
a $\tau\in(3,4)$ and $0<c_{ F}<\infty$ such that
%
%
\begin{equation} \label{F-bound-tau3,4b} \lim_{x\rightarrow\infty
}x^{\tau-1}[1-F(x)]= c_{ F},
\end{equation}
and where $[1-F]^{-1}$ is the generalized inverse function of
$1-F$ defined, for $u\in(0,1)$, by
%
%
\begin{equation} \label{invverd} [1-F]^{-1}(u)=\inf\{ s\dvtx
[1-F](s)\leq u\}.
\end{equation}
By convention, we set $[1-F]^{-1}(1)=0$. We often make use of the fact
that, with $U$ uniform on $[0,1]$,
the random variable $[1-F]^{-1}(U)$ has distribution function~$F$.

An interpretation of the choice in~(\ref{choicewi}) is that the weight of
a vertex $\Ver$ chosen uniformly in $[n]$ has distribution function $F_n$
given by
%
%
\begin{eqnarray} \label{Fn-def}
F_n(x)&=&\prob(w_{\Ver}\leq x)=\frac
1n\sum_{j\in[n]} \indic{w_j\leq x}\nonumber\\
&=&\frac1n\sum_{j\in[n]} \indic
{[1-F]^{-1}({j/n})\leq x}
=\frac1n\sum_{i=0}^{n-1} \indic
{[1-F]^{-1}(1-{i/n})\leq x}\nonumber\\[-8pt]\\[-8pt]
&=&\frac1n\sum
_{i=0}^{n-1}\indic{F^{-1}({i/n}) \leq x}
=\frac1n\sum
_{i=0}^{n-1} \indic{{i/n}\leq F(x)} \nonumber\\
&=&\frac1n\bigl(\lfloor
n F(x)\rfloor+1\bigr)\wedge1,\nonumber
\end{eqnarray}
where, throughout this paper and for $x,y\in\mathbb{R}$, we write
$(x\vee y)=\max(x,y)$ and
$(x \wedge y)=\min(x,y)$. By~(\ref{Fn-def}), $F_n\rightarrow F$
uniformly. As a result, a uniformly
chosen vertex has a weight which is close in distributional sense to $F$.

For the setting in~(\ref{pij-NR}) and~(\ref{choicewi}), by
\cite{BolJanRio07}, Theorem 3.13, the
number of vertices with degree $k$, denoted by $N_k$,
satisfies
%
%
\begin{equation} \label{N-k-conv}
N_k/n\convp\expec\biggl[\mathrm
{e}^{-W} \frac{W^k}{k!}\biggr],\qquad k\geq0,
\end{equation}
where $\convp$ denotes convergence in probability, and where $W$ has
distribution function $F$ appearing in~(\ref{choicewi}). We recognize
the limiting distribution as a \textit{mixed Poisson
distribution with mixing distribution $F$}; that is, conditionally on
$W=w$, the distribution is Poisson with mean $w$. Equation
(\ref{N-k-conv}) also implies that the distribution of the degree of a
uniformly chosen vertex in $[n]$ converges to a mixed Poisson
distribution with mixing distribution $F$. This can be understood by
noting that the weight of a uniformly chosen vertex is, by
(\ref{Fn-def}), close in distribution to $F$. In turn, when a vertex
has weight $w$, then, by~(\ref{pij-NR}), its degree is close to Poisson
with parameter $w$. Since a Poisson random variable with large
parameter $w$ is closely concentrated around its mean $w$, we see that
the tail behavior of the degrees in our random graph is close to that
of the distribution $F$. As a result, when~(\ref{F-bound-tau3,4b})
holds, and with $D_n$ the degree of a uniformly chosen vertex in~$[n]$,
$\limsup_{n\rightarrow\infty} \expec [D^a_n]<\infty$ when $a<\tau-1$
and $\limsup_{n\rightarrow\infty} \expec [D^a_n]=\infty$ when
$a\geq\tau-1$. In particular, the degree of a uniformly chosen vertex
in $[n]$ has finite second, but infinite third moment when
(\ref{F-bound-tau3,4b}) holds with $\tau\in(3,4)$.

We shall frequently make use of the fact that
(\ref{F-bound-tau3,4b}) implies that, as $u\downarrow0$,
%
%
\begin{equation} \label{[1-F]inv-bd} [1-F]^{-1}(u)=(c_{
F}/u)^{1/(\tau-1)}\bigl(1+o(1)\bigr).
\end{equation}
Under the key assumption in~(\ref{F-bound-tau3,4b}), we have that
the third moment of the degrees tends to infinity; that is,
with $W\sim F$, $\expec[W^3]=\infty$. Define
%
%
\begin{equation} \label{nu-def} \nu= \expec[W^2]/\expec[W],
\end{equation}
so that, again by~(\ref{F-bound-tau3,4b}), $\nu<\infty$.
Then, by~\cite{BolJanRio07}, Theorem 3.1 (see also
\cite{BolJanRio07}, Section~16.4, for a detailed discussion
on rank-1 inhomogeneous random graphs,
of which our random graph is an example),
when $\nu> 1$, there is one giant component of size proportional to $n$,
while all other components are of smaller size $o(n)$, and when
$\nu\leq1$, the largest connected component contains a proportion of
vertices that converges to zero in probability. Thus, the critical value
of the model is $\nu=1$. The main aim of this paper is to investigate what
happens \textit{close to} the critical point, that is, when $\nu=1$.

A simple example of our model arises when we take
%
%
\begin{equation} \label{F-example}
F(x)=\cases{ 0, &\quad for $x<a$,\cr
1-(a/x)^{\tau-1}, &\quad for $x\geq a$,}
\end{equation}
in which case $[1-F]^{-1}(u)= a(1/u)^{1/(\tau-1)}$, so
that $w_j=a(n/j)^{1/(\tau-1)}$ and
%
%
\begin{equation} \expec[W]=\frac{a(\tau-1)}{\tau-2},\qquad \expec
[W^2]=\frac{a^2(\tau-1)}{\tau-3}.
\end{equation}
The critical case thus arises when
%
%
\begin{equation} \nu=\expec[W^2]/\expec[W]=\frac{a(\tau-2)}{\tau
-3}=1,
\end{equation}
that is, when $a=(\tau-3)/(\tau-2)$.

With the definition of the weights in~(\ref{choicewi}),
we shall write $\GG_n^0(\bfw)$ for the
graph constructed with the probabilities in~(\ref{pij-NR}), while, for
any fixed $\lambda\in\Rbold$, we
shall write $\GG_n^\lambda(\bfw)$ when we use the weight sequence
$\bfwit(\lambda)=(w_i(\lambda))_{i\in[n]}$ defined by
%
%
\begin{equation} \label{w-lambda-def} \bfwit(\lambda)=\bigl(1+\lambda
n^{-(\tau-3)/(\tau-1)}\bigr)\bfw.
\end{equation}

We shall assume that $n$ is so large that $1+\lambda n^{-(\tau
-3)/(\tau-1)}\geq0$,
so that $w_i(\lambda)\geq0$ for all $i\in[n]$.
This setting was first explored in~\cite{Hofs09a},
where, for the largest connected component $\Cmax$ and all
$\lambda\in{\mathbb R}$, it is proved that both
$n^{-(\tau-2)/(\tau-1)}|\Cmax|$ and $n^{(\tau-2)/(\tau-1)}/|\Cmax|$
are \textit{tight} sequences of random variables. In this paper, we
bring the discussion of the critical behavior of such inhomogeneous
random graphs
substantially further, by identifying the scaling limit of
$(n^{-(\tau-2)/(\tau-1)}|\cluster_{(i)}|)_{i\geq1}$,
where $(\cluster_{(i)})_{i\geq1}$
denote the connected components ordered in size, that is, $|\Cmax
|=|\cluster_{(1)}|\geq|\cluster_{(2)}|\geq\cdots.$

Interestingly, as proved in
\cite{BhaHofLee09a,Hofs09a,Turo09}, when $\tau>4$, so that
$\expec[W^3]<\infty$, the scaling limits of the random graphs studied here
are (apart from a trivial scaling constant) \textit{equal} to the
scaling limit
of the ordered connected components in the Erd\H{o}s--R\'enyi random
graph, as
first identified by Aldous in~\cite{Aldo97}. This suggests that the
\textit{high-weight vertices} play a crucial role in our setting,
a fact that shall feature extensively throughout our proof.
The importance of the high-weight vertices also partly explains why
we restrict our setting to~(\ref{choicewi}) and~(\ref{F-bound-tau3,4b}),
which give us sharp asymptotics of the weights of the high-weight vertices
in the heavy-tailed setting we study here.
We shall comment on extensions of our results in more detail in
Section~\ref{sec-disc3,4} below.

Before stating our main results, we introduce some notation. For a vertex
$i\in[n]$, we write $\cluster(i)$ for the vertices in the connected component
or \textit{cluster} of $i$. Further, let
%
%
\begin{equation} \Ccal_{\leq}(i) =
\cases{\Ccal(i), &\quad if $i\leq j$ $\forall j\in\Ccal(i)$,\cr
\varnothing, &\quad otherwise.}
\end{equation}
Then, clearly, $|\Cmax|=\max_{i\in[n]} |\Ccal(i)|=\max_{i\in[n]}
|\Ccal_{\leq}(i)|$,
and $(|\cluster_{(i)}|)_{i\geq1}$ is equal to the sequence
$(|\Ccal_{\leq}(i)|)_{i\geq1}$ ordered in size.
We further define the \textit{cluster weight}
of vertex $i$ to be
%
%
\begin{equation} \label{weight-cluster-def} \weight(i)=\sum_{j\in
\cluster(i)} w_j,
\end{equation}
and let $\weight_{\leq}(i)$ be as in~(\ref{weight-cluster-def}),
where the sum is restricted to $\cluster_{\leq}(i)$.
We again let $(\weight_{(i)})_{i\geq1}$ be equal to the sequence
$(\weight_{\leq}(i))_{i\geq1}$ ordered in size.

Throughout this paper, we shall make use of the following standard notation.
We let $\convd$ denote convergence in distribution, and
$\convp$ convergence in probability. For a sequence of random variables
$(X_n)_{n\geq1}$, we write $X_n=\Op(b_n)$ when $|X_n|/b_n$
is a tight sequence of random variables as $n\rightarrow\infty$,
and $X_n=o_{\prob}(b_n)$ when $|X_n|/b_n\convp0$ as
$n\rightarrow\infty$.
For a nonnegative function $n\mapsto g(n)$,
we write $f(n)=O(g(n))$ when $|f(n)|/g(n)$ is uniformly bounded, and
$f(n)=o(g(n))$ when $\lim_{n\rightarrow\infty} f(n)/g(n)=0$.
Furthermore, we write $f(n)=\Theta(g(n))$ if $f(n)=O(g(n))$ and $g(n)=O(f(n))$.
Finally, we write that a sequence of events $(E_n)_{n\geq1}$
occurs \textit{with high probability} (\textit{whp}) when $\prob
(E_n)\rightarrow1$.

Now we are ready to state our main results. We start in Section \ref
{sec-tau3,4}
by describing the scaling limit of the ordered clusters, and in Section
\ref{sec-furth-res} we discuss further properties of the scaling limit.

\subsection{\texorpdfstring{The scaling limit for $\tau\in(3,4)$}{The scaling limit for tau in (3,4)}}
\label{sec-tau3,4}

In this section, we investigate the scaling limit of
the connected components ordered in size.
Our first main result is as follows:
\begin{Theorem}[{[Weak convergence of the ordered critical clusters for
$\tau\in(3,4)$]}]
\label{thm-WC-3,4}
Fix the Norros--Reittu random graph with weights
$\bfw(\lambda)$ defined in~(\ref{choicewi}) and~(\ref{w-lambda-def}).
Assume that $\nu=1$ and that~(\ref{F-bound-tau3,4b}) holds. Then,
for all $\lambda\in\Rbold$,
%
%
\begin{equation}
\bigl(n^{-(\tau-2)/(\tau-1)}\bigl|\cluster_{
(i)}\bigr|\bigr)_{i\geq1} \convd(\gamma_i(\lambda))_{i\geq1}
\end{equation}
in the product topology, for some nondegenerate limit $(\gamma
_i(\lambda))_{i\geq1}$.
\end{Theorem}

We next study the \textit{joint} convergence of the clusters for different
values of \mbox{$\lambda\in\mathbb{R}$}. By increasing $\lambda$,
more and more edges are added to the system. These extra edges potentially
create connections between \textit{disjoint} clusters, thus merging them.
As a result, we can interpret $\lambda$ as a \textit{time variable},
and as time increases, clusters are being merged. This resembles a
\textit{coalescence process}, as studied in~\cite{Bert06}. We now make
this connection
precise. Before being able to do so, we introduce some necessary notation.

We first give a quick overview of Aldous's standard multiplicative
coalescent and how it relates to the limiting random variables in
Theorem~\ref{thm-WC-3,4}, seen as functions of the parameter
$\lambda$.
It will not be possible to give a full description of the process and its
many fascinating properties here, and we refer the interested reader to
the paper~\cite{AldLim98}, the survey paper~\cite{Aldo99} and the
book~\cite{Bert06}.

Write $\ell^2_{\searrow}$ for the metric space of infinite real-valued
sequences $\bfitx= (x_1,\break x_2,\ldots)$ with $x_1\geq x_2\geq\cdots
\geq0$
and $\sum_{i=1}^{\infty} x_i^2<\infty$,
with the $\ell^2$-norm as the metric.
The standard multiplicative coalescent is described as the Markov
process with
states in $\ell_{\searrow}^2$ whose dynamics is as follows: for
each pair of clusters $(x,y)$, the pair merges at rate $xy$. Thus, the
multiplicative coalescent is a continuous-time Markov process of
the masses of an infinite number of particles, where two particles merge
at a rate equal to the \textit{product} of their masses.

In~\cite{Aldo97}, Aldous showed that there is a Feller process on the
space $\ell^2_{\searrow}$ defined for all times
$-\infty< t < \infty$ starting from
infinitesimally small masses at time $-\infty$, and following the above
merging dynamics. The distribution of the coalescent
process at any time $t$ is the same as the limiting ordered cluster sizes
of an Erd\H{o}s--R\'enyi random graph with edge probabilities
$p_n = (1 + tn^{-1/3})/n$.

Aldous and Limic~\cite{AldLim98} explicitly characterize the entrance
boundary at $-\infty$ of the above Markov process, in the sense that
they prove that every \textit{extreme} version of the above Markov process
is characterized by a diffusion parameter $\kappa$, a translation
parameter $\beta$, and a vector of ``limiting largest weights''
$\bfc= (c_1, c_2, \ldots)$ that describe the asymptotic decay of
the masses of the particles at time $-\infty$. In this terminology,
the multiplicative coalescent can
be described as the ordered lengths of excursions beyond past minima of
the process
%
%
\begin{equation} W^{\kappa, \beta, \bfc}(s)=\kappa^{1/2}W(s)+\beta
s-\tfrac12\kappa s^2+V^{\bfc}(s),
\end{equation}
where $(W(s))_{s\geq0}$ is a standard Brownian motion, while
%
%
\begin{equation} V^{\bfc}(s)=\sum_{j=1}^{\infty} c_j \bigl(\indic
{E_j\leq s}-c_js\bigr)
\end{equation}
with $(E_j)_{j\geq1}$ independent exponential random variables, $E_j$
having mean $1/c_j$. Then, the $(\kappa, \beta, \bfc
)$-multiplicative coalescent is the set of ordered lengths of
excursions from zero of
the reflected process
%
%
\begin{equation} B^{\kappa, \beta, \bfc}(s)=W^{\kappa, \beta, \bfc
}(s)-\min_{0\leq s'\leq s} W^{\kappa, \beta, \bfc}(s').
\end{equation}
Part of the proof in~\cite{AldLim98} is the fact that these ordered
excursions can be defined properly.

The following theorem draws a connection between the components of the
graph for a
fixed $\lambda$ and the sizes of clusters at the same time in a
multiplicative coalescent with a particular entrance boundary, scale and
translation parameter. For this, define the sequence
%
%
\begin{equation} \label{bfc-def} {\mathbf c} = \bigl(c i^{-1/{(\tau
-1)}}\bigr)_{i\geq1}\qquad \mbox{with } c=c_{ F}^{1/(\tau
-1)}.
\end{equation}
Then, we have the following theorem:
\begin{Theorem}[(Relation to multiplicative coalescents)]
\label{thm-mult-coal}
Assume that the conditions in Theorem~\ref{thm-WC-3,4} hold.
Consider the
sequence-valued random variables ${\mathbf X}^*(\lambda) = (\gamma
_1(\expec[W]\lambda), \gamma_2(\expec[W]\lambda), \ldots)$ with
$(\gamma_i(\lambda))_{i\geq1}$ as in Theorem~\ref{thm-WC-3,4}.
Then ${\mathbf X}^*(\lambda)$ has the same distribution as a
multiplicative coalescent at time $\lambda$ with entrance boundary
$\bfc/\expec[W]$, diffusion constant $\kappa=0$ and centering
constant $\beta=-\zeta/\expec[W]$, where $\zeta$ is identified
explicitly in~(\ref{zeta-form}) below.
More precisely, there exists a simultaneous coupling of the clusters
$(|\cluster_{(i)}(\lambda)|)_{i\geq1}$, where
$|\cluster_{(i)}(\lambda)|$ is the $i$th largest cluster when
the weights
are equal to $\bfwit(\lambda)$, such that,
for every vector $(\lambda_1,\lambda_2,\ldots,\lambda_k)$,
%
%
\begin{equation} \bigl(n^{-(\tau-2)/(\tau-1)}\bigl(\bigl|\cluster_{
(i)}(\lambda_l)\bigr|\bigr)_{i\geq 1}\bigr)_{l=1}^k \convd\bigl({\mathbf X}^*(\lambda
_l/\expec[W])\bigr)_{l=1}^k.
\end{equation}
In particular, with $c_j$ defined as in~(\ref{bfc-def}),
%
%
\begin{equation} |\lambda| \gamma_j(\lambda)\convp c_j
\qquad\mbox{as } \lambda\to-\infty\mbox{ for each } j\geq1.
\end{equation}
\end{Theorem}

Theorem~\ref{thm-mult-coal} proves that the finite-dimensional
distributions of the rescaled cluster sizes converge to those
of a multiplicative coalescent. While we believe that also
\textit{process} convergence holds, viewing the processes
as elements of an appropriate function space,
we have no proof for this fact.
See Section~\ref{secmult-coal} for a full proof
of Theorem~\ref{thm-mult-coal}. The setting in this paper is
the first example where the multiplicative coalescent
with $\kappa=0$ arises in random graph
theory. Indeed, all random graph examples in
\cite{AldLim98} have largest component sizes
of the order $n^{2/3}$, like for the Erd\H{o}s--R\'enyi random graph
studied in
\cite{Aldo97}. Our example links the multiplicative coalescent also to
random graphs
with the largest critical connected components of the order $n^{(\tau
-2)/(\tau-1)}$
instead of $n^{2/3}$.

A crucial part of the proof of Theorem~\ref{thm-mult-coal}
is the analysis of the \textit{subcritical} phase of our model.
The asymptotics of the rescaled ordered cluster sizes
in the subcritical regime acts as the \textit{entrance boundary} of
the multiplicative coalescent, as explained in more detail
in~\cite{AldLim98}, Proposition 7.
This entrance boundary is identified in the following theorem, which is
of independent interest. In the statement of Theorem~\ref{thm-sub-crit},
the lower bound on $\lambda_n$ appears only
to ensure that $w_i(\lambda_n)=(1+\lambda_n n^{-(\tau-3)/(\tau
-1)})w_i\geq0$
for every $i\in[n]$.
\begin{Theorem}[(Subcritical phase)]
\label{thm-sub-crit}
Assume that the conditions in Theorem~\ref{thm-WC-3,4} hold, but now
take $\lambda=\lambda_n\rightarrow-\infty$
as $n\rightarrow\infty$ such that $\lambda_n\geq-n^{-(\tau
-3)/(\tau-1)}$. Then, for each $j\in{\mathbb N}$,
with $c_j$ defined as in~(\ref{bfc-def}),
%
%
\begin{equation} \label{sub-conv} |\lambda_n| n^{-(\tau-2)/(\tau
-1)}\bigl|\Ccal_{(j)}\bigr|\convp c_j.
\end{equation}
\end{Theorem}

Theorem~\ref{thm-sub-crit} is proved in Section~\ref{sec-sub-crit}.
Interestingly, the limit in~(\ref{sub-conv}) is \textit{deterministic}
[recall also~(\ref{bfc-def})].
The rough idea for this is as follows. As $\lambda=\lambda
_n\rightarrow-\infty$, the
random graph becomes more and more subcritical. Now, if we look at
$\Ccal(j)$, the cluster of vertex $j$, then we can view it
as the union of approximately $w_j$ (which is
roughly the degree of vertex $j$) almost \textit{independent} clusters.
These clusters
are close to total progenies\vspace*{1pt} of branching processes having mean offspring
$\nu_n(\lambda_n)\approx1+\lambda_n n^{-(\tau-3)/(\tau-1)}$. The
expected total progeny
of a branching process with mean offspring $\nu$ equals $1/(1-\nu)$.
As a result,
the expected cluster size of vertex $j$ is close to
%
%
\begin{equation} \frac{w_j}{1-\nu_n(\lambda_n)}\approx\frac
{w_j}{|\lambda_n|n^{-(\tau-3)/(\tau-1)}} =\frac{n^{(\tau-2)/(\tau
-1)}}{|\lambda_n|} c_j \bigl(1+o(1)\bigr).
\end{equation}
In our setting, $c_j=c_{ F}^{1/(\tau-1)}j^{-1/(\tau-1)}$, so that
$j\mapsto c_j$ is strictly
decreasing. Thus, we must also have that $|\Ccal(j)|=|\cluster_{
(j)}|$ \textit{whp}.
The proof of Theorem~\ref{thm-sub-crit}
makes this argument precise, by investigating the deviation
from a branching process, a technique that is also crucially used in
\cite{Hofs09a} to study tightness of the sequence of
random variables $|\Cmax|n^{-(\tau-2)/(\tau-1)}$.
A result similar to Theorem~\ref{thm-sub-crit} is proved for the
near-critical phase of
the configuration model in~\cite{HofJanLuc09}, but the proof we give
here is
entirely different.

We also obtain that the ordered cluster \textit{weights}
as defined in~(\ref{weight-cluster-def}) satisfy the same scaling results
as described above.
\begin{Theorem}[(Scaling limit of cluster weights)]
\label{thm-scal-lim-cluster-weights}
Theorems~\ref{thm-WC-3,4},~\ref{thm-mult-coal}\break and~\ref{thm-sub-crit}
also hold for the ordered cluster weights
$(n^{-(\tau-2)/(\tau-1)}\weight_{(i)})_{i\geq1}$,
with \textit{identical} scaling limits as in
Theorems~\ref{thm-WC-3,4},~\ref{thm-mult-coal} and~\ref{thm-sub-crit}.
\end{Theorem}

As explained in more detail in Section~\ref{sec-cluster-BP} below,
Theorem~\ref{thm-scal-lim-cluster-weights} can be heuristically
understood by noting that the average weight of a vertex in a cluster
is close to $\nu=1$, and therefore it contributes the same to the
weight of the cluster as it does to the cluster size. In fact, the
proof will
show that $n^{-(\tau-2)/(\tau-1)}\weight_{(i)}$ and $n^{-(\tau
-2)/(\tau-1)}|\cluster_{(i)}|$
converge to the \textit{same} limit.
The proof of Theorem~\ref{thm-scal-lim-cluster-weights}
shall be given simultaneously with the proofs of
Theorems~\ref{thm-WC-3,4},~\ref{thm-mult-coal} and~\ref{thm-sub-crit},
respectively, adapted so as to deal with cluster weights or cluster sizes.
Sometimes, it is more convenient to study cluster sizes
(e.g., since cluster explorations can more naturally
be formulated in terms of the number of vertices than their weight),
in some cases it is more convenient to work with cluster weights
(e.g., since the cluster\vadjust{\goodbreak} weights can be described in terms
of multiplicative coalescents, a fact that is crucial in the
proof of Theorem~\ref{thm-mult-coal}).

\subsection{Properties of large critical clusters}
\label{sec-furth-res}
We shall also derive some related interesting properties of the
limiting largest clusters. In the following theorem,
we consider the connectivity
structure of the high-weight vertices:

\begin{Theorem}[(Connectivity of high-weight vertices)]
\label{thm-prop-3,4}
Under the assumptions in Theorem~\ref{thm-WC-3,4},
for every $i,j\geq1$ \textit{fixed},
%
%
\begin{equation} \label{lim-qij}
\lim_{n\rightarrow\infty} \prob
\bigl(j\in\cluster(i)\bigr)=q_{ij}(\lambda)\in(0,1)
\end{equation}
and
%
%
\begin{equation} \label{lim-qi} \lim_{n\rightarrow\infty} \prob
(i\in\Cmax)=q_{i}(\lambda)\in(0,1).
\end{equation}
\end{Theorem}

Theorem~\ref{thm-prop-3,4} states that the high-weight vertices
play an essential role in the critical regime. Indeed, we shall see that
in the subcritical regime, with high probability, $\Cmax=\cluster(1)$,
so that $\prob(1\in\Cmax)=1-o(1)$, while $\prob(i\in\Cmax)=o(1)$
for $i>1$.
In the supercritical regime, instead, $\prob(i\in\Cmax)=1-o(1)$
for every $i\geq1$ fixed. Thus, the critical regime is precisely
the regime where the high-weight vertices start to form connections.
Informally, this can be phrased as ``power to the wealthy.''
Theorem~\ref{thm-prop-3,4} should be contrasted with the situation
when $\expec[W^3]<\infty$ studied in~\cite{BhaHofLee09a,Turo09},
where the probability that any specific vertex is an element of $\Cmax$
is negligible, and, instead, the largest cluster is born
out of many trials each having a small probability.
This can be informally phrased as ``power to the masses.''

The following theorem, which is a crucial ingredient
in the proof of Theorem~\ref{thm-WC-3,4},
essentially says that, for each fixed $\lambda$, the
maximal size components are those attached to the largest weight
vertices:
\begin{Theorem}[(Large clusters contain a high-weight vertex)]
\label{prop-max-clusters-high-weight}
Assume that the conditions in Theorem~\ref{thm-WC-3,4} hold. Then:

\begin{longlist}[(a)]
\item[(a)]
for any $\vep\in(0,1)$, there exists a $K=K(\vep)\geq1$, such that,
for all $n$,
%
%
\begin{equation} \prob\Bigl(\max_{i\geq K} \bigl|\Ccal_{\leq
}(i)\bigr|\geq\vep n^{(\tau-2)/(\tau-1)}\Bigr) \leq\vep;
\end{equation}
\item[(b)] for any $m\geq1$,
%
%
\begin{eqnarray}
&&\lim_{\vep\downarrow0} \liminf_{n\rightarrow
\infty}\prob\bigl(\bigl(\bigl|\Ccal_{\leq}(i)\bigr|\bigr)_{i\in[n]} \mbox{
contains $m$ components of size} \nonumber\\[-8pt]\\[-8pt]
&&\hspace*{184pt} \geq\vep n^{(\tau
-2)/(\tau-1)}\bigr) =1.\nonumber
\end{eqnarray}
\end{longlist}
\end{Theorem}

\subsection{Overview of the proofs}
\label{sec-overview}

In this section, we give an overview of the proofs of our main results.
We start by explaining the proof of Theorem~\ref{thm-WC-3,4}, along
the way also explaining the key ideas behind Theorems~\ref{thm-prop-3,4}
and~\ref{prop-max-clusters-high-weight}. After this, we shall discuss the
proofs of Theorem~\ref{thm-mult-coal} and~\ref{thm-sub-crit}.\vadjust{\goodbreak}

We note that, since $u\mapsto[1-F]^{-1}(u)$ is nonincreasing,
$\bfw$ is ordered in size, that is, $w_1\geq w_2\geq w_3\geq\cdots.$
We start by \textit{exploring} the clusters from the largest weight
vertices onward.
Here, by a cluster exploration, we mean the recursive investigation of
the neighbors of the vertices already found to be in the cluster. This
cluster exploration shall be described in detail in Section~\ref
{sec-cluster-BP}.
The rough idea is as follows. We start with a vertex $i$, and
wish to find all the vertices that are in its cluster.
For this,
we sequentially take each vertex in the cluster being currently
explored and find its direct neighbors, that have not yet been found
by the exploration process.
Call a vertex \textit{active}
when it
is found to be in the cluster, but has not yet been explored.
A vertex is called \textit{explored} when its neighbors have been
investigated and \textit{neutral} when it has not yet
appeared in the exploration process. Then, in the exploration process
at time $t$, we take a vertex, turn it from active to explored,
and explore it, that is, see which neutral neighbors it has. Turn the
status of its neutral neighbors to active. Let $Z_l$ denote the number
of active
neighbors after the exploration of the $l$th active vertex.
When $Z_l=0$ for the first time, then there are no more active vertices,
so all elements of the cluster have been found. (The description
in Section~\ref{sec-cluster-BP} is slightly different than the one
described here,
as it studies \textit{potential} elements of the cluster instead.)

We note that the high-weight vertices have weights of the order
$w_j\sim(c_{ F}n/j)^{1/(\tau-1)}$, so, when we start
with a high-weight vertex, initially, the number of active vertices shall
be of the order $n^{1/(\tau-1)}$. When our exploration process
hits \textit{another} high-weight vertex, then the number of active vertices
gets a large push of the order $n^{1/(\tau-1)}$ upward. It is these
upward pushes that change the number of active vertices in a substantial
way, and, therefore, the high-weight vertices play a crucial role
in the critical behavior of our random graph. In turn, this suggests
that the largest clusters contain at least one high-weight vertex,
as indicated by Theorem~\ref{prop-max-clusters-high-weight}.
Due to the critical nature of our random graph, it turns
out that the \textit{average} number and weight of active vertices
being added
in each exploration is close to one, so that, due to the removal
of the vertex which is being explored, the exploration process
has increments that have a mean close to zero.

In Section~\ref{sec-WC-C1}, we start by identifying the scaling
limit of $n^{-(\tau-2)/(\tau-1)}\*|\cluster_{\leq}(1)|=n^{-(\tau
-2)/(\tau-1)}|\cluster(1)|$.
The weak limit of $n^{-(\tau-2)/(\tau-1)}|\cluster(1)|$ is given in terms
of the hitting time of zero of an exploration process exploring the
cluster of vertex~1
(the vertex with the highest weight).
See Theorems~\ref{thm-WC-C1} and~\ref{thm-weak-lim-S}.
The scaling limit of the exploration process of a cluster exists (see
Theorem~\ref{thm-weak-lim-S}), and can be viewed as a ``thinned''
L\'evy process. Therefore, the convergence
in distribution of $n^{-(\tau-2)/(\tau-1)}|\cluster(1)|$ in
Theorem~\ref{thm-WC-C1} is equivalent to the convergence of the
first hitting time
of zero of the exploration process to the one of this thinned L\'evy process.
In proving this, we perform a careful analysis of hitting times of
a spectrally
positive L\'evy process that stochastically dominates the thinned L\'
evy process.

Following the proof of convergence of $n^{-(\tau-2)/(\tau
-1)}|\cluster(1)|$ in Theorem~\ref{thm-WC-C1}, we shall prove the
convergence in distribution of
$(n^{-(\tau-2)/(\tau-1)}|\Ccal_{\leq}(i)|)_{i\in[n]}$ in Theorem
\ref{thm-WC-C[K]}. This proof makes crucial use of the estimates in the
proof of Theorem~\ref{thm-WC-C1}, and allows us to extend the result in
Theorem~\ref{thm-WC-C1} to the (joint) convergence of several rescaled
clusters by an inductive argument. 
The largest
$m$ clusters are given by the largest $m$ elements of the vector
$(|\Ccal_{\leq}(i)|)_{i\in[n]}$, so that this completes the proof of
Theorem~\ref{thm-WC-3,4}. The conclusion of this argument shall be
carried out in Section~\ref{sec-pf-prop-max-clusters-high-weight}.

In Section~\ref{sec-sub-crit}, we prove Theorem~\ref{thm-sub-crit} by
a second moment argument, using the fact that
the \textit{subcritical} phase of our random graph is closely
related to (and even stochastically dominated by) a branching process.
In Section~\ref{secmult-coal}, we use the results
proved in Section~\ref{sec-sub-crit}, jointly with the results in
\cite{AldLim98}, to prove Theorem~\ref{thm-mult-coal}.
We now discuss in a bit more detail how one can understand the
appearance of multiplicative coalescents in the random graphs we study here.

We make crucial use of~\cite{AldLim98}, Proposition 7, whose
application we now explain. Fix a sequence $\lambda_n\rightarrow
-\infty$.
For each fixed $t$, consider the construction of the inhomogeneous random
graph as in~(\ref{pij-NR}) but with the weight sequence
$\bar{\bfw}(t)=(\bar{w}_j(t))_{1\leq j\leq n}$ given by
%
%
\begin{equation} \label{eqnweight-wt} \bar w_j(t)=w_j\bigl(1+(t+\lambda
_n)\ell_n n^{-2(\tau-2)/(\tau-1)}\bigr).
\end{equation}
Let
%
%
\begin{equation} {\mathbf X}^{(n)}(t)=\bigl(n^{-(\tau
-2)/(\tau-1)}\weight_{(i)}(t)\bigr)_{i\geq1}
\end{equation}
denote the ordered version of cluster weights when the vertex weights
are given by
$\bar{\bfw}(t)$.

Note that the above process, when taking $t=-\lambda_n+\lambda/\expec
[W]$, is closely related to
the ordered cluster weights of our random graph with weights
$w_j(\lambda)=w_j(1+\lambda n^{-(\tau-3)/(\tau-1)})$,
since $\ell_n=\expec[W] n (1+o(1))$. We then note that ${\mathbf
X}^{(n)}$ can be constructed
so that, viewed as a function in $t$, it is a multiplicative coalescent.
\begin{Lemma}[(Discrete multiplicative coalescent)]
\label{lem-mult-coal}
We can construct the process ${\mathbf X}^{(n)}=({\mathbf
X}^{(n)}(t))_{t\geq0}$ such that, for each
fixed $t$, ${\mathbf X}^{(n)}(t)$ has the distribution
of the ordered rescaled weighted component
sizes of the random graph with weight sequence given by~(\ref{eqnweight-wt})
and such that, for each fixed $n$, the process viewed as a process in
$t$ is a
multiplicative coalescent. The initial state denoted by $\bfitx^{
(n)}(0)$ has the
same distribution as the ordered cluster weights of a random graph with edge
probabilities as in~(\ref{pij-NR}) and weight sequence
%
%
\begin{equation} \bar w_j(0)=w_j\bigl(1+\lambda_n\ell_n n^{-2(\tau
-2)/(\tau-1)}\bigr).
\end{equation}
\end{Lemma}
\begin{pf}
For each unordered pair $(i,j)$, let $\xi_{ij}$ be an exponential
random variable with rate $w_iw_j/\ell_n$, where $(\xi_{ij})_{(i,j)}$
are independent. For fixed $t$, define the graph $\bar{\GG}_n^t$ to
consist of all edges $(i,j)$ for which
%
%
\begin{equation} \xi_{ij} \leq1+\frac{(\lambda_n+t)\ell
_n}{n^{2(\tau-2)/(\tau-1)}}.
\end{equation}
Then, by construction, for all $t\geq0$, the rescaled weighted
component sizes of $\bar{\GG}_n^t(\bfwit)$ have the same
distribution as ${\mathbf X}^{(n)}(t)$. Further, for any time
$t$ we note that two \textit{distinct} clusters $\cluster_1$ and
$\cluster_2$ having weights $\weight_{(i)}(t)$ and $\weight
_{(j)}(t)$, respectively, coalesce at rate
%
%
\begin{eqnarray}
&&
\ell_n n^{-2(\tau-2)/(\tau-1)} \sum_{s_1\in
\cluster_1, s_2\in\cluster_2} \frac{w_{s_1}w_{s_2}}{\ell_n} \nonumber\\[-8pt]\\[-8pt]
&&\qquad=
\bigl(n^{-(\tau-2)/(\tau-1)}\weight_{(i)}(t)\bigr)
\bigl(n^{-(\tau-2)/(\tau-1)} \weight_{(j)}(t)\bigr)\nonumber
\end{eqnarray}
as required.
\end{pf}

In effect, Theorems~\ref{thm-WC-3,4} and~\ref{thm-mult-coal} give us
two \textit{distinct} proofs of the statement that the ordered cluster
weights converge, and we now discuss the advantages of these two different
proofs. Theorem~\ref{thm-WC-3,4} proves that for any \textit{fixed}
$\lambda$,
$(n^{-(\tau-2)/(\tau-1)}\weight_{(i)})_{i\geq1}$ converges
in distribution. Further, by the fact that this vector is obtained
by sequentially investigating the clusters of the high-weight vertices,
it allows us to prove properties about the high-weight vertices that are
part of the largest clusters, as in Theorems
\ref{thm-prop-3,4} and~\ref{prop-max-clusters-high-weight}. Finally,
it allows us to show that the ordered cluster \textit{sizes} have the same
scaling limit as the ordered cluster \textit{weights} (see Theorem \ref
{thm-scal-lim-cluster-weights}), a feature that is also crucial in the
proofs of Theorems~\ref{thm-mult-coal} and~\ref{thm-sub-crit}.

Theorem~\ref{thm-mult-coal}, instead, shows that the
\textit{process} of the ordered cluster sizes or weights
converges in distribution. This means that there exists a
\textit{stochastic process} that describes the joint
convergence of the ordered cluster sizes or weights
for different values of $\lambda$ simultaneously.
Due to the fact that the proof of Theorem~\ref{thm-mult-coal}
relies on~\cite{AldLim98}, Proposition 7, however, we obtain
less information about the vertices that are part of
the large critical clusters. The combination of the two proofs
provides us with a detailed and full understanding of
the scaling limit of the ordered cluster sizes or weights.


\subsection{Discussion}
\label{sec-disc3,4}

\subsubsection*{Comparison to the case of weights with finite third moments}
In~\cite{Aldo97,BhaHofLee09a,Turo09}, the scaling limit was considered
when $\expec[W^3]<\infty$. In this case, the scaling limit
turns out to be (a trivial rescaling of) the scaling limit for the
Erd\H{o}s--R\'enyi random graph as found by Aldous in~\cite{Aldo97}.
Thus, the
setting for $\tau\in(3,4)$ is fundamentally different. When $\expec
[W^3]<\infty$,
the probability that $1\in\Cmax$ is negligible, while in our setting this
is not true, as shown in Theorem~\ref{thm-prop-3,4}.

\subsubsection*{Other weights} Our proof reveals that the precise limits
of $w_in^{-1/(\tau-1)}$, for fixed $i\geq1$, arise in the scaling limit.
We make crucial use of the fact that, by~(\ref{[1-F]inv-bd})
$c_i=\lim_{n\rightarrow\infty}w_in^{-1/(\tau-1)}=(c_{
F}/i)^{1/(\tau-1)}$.
However, we believe that our results can be appropriately adapted to
the situation
that $\lim_{n\rightarrow\infty}w_in^{-1/(\tau-1)}$
exists for every $i\geq1$ and is \textit{asymptotically} equal to
$ai^{-1/(\tau-1)}$
for some $a>0$. This suggests that, by varying the
precise values of high weights, there are many possible scaling limits.
It would be of interest to investigate this further.

Also, we restrict
to tail distributions $1-F(x)$ that are, for large $x\geq0$,
asymptotic to an inverse power of $x$; see~(\ref{F-bound-tau3,4b}).
It would be of interest to investigate the
scaling behavior when~(\ref{F-bound-tau3,4b}) is replaced with the assumption
that $1-F(x)$ is regularly varying with exponent $1-\tau$, that is,
$[1-F](x)=x^{-(\tau-1)}\ell(x)$ for some $x\mapsto\ell(x)$ which is slowly
varying at $\infty$. In this case,
we believe that the asymptotic sizes of the largest critical
clusters are given by $\ell^*(n)n^{(\tau-2)/(\tau-1)}$
for some suitable slowly varying function $n\mapsto\ell^*(n)$
that can be described in terms of $x\mapsto\ell(x)$.
For more details, see~\cite{Hofs09a}, Section~1.3,
where also the critical cases $\tau=3$ and $\tau=4$ are
discussed.

\subsubsection*{I.i.d. weights}
In our analysis, we make crucial use of the choice for $w_i$ in~(\ref
{choicewi}).
In the literature, also the setting where $(W_i)_{i\in[n]}$ are
independent and
identically distributed (i.i.d.) random
variables with distribution function $F$ has been considered. We expect
the behavior in
this model to be \textit{different}. Indeed, let $w_{i}=W_{(i)}$, where
$W_{(i)}$ are the order statistics of the i.i.d. sequence
$(W_i)_{i\in[n]}$.
It is well known that
%
%
\begin{equation} n^{-1/(\tau-1)} W_{(i)}\convd\xi_i\equiv
(E_1'+\cdots+E_i')^{-1/(\tau-1)},
\end{equation}
where $(E_i')_{i=1}^{\infty}$ are i.i.d. exponential random variables
with mean 1.
In particular, when $\tau\in(3,4)$, $\expec[\xi_1^a]<\infty$
whenever $a<\tau-1$.
The extra randomness of the order statistics has an effect on the scaling
limit, which is thus \textit{different}. In most cases, the two
settings have
the \textit{same} behavior (see, e.g.,~\cite{BhaHofLee09a}, where
this is
shown to hold for weights for which $\expec[W^3]<\infty$, where $W$
has distribution
function~$F$). See~\cite{Jose10} for the identification of the scaling
limit of
the largest cluster sizes in the critical configuration model with
i.i.d. degrees, which
is markedly \textit{different} from ours.
We believe that the same applies to the Norros--Reittu model with i.i.d.
weights.

\subsubsection*{High-weight vertices} The fact that the vertex $i$ is in
the largest connected component with nonvanishing probability as
$n\rightarrow\infty$
(see Theorem~\ref{thm-prop-3,4}) is remarkable. In our setting, a
\textit{uniformly}
chosen vertex in $[n]$ is an element of $\Cmax$ with negligible probability.
The point\vspace*{1pt} is that vertex $i$ has weight $w_i$, which, for $i$ fixed,
is close to $(c_{ F}/i)^{1/(\tau-1)}n^{1/(\tau-1)}$,
while a uniformly chosen vertex has a bounded weight.
Thus, Theorem~\ref{thm-prop-3,4}
can be interpreted as saying that the highest-weight vertices characterize
the largest components. In the subcritical case (see, e.g.,
the results by Janson in~\cite{Jans08b} or Theorem~\ref{thm-sub-crit}),
the largest connected component is the one of
the vertex with the highest weight, and the critical situation arises when
the highest-weight vertices start connecting to each other.


\subsubsection*{Connection to the multiplicative coalescent}
%
%
The mental picture associated with the entrance
boundary of the coalescent here seems to be different from~\cite{AldLim98},
where in spirit many of the component sizes are of order $n^{2/3}$.
Here the
entrance boundary describes the sizes of the maximal components
rescaled by $n^{-(\tau-2)/(\tau-1)}$
in the $\lambda\to-\infty$ regime, whilst in~\cite{AldLim98} they
arise as
limits of random graphs similar to critical Erd\H{o}s--R\'enyi random graphs,
where, in addition to the random edges, there are initially a number of large
``planted'' components of sizes $\lfloor c_i n^{1/3}\rfloor$;
see~\cite{AldLim98}, Section 1.3. However, the results of~\cite{AldLim98}
are crucial in identifying the distribution of the limiting component sizes
for fixed $\lambda$. It would be interesting to see if the stochastic calculus
techniques developed in~\cite{AldLim98} can be further modified to
give useful
information about the surplus of edges in the maximal components [the
surplus of a
component $\CC$ with $E(\CC)$ edges and $V(\CC)$ vertices is equal to
$E(\CC) - (V(\CC)-1)$ and denotes the minimal number of edges that
must be removed
from the component to make it a tree].

\section{The scaling limit of the cluster of vertex 1}
\label{sec-WC-C1}

In this section, we identify the scaling limit of $|\cluster(1)|$.
We note from~(\ref{choicewi}) that the weight of vertex 1 is \textit{maximal},
that is, $w_1\geq w_2\geq\cdots\geq w_n$. When $\tau>4$,
the probability that vertex 1 belongs to $\Cmax$ is negligible.
When $\tau\in(3,4)$, instead, we shall see that
vertex 1 is in $\Cmax$ with \textit{positive} probability,
so that it is quite reasonable to start exploring the cluster of vertex
1 first,
since $|\cluster(1)|$ stochastically dominates $|\cluster(j)|$ for
all $j\in[n]$.
Theorem~\ref{thm-WC-C1} below states that $|\cluster(1)|$ is of order
$n^{(\tau-2)/(\tau-1)}$. By~\cite{Hofs09a}, Theorem~1.2, the same
is valid for $|\Cmax|$, which confirms the above heuristic.
\begin{Theorem}[{[Weak convergence of the cluster of vertex 1 for $\tau
\in(3,4)$]}]
\label{thm-WC-C1}
Fix the Norros--Reittu random graph with weights
$\bfw(\lambda)$ defined in~(\ref{w-lambda-def}).
Assume that $\nu=1$ and that~(\ref{F-bound-tau3,4b}) holds. Then,
for all $\lambda\in\Rbold$,
%
%
\begin{equation}
n^{-(\tau-2)/(\tau-1)}|\cluster(1)|\convd H_1(0)
\end{equation}
for some nondegenerate limit $H_1(0)$.
\end{Theorem}

Theorem~\ref{thm-WC-C1} is proved in Section~\ref{sec-pf-thm-WC-C1}.
We now start by discussing \textit{cluster explorations} and their
relation to branching processes, which play an essential role in our proofs.

\subsection{Cluster explorations and their relation to branching processes}
\label{sec-cluster-BP}

We fix the weight sequence to be $\bfw(\lambda)$ defined in (\ref
{w-lambda-def}),
and we shall denote the weight of vertex $i$ [or the $i$th coordinate of
$\bfw(\lambda)$] by $w_i(\lambda)$.\vadjust{\goodbreak}

In order to prove Theorem~\ref{thm-WC-C1}, we make heavy use of the
\textit{cluster
exploration}, which is described in detail in~\cite{NorRei06} and
\cite{Hofs09a}.
The model in~\cite{NorRei06} is a \textit{random multigraph}, that is,
a random graph potentially having self-loops and multiple edges.
Indeed, for each $i,j\in[n]$, we let the number of edges between
vertex $i$ and $j$ be $\Poi(w_i(\lambda)w_j(\lambda)/\ell_n(\lambda
))$, where, for $\mu\geq0$,
we let $\Poi(\mu)$ denote a Poisson random variable with mean $\mu$,
and we
define
%
%
\begin{equation} \label{ell-n-lambda-def} \ell_n(\lambda)=\sum
_{i\in[n]} w_i(\lambda)=\ell_n \bigl(1+\lambda n^{-(\tau-3)/(\tau
-1)}\bigr).
\end{equation}
The number of edges between different pairs of vertices are \textit
{independent}.
To retrieve our random graph model, we merge multiple edges and erase
self-loops. Then, the probability that an edge exists between two vertices
$i,j\in[n]$ is equal to
%
%
\begin{equation}
p_{ij}=\prob\bigl(\Poi\bigl(w_i(\lambda)w_j(\lambda
)/\ell_n(\lambda)\bigr)\geq1\bigr)=1-{\mathrm e}^{-w_i(\lambda
)w_j(\lambda)/\ell_n(\lambda)}
\end{equation}
as required. Further, the number of potential edges from a vertex $i$
has a Poisson
distribution with mean $w_i(\lambda)$. We shall work with the above
Poisson random graph instead,
and we shall refer to the Poisson random variable $\Poi(w_i(\lambda
))$ as the number
of \textit{potential neighbors} of vertex $i$. When we find what the
vertices are that correspond
to these $\Poi(w_i(\lambda))$ potential neighbors, that is, when we
determine their \textit{marks},
then we can see how many real neighbors there are. Here by a ``mark'' we
mean a random
variable $M$ with distribution
%
%
\begin{equation} \label{dist-mark} \prob(M=m)= w_m(\lambda)/\ell
_n(\lambda)=w_m/\ell_n,\qquad 1 \leq m \leq n.
\end{equation}
The variable $M$ corresponds to the actual \textit{vertex label}
associated to the potential vertex. A potential vertex
arising in our exploration process is an actual vertex when
its mark has not arisen in the exploration up to that point.
We now describe this cluster exploration in detail.

We denote by $(Z_l)_{l\geq0}$ the exploration process in the
breadth-first search,
where $Z_0=1$ and where $Z_1$ denotes the number of potential neighbors
of the initial vertex (which is in the case of Theorem~\ref{thm-WC-C1}
equal to vertex 1, and which we shall often take to be vertex $i$).
The variable $Z_l$ has the interpretation of the
number of potential neighbors of the first $l$ explored potential vertices
in the cluster whose neighbors have not yet been explored. As a result,
we explore by taking one vertex of the ``stack'' of size $Z_l$, drawing
its mark and checking
whether it is a real vertex, followed by drawing its number of
potential neighbors. Thus, we set $Z_0=1, Z_1=\Poi(w_i(\lambda))$,
and note that, for $l\geq2$,
$Z_l$ satisfies the recursion relation
%
%
\begin{equation} \label{Zl-recur} Z_l=Z_{l-1}+X_l-1,
\end{equation}
where $X_l$ denotes the number of potential neighbors of the
$l$th potential vertex that is explored.
More precisely, when we explore a potential vertex, we start by drawing
its mark in an i.i.d. way with distribution~(\ref{dist-mark}). When
we have already
explored a vertex with the same mark as the one drawn, we turn the
status of
the vertex to be explored to inactive, the potential vertex does not
become a real vertex, and
proceed with the next potential vertex. When, instead, it receives
a mark which we have not yet seen, then the potential vertex becomes a real
vertex, its mark $M_l\in[n]$ indicating which vertex in $[n]$
the $l$th explored vertex corresponds to, so that $M_l\in\cluster
(i)$. We then
draw $X_l=\Poi(w_{M_l})$, and $X_l$ denotes the number of potential
vertices incident to the real vertex $M_l$. Again, upon exploration,
these potential vertices might become real vertices, and this occurs precisely
when their mark corresponds to a vertex in $[n]$ that has not appeared
in the
cluster exploration so far. We call the above procedure of drawing a mark
for a potential vertex to investigate whether it corresponds to a real vertex
a \textit{vertex check}.

In~\cite{NorRei06}, Proposition 3.1 (see also
\cite{Hofs09a}, Section 3.2, in particular, Proposition~3.4),
the cluster exploration was described
in terms of a thinned marked mixed Poisson branching process.
This description implies that the distribution of $X_l$ (for $2\leq
l\leq n$)
is equal to $\Poi(w_{M_l}(\lambda))J_l$, where (a) the marks
$(M_l)_{l=2}^{\infty}$
are i.i.d. random variables with distribution~(\ref{dist-mark}); and
(b) $J_l=\indic{M_l\notin\{i\}\cup\{M_2, \ldots, M_{l-1}\}}$ is
the indicator
that the mark $M_l$ has not been found before and is not 1. Here, the mark
$M_l$ is the label of the potential element of the cluster that we are
exploring, and, clearly, if a vertex has already been observed to be part
of $\cluster(i)$, and its\vspace*{1pt} neighbors have been explored, then we should
not do so again. We sometimes write $J_j^{(i)}$, $X_l^{(i)}$
and $Z_l^{(i)}$ to explicitly indicate the vertex whose cluster
we are exploring, and omit the superscript when no confusion can arise.

We conclude that we arrive at, for $l\geq2$,
%
%
\begin{eqnarray} \label{Zl-recur-MP}
Z_l=Z_{l-1}+X_l-1\hspace*{120pt}\nonumber\\[-8pt]\\[-8pt]
&&\eqntext{\mbox
{where } X_l=\Poi(w_{M_l}(\lambda))J_l
\quad\mbox{and}\quad J_l=\indic{M_l\notin\{i\}\cup\{M_2, \ldots, M_{l-1}\}
}.}
\end{eqnarray}
Then, the number of vertex checks that have been performed
when exploring the cluster of vertex $i$ equals $V(i)$,
which is given by
%
%
\begin{equation} \label{V1-def} V(i)=\inf\{l\dvtx Z_l=0\}
\end{equation}
since the first time at which there are no more potential vertices to
be checked,
all vertices in the cluster have been checked.

Further, the number of \textit{real vertices} found to be part of
$\cluster(i)$
after $l$ vertex checks equals
%
%
\begin{equation} |\cluster(i;l)|=1+\sum_{j=2}^l J_j,
\end{equation}
that is, all the potential vertices, except for those that have a mark
that has appeared previously. Therefore, we conclude that
%
%
\begin{equation} \label{cluster1-comp} |\cluster(i)|=1+\sum
_{j=2}^{V(i)} J_j =V(i)-\sum_{j=2}^{V(i)}(1-J_j).
\end{equation}
It turns out that the second contribution is an error term
(see Lemma~\ref{lem-multiple-hits} below), so that
the cluster size of 1 asymptotically corresponds to the first hitting
time of $0$ of
$l\mapsto Z_l$. We prove Theorem~\ref{thm-WC-C1} by applying the
above to $i=1$.

Throughout the paper, we abbreviate
%
%
\begin{equation} \label{defs-al-rho-eta}\qquad
\alpha=1/(\tau-1),\qquad
\rho=(\tau-2)/(\tau-1), \qquad\eta=(\tau-3)/(\tau-1).
\end{equation}

\subsection{Branching process computations}
\label{sec-BPs}

In this section, we discuss some useful facts about branching processes.
Note that if, in the recursion arising in the exploration of the cluster
in~(\ref{Zl-recur-MP}), we ignore the $J_l$'s (i.e., we
ignore the effect of marks that have already been used),
then we arrive at the recursion
%
%
\begin{equation} \label{ZBP-def} \ZBP_l=\ZBP_{l-1}+\XBP_l-1,
\end{equation}
where now
%
%
\begin{equation} \label{XBP-def} \XBP_l=\Poi(w_{M_l}(\lambda)),
\end{equation}
and where $(\Poi(w_{M_l}(\lambda)))_{l\geq2}$ are i.i.d. random variables,
while $M_1=i$. This recursion is the random walk description in the
exploration of the total progeny of a branching process. Indeed, let
%
%
\begin{equation} \label{T-def} T(i)=\inf\bigl\{l\dvtx\ZBP_l=0\bigr\}
\end{equation}
be the first hitting time of 0 of the process $(\ZBP_l)_{l\geq0}$.
Then, by the random walk description of a branching process
(see, e.g.,~\cite{Hofs08}, Section 3.3), $T(i)$ has the same
distribution as
the total progeny of a branching process in which the root has
offspring distribution $\Poi(w_i(\lambda))$, while the offspring of all
other individuals is i.i.d. with mixed Poisson offspring distribution
$\Poi(w_{M}(\lambda))$, where $M$ is the mark distribution in (\ref
{dist-mark}).
In the setting in Section~\ref{sec-cluster-BP}, we have $i=M_1=1$, so
that we
start from the root having mark 1, but in this section,
we shall generalize as well to $M_1=i$, where $i\in[n]$.
Further, we shall also denote the total progeny of the branching process
with offspring distribution $\Poi(w_{M}(\lambda))$ by $T$. In this
section, we
investigate properties of such branching processes.

The connection to branching processes [in particular, the
\textit{stochastic domination} of the cluster sizes by
branching processes due to~(\ref{Zl-recur-MP})]
plays a crucial role in~\cite{Hofs09a}, where
this comparison was used in order to prove
that $n^{-\rho}|\Cmax|$ and $n^{\rho}/|\Cmax|$
are tight sequences of random variables. There, only \textit{bounds}
on the maximal cluster size were shown, while, in this paper, we identify
the \textit{scaling limit} of all large clusters.

The difference between the branching process recursion relation
in~(\ref{ZBP-def}) and~(\ref{XBP-def}), and the corresponding one for
the cluster exploration in~(\ref{Zl-recur-MP}) resides in the random
variables $(J_l)_{l\geq1}$. Indeed, when $J_l=0$, then
$X_l=\Poi(w_{M_l}(\lambda))J_l=0$, while
$\XBP_l=\Poi(w_{M_l}(\lambda))$ is unaffected.
Therefore, we can think of this procedure as a \textit{thinning}
of our branching process. Indeed, when the mark of the $l$th
potential vertex has been seen before, then, in the cluster
exploration, we remove this vertex and \textit{all of its offspring}.
Thus, the recursions in~(\ref{ZBP-def}),~(\ref{XBP-def}) and (\ref
{Zl-recur-MP})
give us a \textit{simultaneous coupling} of the cluster exploration
process and
the branching process such that any deviation between the two arises
from the thinning
of the potential vertices and the subsequent removal of the branching
process tree
that is attached to the thinned potential vertices. This description
shall prove
to be crucial in the comparison of cluster sizes and branching process
total progenies used in the proofs of Theorems~\ref{thm-mult-coal} and
\ref{thm-sub-crit}.

We continue to investigate the critical behavior of the
branching processes at hand.
We denote
%
%
\begin{equation} \label{nun-def} \nu_n(\lambda)=\frac{1}{\ell
_n}\sum_{j\in[n]} w_jw_j(\lambda),
\end{equation}
and we write $\nu_n=\nu_n(0)$. Then, we note that
%
%
\begin{equation} \label{nun-rep} \nu_n(\lambda)=\expec[\Poi
(w_{M}(\lambda))],
\end{equation}
so that $\nu_n(\lambda)$ is the \textit{mean offspring} of the
branching process, and
$\nu_n(\lambda)\rightarrow1$ corresponds
to our branching process being \textit{critical}.
Further,
%
%
\begin{eqnarray}
&&
\expec\bigl[\Poi(w_{M}(\lambda))\bigl(\Poi(w_{M}(\lambda
))-1\bigr)\bigr]\nonumber\\[-8pt]\\[-8pt]
&&\qquad =\expec[w_{M}^2(\lambda)]=\frac{1}{\ell_n}\sum_{j\in
[n]}w_j w_j(\lambda)^2\rightarrow\infty,\nonumber
\end{eqnarray}
so that our branching process has asymptotically \textit{infinite variance}
in the setting in~(\ref{F-bound-tau3,4b}).
We now give detailed asymptotics for the mean $\nu_n=\nu_n(0)$
of the above branching process. From this asymptotics, we can easily
deduce the asymptotics of $\nu_n(\lambda)=\nu_n (1+\lambda n^{-\eta})$
[recall~(\ref{defs-al-rho-eta}),~(\ref{nun-def}) and~(\ref{w-lambda-def})].
%
%
%
\begin{Lemma}[(Sharp asymptotics of $\nu_n$)]
\label{lem-nun-asym}
Let the distribution function $F$ satisfy~(\ref{F-bound-tau3,4b}),
and let $\nu_n=\nu_n(0)$ be given by~(\ref{nun-def}) and $\nu$ by
(\ref{nu-def}). Then,
with $\eta$ given in~(\ref{defs-al-rho-eta}),
%
%
\begin{equation} \label{nun-asym} \nu_n= \nu+\zeta n^{-\eta
}+o(n^{-\eta}),
\end{equation}
where
%
%
\begin{equation} \label{zeta-form} \zeta=-\frac{c_{ F}^{2/(\tau
-1)}}{\expec[W]} \sum_{i=1}^{\infty} \biggl[\int_{i-1}^{i}
u^{-2\alpha}\,du-i^{-2\alpha}\biggr]\in(-\infty,0).
\end{equation}
\end{Lemma}
\begin{pf}
By~\cite{Hofs09a}, Corollary 3.2, $\ell_n= \sum_{i\in[n]} w_i = n
\expec[W]+O(n^{\alpha})$,
where it is also proved that $\nu_n-\nu=O(n^{-\eta})$. The sharper
asymptotics for $\nu_n$ in~(\ref{nun-asym}) is obtained by a more careful analysis of the arising sum.
We note that, by the remark below~(\ref{invverd}),
%
%
\begin{equation} \nu=\frac{\int_0^1 [1-F]^{-1}(u)^2\,du}{\int_0^1
[1-F]^{-1}(u)\,du}.
\end{equation}
By the asymptotics of $\ell_n$ above, we have that
%
%
\begin{equation} \nu_n=\frac{\sum_{i\in[n]} w_i^2}{n\expec
[W]}+o(n^{-\eta}).
\end{equation}
We shall make use of the fact that, when $f$ is nonincreasing,
%
%
\begin{equation} f(i)\leq\int_{i-1}^i f(u)\,du \leq f(i-1).
\end{equation}
Applying this to $f(u)=[1-F]^{-1}(u)^2$, which is nonincreasing, we obtain
in particular that, for any $K\geq1$,
%
%
\begin{eqnarray}
&&\int_{K/n}^1 [1-F]^{-1}(u)^2\,du - \frac{1}{n}
w_{K/n}^2\nonumber\\[-8pt]\\[-8pt]
&&\qquad\leq\frac{1}{n}\sum_{i=K+1}^n w_i^2\leq\int_{K/n}^1
[1-F]^{-1}(u)^2\,du.\nonumber
\end{eqnarray}
Now,
%
%
\begin{equation} \frac{1}{n} w_{K/n}^2= \frac{c_{ F}}{n}
(n/K)^{2\alpha}\bigl(1+o(1)\bigr) = \Theta(K^{-2\alpha}n^{-\eta}).
\end{equation}
Thus we conclude that
%
%
\begin{eqnarray}
\nu-\nu_n&=&\frac{1}{\expec[W] n} \sum_{i=1}^K\int
_{(i-1)/n}^{i/n} [1-F]^{-1}(u)^2\,du -\frac{1}{\expec[W] n} \sum
_{i=1}^K w_i^2 \nonumber\\[-8pt]\\[-8pt]
&&{}+\Theta(K^{-2\alpha}n^{-\eta})+o(n^{-\eta}).\nonumber
\end{eqnarray}
Next, by~(\ref{F-bound-tau3,4b}), for every $K\geq1$ fixed,
%
%
\begin{equation}
\frac{1}{n}\sum_{i=1}^K w_i^2=n^{-\eta} \sum
_{i=1}^K (c_{ F}/i)^{2\alpha}+o(n^{-\eta})
\end{equation}
and
\begin{equation}
\frac{1}{n} \sum_{i=1}^K\int_{(i-1)/n}^{i/n} [1-F]^{-1}(u)^2\,du
= n^{-\eta
} \sum_{i=1}^K \int_{i-1}^{i} (c_{ F}/u)^{2\alpha
}\,du+o(n^{-\eta}).\hspace*{-25pt}
\end{equation}
Combining these two estimates yields
%
%
\begin{equation}\quad n^{\eta}[\nu-\nu_n] =\frac{c_{ F}^{2\alpha
}}{\expec[W]} \sum_{i=1}^K \biggl[\int_{i-1}^{i} u^{-2\alpha
}\,du-i^{-2\alpha}\biggr] +\Theta(K^{-2\alpha})+o(1).
\end{equation}
Letting first $n\rightarrow\infty$ followed by $K\rightarrow\infty
$, we conclude that
%
%
\begin{equation} \lim_{n\rightarrow\infty} n^{\eta}[\nu_n-\nu]
=\zeta
\end{equation}
as required. The fact that $\zeta>-\infty$ follows from the fact that,
for $i\geq2$,
%
%
\begin{equation} 0\leq\int_{i-1}^{i} u^{-2\alpha}\,du-i^{-2\alpha
}\leq(i-1)^{-2\alpha}-i^{-2\alpha},
\end{equation}
which is a summable sequence.
\end{pf}

We conclude that, in the critical regime where $\nu=1$, we have
%
%
\begin{equation} \label{tildenun-asymp} \nu_n(\lambda)=1+\theta
n^{-\eta}+o(n^{-\eta}),
\end{equation}
where $\theta=\lambda+\zeta$. The parameter $\theta\in{\mathbb R}$
indicates
the location inside the critical window formed by
the weights $\bfw(\lambda)$.
Indeed, in the asymptotics for $\nu_n(\lambda)$ in~(\ref{tildenun-asymp}),
the fact that $\theta=\zeta+\lambda$ arises from $\nu_n(\lambda
)=(1+\lambda n^{-\eta})\nu_n$,
together with the sharp asymptotics of $\nu_n$ in~(\ref{nun-asym}).
The value of $\zeta$ is constant and
does not depend on $\lambda$, while the value of $\lambda$ indicates
the location inside the scaling window,
so we can, alternatively, measure the location inside the scaling
window by $\theta\in{\mathbb R}$.

We continue by computing first and second moments of
total progenies and their weights, where, for our marked mixed Poisson
branching processes,
we define the \textit{weight of the branching process total progeny} to be
%
%
\begin{equation} \label{weight-BP}
w_T=\sum_{l=1}^T w_{M_l}
\end{equation}
and similar for $w_{T(i)}$.
Then we can compute the following moments, the proof of which
is standard and shall be omitted:
\begin{Lemma}[(Branching process characteristics)]
\label{lem-mom-cluster-weight}
\begin{longlist}[(a)]
\item[(a)]
%
%
\begin{equation} \expec[T]=\frac{1}{1-\nu_n},\qquad \expec
[T^2]=\frac{1+\nu_n}{(1-\nu_n)^2}+\frac{1}{(1-\nu_n)^3}\sum
_{j\in[n]} \frac{w_j^3}{l_n}.
\end{equation}

\item[(b)]
%
%
\begin{equation} \expec[w_{T}]=\frac{\nu_n}{1-\nu_n},\qquad
\expec[w_{T}^2]=\frac{1}{(1-\nu_n)^3}\sum_{j\in[n]} \frac
{w_j^3}{l_n}.
\end{equation}

\item[(c)]
%
%
\begin{eqnarray}
\expec[T(i)]&=&1+\frac{w_i}{1-\nu_n},\nonumber\\[-8pt]\\[-8pt]
\expec
[T(i)^2]&=&\biggl(1+\frac{w_i}{1-\nu_n}\biggr)^2+\frac{w_i(1+\nu
_n)}{(1-\nu_n)^2}+\frac{w_i}{(1-\nu_n)^3}\sum_{j\in[n]} \frac
{w_j^3}{l_n}.\nonumber
\end{eqnarray}

\item[(d)]
%
%
\begin{equation}\hspace*{28pt}
\expec\bigl[w_{T(i)}\bigr]=\frac{w_i}{1-\nu_n},\qquad
\expec\bigl[w_{T(i)}^2\bigr]=\biggl(\frac{w_i}{1-\nu_n}\biggr)^2+ \frac
{w_i}{(1-\nu_n)^3}\sum_{j\in[n]} \frac{w_j^3}{l_n}.
\end{equation}
%
\end{longlist}
\end{Lemma}


\subsection{Scaling limit of the cluster exploration process}
\label{sec-scal-lim-cluster-explor}

Theorem~\ref{thm-WC-C1} will follow from the fact that
we can identify the scaling limit of the process $(Z_l)_{l\geq0}$.
To do so, we let
%
%
\begin{equation} \label{Ztn-def} \ZZ_t^{(n)}=n^{-1/(\tau-1)}
Z_{tn^{(\tau-2)/(\tau-1)}} =n^{-\alpha} Z_{tn^{\rho}},
\end{equation}
where we recall the abbreviations in~(\ref{defs-al-rho-eta}).
By convention, for $t\geq0$ and for a discrete-time process
$(S_l)_{l\geq0}$, we let $S_t=S_{\lfloor t\rfloor}$.

The intuition behind~(\ref{Ztn-def}) is as follows. First, since
the largest connected components are of order $n^{\rho}$ as
proved in~\cite{Hofs09a}, Theorem 1.2, and the successive elapsed times
between hits of zero
of the process $(Z_l)_{l\geq0}$ correspond to the cluster sizes, the
relevant time scale is
$tn^{\rho}$. Further, by Theorem~\ref{prop-max-clusters-high-weight},
we see
that the large clusters correspond to the clusters of the high-weight vertices.
The maximal weight is of the order $n^{\alpha}$, so that this needs
to be the relevant scale on which the process $Z_l$ runs.
The proof below makes this intuition precise.

In order to define the scaling limit, we introduce a nonnegative
continuous-time process $(\SSS_t)_{t\geq0}$.
For some $a>0$, we let $(\II_i(t))_{i=1}^{\infty}$ denote independent
increasing
indicator processes defined by
%
%
\begin{equation} \label{Iit-def} \II_i(s)=\indic{\Exp(a i^{-\alpha
})\in[0,s]},\qquad s\geq0,
\end{equation}
so that
%
%
\begin{equation} \prob\bigl(\II_i(s)=0\mbox{ }\forall s\in[0,t]\bigr)={\mathrm
e}^{-ati^{-\alpha}}.
\end{equation}
We further let, for some $b>0$ and $c\in{\mathbb R}$, and $a$ as in
(\ref{Iit-def}),
%
%
\begin{equation} \label{SSt-def-abc} \SSS_t=b-abt+ct+\sum
_{i=2}^{\infty} b i^{-\alpha} [\II_i(t)- at i^{-\alpha}]
\end{equation}
for all $t\geq0$. We call $(\SSS_t)_{t\geq0}$ a
\textit{thinned L\'evy process}, a name we shall explain in more detail
after the theorem. To make the dependence on $(a,b,c)$ explicit, we now
denote $\SSS_t=\SSS_t(a,b,c)$.
Then, we have the obvious scaling relation
%
%
\begin{equation} \label{SS-scal-rel} \SSS_t(a,b,c)=b\SSS
_{at}\bigl(1,1,c/(ab)\bigr),
\end{equation}
where
%
%
\begin{eqnarray} \SSS_t(1,1,\beta) &=&1+(\beta-1)t+\sum_{i=2}^{\infty
} i^{-\alpha} [\II_i(t)-t i^{-\alpha}],\nonumber\\[-8pt]\\[-8pt]
\II_i(t)&=&\indic
{\Exp(i^{-\alpha})\in[0,t]}.\nonumber
\end{eqnarray}
%

The main result concerning the scaling limit of the exploration process
is the following theorem:
\begin{Theorem}[(The scaling limit of $Z_l$)]
\label{thm-weak-lim-S}
As $n\rightarrow\infty$, under the conditions of Theorem~\ref{thm-WC-3,4},
%
%
\begin{equation} \bigl(\ZZ_t^{(n)}\bigr)_{t\geq0}\convd(\SSS_t)_{t\geq
0},
\end{equation}
where $a=c_{ F}^{\alpha}/\expec[W]$,
$b=c_{ F}^{\alpha}$, $c=\theta-ab$,
in the sense of convergence in the $J_1$-Skorokhod topology on the
space of c\`adl\`ag functions on ${\mathbb R}^+$.
\end{Theorem}

It is worthwhile to note that while the convergence in Theorem \ref
{thm-weak-lim-S} only
has implications for our random graph for $t\leq H_1(0)$, which is the hitting
time of zero of the process $(\SSS_t)_{t\geq0}$, the processes
$(\ZZ_t^{(n)})_{t\geq0}$ and $(\SSS_t)_{t\geq0}$ are well
defined also for larger $t$, and convergence holds for \textit{all}
$t$. This is, in fact, useful in the proof.

The proof of Theorem~\ref{thm-weak-lim-S} shall be given in
Section~\ref{sec-pf-thm-weak-lim-S} below. We now first discuss
the limiting process $(\SSS_t)_{t\geq0}$ and its connection
to L\'evy processes. To do this, we
denote by $(\RR_t)_{t\geq0}$ the process given by
%
%
\begin{equation} \label{RRt-def-abc} \RR_t=b-abt+ct+\sum
_{i=2}^{\infty} b i^{-\alpha} [N_i(t)- at i^{-\alpha}],
\end{equation}
where $(N_i)_{t\geq0}$ are independent Poisson processes with rates
$ai^{-\alpha}$.
Clearly, the process $(\RR_t)_{t\geq0}$ is a spectrally positive L\'
evy process,
that is, $(\RR_t)_{t\geq0}$ has no negative jumps (see, e.g.,
\cite{Bert96,Kypr06} for more information on L\'evy processes),
with exponent $\psi(\varthetas)$ [for which $\E({\mathrm
e}^{-\varthetas(\RR_t-\RR_0)})={\mathrm e}^{-t\psi(\varthetas)}$]
given by
%
%
\begin{equation} \label{Psi-def} \psi(\varthetas)=(c-ab)\varthetas
+\sum_{i=2}^{\infty} a i^{-\alpha}[1-{\mathrm e}^{-\varthetas
b i^{-\alpha}}-b\varthetas i^{-\alpha}].
\end{equation}
Alternatively, the exponent $\psi(\varthetas)$ can be expressed as
%
%
\begin{eqnarray} \label{psi-def}
\psi(\varthetas)&=&(c-ab)\varthetas
-\varthetas\int_1^\infty x\Pi(dx)\nonumber\\[-8pt]\\[-8pt]
&&{}+\int_{0}^{\infty} \bigl(1-{\mathrm
e}^{-\varthetas x}-\varthetas x \indic{x<1}\bigr)\Pi(dx),\nonumber
\end{eqnarray}
where the L\'evy measure $\Pi$ is defined by
%
%
\begin{equation} \label{Pi-def} \Pi(dx)=\sum_{i=2}^{\infty} a
i^{-\alpha}\delta_{x, b i^{-\alpha}}.
\end{equation}
Since $\Pi(b,\infty)=0$, the jumps of $(\RR_t)_{t\geq0}$ are
bounded by $b$.
Further,
%
%
\begin{equation} \int_0^\infty(1\wedge x^2) \Pi(dx)\leq\int
_0^\infty x^2 \Pi(dx)=a\sum_{i=2}^{\infty} \biggl(\frac
{b}{i^{\alpha}}\biggr)^3 =ab^3 \sum_{i=2}^{\infty} i^{-3\alpha
}<\infty,\hspace*{-28pt}
\end{equation}
since $\tau\in(3,4)$ so that $3\alpha=3/(\tau-1)>1$. Therefore, the process
$(\RR_t)_{t\geq0}$ is a well-defined L\'evy process.

We may reformulate~(\ref{SSt-def-abc}) as
%
%
\begin{equation} \label{SSt-def-abc-rep} \SSS_t=b-abt+ct+\sum
_{i=2}^{\infty} b i^{-\alpha} \bigl[\indic{N_i(t)\geq1}-at i^{-\alpha
}\bigr],
\end{equation}
so that the process $(\SSS_t)_{t\geq0}$ does not include multiple
counts of
the independent processes $(N_i(t))_{t\geq0}$. This is the reason that
we call the process $(\SSS_t)_{t\geq0}$ a \textit{thinned} L\'evy process.
In~\cite{AldLim98}, this process is called a L\'evy process \textit
{without repetitions}.
Naturally, we have that the descriptions in~(\ref{RRt-def-abc}) and
(\ref{SSt-def-abc-rep}) satisfy that, a.s., for all $t\geq0$,
%
%
\begin{equation}
\label{stoch-dom-S-R} \SSS_t\leq\RR_t,
\end{equation}
which allows us to make use of L\'evy process methodology in our
proofs. We do note that $\RR_t$ is a rather poor approximation for
$\SSS_t$, particularly on large time scales, because the \textit
{thinning} becomes more important as time progresses.

\section{\texorpdfstring{Proofs of Theorems \protect\ref{thm-WC-C1} and \protect\ref{thm-weak-lim-S}}
{Proofs of Theorems 2.1 and 2.4}}
\label{sec-pf-thm-weak-lim-S}
In this section, we prove Theorems~\ref{thm-WC-C1} and~\ref{thm-weak-lim-S}.
We start by proving Theorem~\ref{thm-weak-lim-S} in Section \ref
{subsec-pf-thm-weak-lim-S},
and make use of Theorem~\ref{thm-weak-lim-S} to prove Theorem \ref
{thm-WC-C1}
in Section~\ref{sec-pf-thm-WC-C1}.

\subsection{\texorpdfstring{Proof of Theorem \protect\ref{thm-weak-lim-S}}{Proof of Theorem 2.4}}
\label{subsec-pf-thm-weak-lim-S}
Instead of $(Z_l)_{l\geq0}$, it is convenient to work with
a related process $(S_l)_{l\geq0}$, which is defined as $S_0=1,
S_1=w_1(\lambda)$ and satisfies the recursion relation, for $l\geq2$,
%
%
\begin{equation} \label{Sl-recurs} S_l=S_{l-1}+w_{M_l}(\lambda)J_l-1,
\end{equation}
that is, the Poisson random variables $\Poi(w_{M_l}(\lambda))$
appearing in
the recursion for $Z_l$ in~(\ref{Zl-recur-MP}) are replaced with
their (random) weights $w_{M_l}(\lambda)$. We shall first show that $S_l$
and $Z_l$ are quite close:
\begin{Lemma}
\label{lem-approxSZ}
Uniformly in $m\geq0$,
%
%
\begin{equation}
\sup_{l\leq m}|Z_{l}-S_{l}|=\Op(m^{1/2}).
\end{equation}
\end{Lemma}
\begin{pf}
We have that $(Z_l-S_l)_{l\geq0}$ is a martingale w.r.t.
the filtration $\mathcal{F}_l=\sigma((M_i)_{i=1}^l)$. Therefore, by the
Doob--Kolmogorov inequality (\cite{GriSti01}, Theorem~(7.8.2), page 338)
for any $M>0$,
%
%
\begin{equation}
\prob\Bigl(\sup_{l\leq m} |Z_{l}-S_{l}|>M\sqrt
{m}\Bigr) \leq\frac{1}{mM^2}\expec[|Z_{m}-S_{m}|^2].
\end{equation}
Now,
%
%
\begin{eqnarray}
\expec[|Z_{m}-S_{m}|^2] &=&\expec
\bigl[\expec[|Z_{m}-S_{m}|^2\mid(M_i)_{i=1}^m]\bigr] =\expec
\Biggl[\sum_{l=1}^m w_{M_l}(\lambda)J_l\Biggr]\nonumber\\[-8pt]\\[-8pt]
&\leq&\expec
\Biggl[\sum_{l=1}^m w_{M_l}(\lambda)\Biggr] =m\nu_n(\lambda
)=m\bigl(1+o(1)\bigr)\nonumber
\end{eqnarray}
by~(\ref{tildenun-asymp}). This proves the claim.
\end{pf}

We proceed by investigating the scaling limit of $(S_l)_{l\geq1}$. For
this, we define
%
%
\begin{equation}
\label{Stn-def} \SSS_t^{(n)}=n^{-\alpha}
S_{tn^{\rho}},
\end{equation}
where we recall the rounding convention right below~(\ref{Ztn-def}).

We shall prove that, in the sense of convergence in the $J_1$-Skorokhod
topology on the
space of c\`adl\`ag functions on ${\mathbb R}^+$,
%
%
\begin{equation} \label{SS-conv} \bigl(\SSS_t^{(n)}\bigr)_{t\geq0}\convd
(\SSS_t)_{t\geq0},
\end{equation}
which shall be enough to prove Theorem~\ref{thm-weak-lim-S}. Indeed,
to see
that~(\ref{SS-conv}) implies Theorem~\ref{thm-weak-lim-S}, we note that
by Lemma~\ref{lem-approxSZ}, for every $t=o(n^{(4-\tau)/(\tau-1)})$,
%
%
\begin{equation} \sup_{s\leq t} \bigl|\ZZ_s^{(n)}-\SSS_s^{
(n)}\bigr|=\Op\bigl(\sqrt{t} n^{(\tau-4)/(2(\tau-1))}\bigr) =o_{
\prob}(1).
\end{equation}

We continue with the proof of~(\ref{SS-conv}).
We shall prove that, due to~(\ref{cluster1-comp}) and
Lem\-ma~\ref{lem-approxSZ}, the first hitting time of $\SSS_s^{
(n)}$ of 0
is close to $n^{-\rho} |\cluster_{\leq}(1)|$.
We note that, by~(\ref{Sl-recurs}),
%
%
\begin{eqnarray} \label{Sl-rep} S_l&=&w_1(\lambda)+\sum_{i\in\mathcal
{V}_l^{(n)}} w_i(\lambda) -(l-1) \nonumber\\[-8pt]\\[-8pt]
&=&w_1(\lambda)+\sum_{i=2}^n
w_i(\lambda)\II^{(n)}_i(l) -(l-1),\nonumber
\end{eqnarray}
where
%
%
\begin{equation} \label{Vcal-def} \II^{(n)}_i(l)=\indic{i\in
\mathcal{V}_l^{(n)}}\qquad \mbox{with}\qquad \mathcal
{V}_l^{(n)}=\bigcup_{j=2}^l \{M_j\}.
\end{equation}
Using that
%
%
\begin{equation} \nu_n(\lambda)=\sum_{i\in[n]} \frac{w_i(\lambda)
w_i}{\ell_n},
\end{equation}
we can rewrite $S_l$ as
%
%
\begin{eqnarray} \label{Sl-rewrite}
S_l &=&w_1(\lambda)-\frac
{(l-1)w_1(\lambda) w_1}{\ell_n}+\sum_{i=2}^n w_i(\lambda)\biggl[\II
^{(n)}_i(l)-\frac{(l-1)w_i}{\ell_n}\biggr]\nonumber\\[-8pt]\\[-8pt]
&&{} +\bigl(\nu_n(\lambda
)-1\bigr)(l-1).\nonumber
\end{eqnarray}
Now we take $l=tn^{\rho}$, use that $\nu_n(\lambda)-1=\theta
n^{-\eta}\nu_n+o(n^{-\eta})$
[recall~(\ref{tildenun-asymp}) and~(\ref{defs-al-rho-eta})],
and we recall from~(\ref{choicewi}) and~(\ref{F-bound-tau3,4b}) that,
for $i$ such that
$n/i\rightarrow\infty$,
%
%
\begin{equation} w_i=[1-F]^{-1}(i/n)=b(n/i)^{\alpha}\bigl(1+o(1)\bigr),
\end{equation}
where $b=c_{ F}^{\alpha}$ and $c_{ F}$ is defined
in~(\ref{F-bound-tau3,4b}). As a result, by~(\ref{Sl-rewrite}),
%
%
\begin{eqnarray}
\label{SStn-rewrite}
\SSS_t^{(n)}&=&n^{-\alpha}
S_{tn^{\rho}}\nonumber\\
&=&b-\frac{b^2}{\expec[W]}t+\sum_{i=2}^n n^{-\alpha
}w_i(\lambda) \biggl[\II^{(n)}_i(tn^{\rho})-n^{-\alpha}\frac
{w_i t}{\mu_n}\biggr]\\
&&{} +\theta t +o(1),\nonumber
\end{eqnarray}
where we write $\mu_n=\ell_n/n=\expec[W]+o(1)$.

We proceed by showing that the sum in~(\ref{SStn-rewrite}) is predominantly
carried by the first few terms. Define
%
%
\begin{equation}
M_l^{(n,K)}= \sum_{i=K}^n n^{-\alpha
}w_i(\lambda) \biggl[\II^{(n)}_i(l)-\frac{(l-1)w_i}{\ell
_n}\biggr].
\end{equation}
We compute the mean and variance of $M_l^{(n,K)}$ for $K$ large.
For the mean, we compute
%
%
\begin{eqnarray}
\expec\bigl[M_l^{(n,K)}\bigr]&=& \sum_{i=K}^n n^{-\alpha
}w_i(\lambda) \biggl[\prob\bigl(\II^{(n)}_i(l)=1\bigr)-\frac
{(l-1)w_i}{\ell_n}\biggr] \nonumber\\[-8pt]\\[-8pt]
&=&\sum_{i=K}^n n^{-\alpha}w_i(\lambda)
\biggl[\biggl(1-\frac{w_i}{\ell_n}\biggr)^{l-1}-1
+\frac{(l-1)w_i}{\ell_n}\biggr].\nonumber
\end{eqnarray}
Thus, since $0\leq1-(1-x)^l-lx\leq(lx)^2/2$, we have that
$\expec[M_l^{(n,K)}]\leq0$ and
%
%
\begin{eqnarray}
\bigl|\expec\bigl[M_l^{(n,K)}\bigr]\bigr|
&=& \sum_{i=K}^n n^{-\alpha
}w_i(\lambda) \biggl[1-\biggl(1-\frac{w_i}{\ell_n}\biggr)^{l-1}-\frac
{(l-1)w_i}{\ell_n}\biggr] \nonumber\\
&\leq&\sum_{i=K}^n n^{-\alpha}w_i(\lambda)
\biggl(\frac{lw_i}{\ell_n}\biggr)^2
\leq C\frac{l^2n^{2\alpha
}}{\ell_n^2}\sum_{i=K}^n i^{-3\alpha} \\
&\leq& C\frac{l^2 n^{2\alpha
}}{\ell_n^2}K^{1-3\alpha},\nonumber
\end{eqnarray}
where, here and in the sequel, $C>0$ denotes a
constant that can change from line to line.
By~(\ref{defs-al-rho-eta}) and the fact that
$\ell_n=\Theta(n)$, we have that $ln^{\alpha}/\ell_n=\Theta(l
n^{\alpha-1})=\Theta(ln^{-\rho})$,
so that, uniformly in $l\leq tn^{\rho}$,
%
%
\begin{equation}
\label{mean-Ml-bd}
\bigl|\expec\bigl[M_l^{(n,K)}\bigr]\bigr| \leq
Ct^2K^{1-3\alpha}.
\end{equation}

To compute the variance of $M_l^{(n,K)}$, we
start by noting that $\II^{(n)}_i(l)$ is the indicator that
$i\in\mathcal{V}_l^{(n)}$, and $\mathcal{V}_l^{(n)}$
contains the
first $l$ marks drawn, where $M_1=1$ and the marks $(M_i)_{i=2}^l$ are
i.i.d.
with distribution given by~(\ref{dist-mark}). Therefore,
$\II^{(n)}_i(l)$ and $\II^{(n)}_j(l)$ are, for different $i,j$,
\textit{negatively correlated}, so that
%
%
\begin{equation}
\Var\bigl(M_l^{(n,K)}\bigr) \leq\sum_{i=K}^n
(n^{-\alpha}w_i(\lambda))^2 \Var\bigl(\II^{
(n)}_i(l)\bigr).
\end{equation}
Since $\II^{(n)}_i(l)$ is an indicator,
%
%
\begin{equation} \Var\bigl(\II^{(n)}_i(l)\bigr) \leq\expec
\bigl[\II^{(n)}_i(l)\bigr]\leq lw_i/\ell_n.
\end{equation}
Therefore, when $l=tn^{\rho}$, and using that $\rho+\alpha=1$
[recall~(\ref{defs-al-rho-eta})]
%
%
\begin{eqnarray}\label{var-Ml-bd}
\Var\bigl(M_l^{(n,K)}\bigr) &\leq&\sum_{i=K}^n
(n^{-\alpha}w_i(\lambda))^2 \frac{w_i tn^{\rho}}{\ell_n}
\leq Ct \sum_{i=K}^n i^{-3\alpha}\nonumber\\[-8pt]\\[-8pt]
&\leq& Ct K^{1-3\alpha} =o(1),\nonumber
\end{eqnarray}
when $K\rightarrow\infty$, since $\tau\in(3,4)$, so that $\alpha
=1/(\tau-1)>1/3$.\vadjust{\goodbreak}

We next observe that $(M_l^{(n,K)})_{l\geq1}$ is a
supermartingale, since
%
%
\begin{eqnarray}\quad
&&\expec\bigl[M_{l+1}^{(n,K)}-M_l^{(n,K)}\mid
\bigl(\II^{(n)}_i(l)\bigr)_{i\in[n]}\bigr]\nonumber\\
&&\qquad =\expec\Biggl[\sum
_{i=K}^n n^{-\alpha}w_i(\lambda) \biggl[\II^{(n)}_i(l+1)-\II
^{(n)}_i(l)-\frac{w_i}{\ell_n}\biggr]\Bigm|\bigl(\II^{
(n)}_i(l)\bigr)_{i\in[n]}\Biggr]\nonumber\\[-8pt]\\[-8pt]
&&\qquad \leq\sum_{i=K}^n n^{-\alpha
}w_i(\lambda) \bigl(1-\II^{(n)}_i(l)\bigr)\biggl(\expec\bigl[\II^{
(n)}_i(l+1)\mid\bigl(\II^{(n)}_i(l)\bigr)_{i\in[n]}\bigr]-\frac
{w_i}{\ell_n}\biggr)\nonumber\\
&&\qquad=0.\nonumber
\end{eqnarray}
Therefore, by the maximal inequality (\cite{GriSti01},
Theorem 12.6.1, page 496),
%
%
\begin{equation} \label{max-ineq}
\prob\Bigl(\max_{l\leq m}
\bigl|M_l^{(n,K)}\bigr|\geq\vep\Bigr) \leq\frac{-\expec[M_0^{
(n,K)}]+\expec[|M_m^{(n,K)}|]}{\vep}.
\end{equation}
We further bound, using Cauchy--Schwarz,
%
%
\begin{equation} \expec\bigl[\bigl|M_m^{(n,K)}\bigr|\bigr]
\leq\bigl|\expec\bigl[M_m^{
(n,K)}\bigr]\bigr|+\sqrt{\Var\bigl(M_l^{(n,K)}\bigr)}.
\end{equation}
Thus by~(\ref{mean-Ml-bd}) and~(\ref{var-Ml-bd}), and uniformly in
$m\leq tn^{\rho}$,
%
%
\begin{equation} \label{max-ineq-res}
\prob\Bigl(\max_{l\leq m}
\bigl|M_l^{(n,K)}\bigr|\geq\vep\Bigr) \leq Ct^2 \vep^{-1} K^{1-3\alpha
}+\vep^{-1} \sqrt{Ct K^{1-3\alpha}}.
\end{equation}
Since $\tau<4$, we obtain that, uniformly in $n$,
we can take $K=K(\vep)$ so large that
$\prob(\max_{l\leq m} |M_l^{(n,K)}|\geq\vep)\leq\vep$.

We denote, with $\mu_n=\ell_n/n$,
%
%
\begin{equation} \label{SStnK-def}\quad
\SSS_t^{(n,K)}=b-\frac
{b^2}{\expec[W]}t+\sum_{i=2}^{K} n^{-\alpha}w_i(\lambda) \biggl[\II
^{(n)}_i(tn^{\rho})-n^{-\alpha}\frac{w_i t}{\mu_n}
\biggr]+\theta t.
\end{equation}
%
Then we obtain the following corollary:
\begin{Corollary}[(Finite sum approximation of $\ZZ^{(n)}$)]
\label{cor-fin-sum}
For every $\vep, \delta,\allowbreak T>0$, there exists $K>0$ and $N\geq1$ such that
for all $n\geq N$,
%
%
\begin{equation} \label{fin-approx-Zn}
\prob\Bigl(\sup_{t\leq T}\bigl|\ZZ
_t^{(n)}-\SSS_t^{(n,K)}\bigr|\geq\delta\Bigr) \leq\vep.
\end{equation}
%
%
\end{Corollary}


The above suggests that it suffices to investigate
$(\II^{(n)}_i(t n^{\rho}))_{i\in[K]}$.
\begin{Lemma}[(Convergence of indicators)]
\label{lem-conv-indic}
As $n\rightarrow\infty$, for all $K\geq1$,
%
%
\begin{equation} \label{convd-indic}
\bigl(\II^{(n)}_i(t n^{\rho
})\bigr)_{i\in[K], t\geq0} \convd(\II_i(t))_{i\in[K], t\geq0}.
\end{equation}
As a consequence, for all $K\geq1$,
%
%
\begin{equation} \label{convdK}
\bigl(\SSS_t^{(n,K)}\bigr)_{t\geq0}\convd
\bigl(\SSS_t^{(\infty,K)}\bigr)_{t\geq0},
\end{equation}
where the limiting process $(\SSS_t^{(\infty,K)})_{t\geq0}$ is
defined as
%
%
\begin{equation} \label{SStK-def} \SSS_t^{(\infty,K)}=b-\frac
{b^2}{\expec[W]}t+\sum_{i=2}^{K} bi^{-\alpha} [\II
_i(t)-ai^{-\alpha}t]+\theta t.
\end{equation}
In both statements, $\convd$ refers to convergence in the
$J_1$-Skorokhod topology on the
space of c\`adl\`ag functions on ${\mathbb R}^+$.
\end{Lemma}
\begin{pf}
Convergence of the process for $t\geq0$ follows
when the process converges for $t\in[0,T]$ for all $T>0$
(see~\cite{Bill99}, Lemma 3, page 173).

Since $(\II^{(n)}_i(t n^{\rho}))_{t\geq0}$ are all indicator processes
of the form
%
%
\begin{equation} \II^{(n)}_i(t n^{\rho})= \indic{T_i\leq t
n^{\rho}},
\end{equation}
where $T_i$ is the first time that mark $i$ is chosen,
it suffices to prove that
%
%
\begin{equation}
(n^{-\rho} T_i)_{i\in[K]} \convd(E_i)_{i\in[K]},
\end{equation}
where $E_i$ are independent exponentials with rate $ai^{-\alpha}$.
For this, in turn, it suffices to prove that, for every sequence $t_1,
\ldots, t_K$,
%
%
\begin{equation} \prob(n^{-\rho} T_i>t_i\mbox{ }\forall i\in[K])
\rightarrow\exp\Biggl({-a\sum_{i=1}^K i^{-\alpha} t_i}\Biggr).
\end{equation}
The latter is equivalent to
%
%
\begin{eqnarray} \prob\bigl(\II^{(n)}_i(t_i n^{\rho})=0 \mbox{ }
\forall i\in[K]\bigr) &\rightarrow&\prob\bigl(\II_i(t_i)=0 \mbox{ } \forall i\in
[K]\bigr) \nonumber\\[-8pt]\\[-8pt]
&=&\exp\Biggl({-a\sum_{i=1}^K i^{-\alpha} t_i}\Biggr).\nonumber
\end{eqnarray}
Now, since the marks are i.i.d., we obtain that
%
%
\begin{eqnarray}
\prob\bigl(\II^{(n)}_i(m_i)=0 \mbox{ } \forall i\in
[K]\bigr) &=& \prod_{l=1}^\infty\prob(M_l \notin\{i\in[K]\dvtx l\leq
m_i\}) \nonumber\\[-8pt]\\[-8pt]
&=&\prod_{l=1}^\infty\biggl(1-\sum_{i\dvtx l\leq m_i} \frac
{w_i}{\ell_n}\biggr).\nonumber
\end{eqnarray}
A Taylor expansion gives that
%
%
\begin{eqnarray}
\prob\bigl(\II^{(n)}_i(m_i)=0 \mbox{ } \forall i\in
[K]\bigr) &=& \exp\Biggl({-\sum_{l=1}^n \sum_{i\dvtx m_i\geq l} \frac
{w_i}{\ell_n}+o(1)}\Biggr)\nonumber\\[-8pt]\\[-8pt]
&=&\exp\biggl({-\sum_{i\in[K]} \frac
{w_im_i}{\ell_n}+o(1)}\biggr).\nonumber
\end{eqnarray}
Applying this to $m_i= t_i n^{\rho}$, for which
%
%
\begin{equation} \frac{m_iw_i}{\ell_n}=\frac{bi^{-\alpha}
t_i}{\expec[W]}\bigl(1+o(1)\bigr),
\end{equation}
we arrive at the claim in~(\ref{convd-indic}) with $a=b/\expec[W]$.
The claim in
(\ref{convdK}) follows from the fact that, by~(\ref{SStnK-def}),
$\SSS_t^{(n,K)}$ is a weighted sum of the $(\II^{(n)}_i(t
n^{\rho}))_{i\in[K]}$,
and the (deterministic) weights converge. Thus, the continuous mapping theorem
gives the claim.
\end{pf}
\begin{pf*}{Proof of Theorem~\ref{thm-weak-lim-S}}
Again we use that convergence of the process for $t\geq0$ follows
when the process converges for $t\in[0,T]$ for all $T>0$
(see~\cite{Bill99}, Lemma 3, page 173). By~(\ref{fin-approx-Zn}),
with probability $1-o(1)$ when first $n\rightarrow\infty$ and
then $K\rightarrow\infty$, the process
$(\ZZ_t^{(n)})_{t\in[0,T]}$ is uniformly\vspace*{1pt} close to $(\SSS_t^{
(n,K)})_{t\in[0,T]}$.
By Lemma~\ref{lem-conv-indic}, the process $(\SSS_t^{
(n,K)})_{t\geq0}$
converges to $(\SSS_t^{(\infty,K)})_{t\geq0}$. Now,
%
%
\begin{equation} \label{SStK-diff} \SSS_t-\SSS_t^{(\infty
,K)}=\sum_{i\geq K+1} bi^{-\alpha} [\II_i(t)-ai^{-\alpha}
],
\end{equation}
and similar techniques as used to prove~(\ref{max-ineq-res}) can be
used to prove that
%
%
\begin{equation} \label{max-ineq-res-rep}\quad
\prob\Bigl(\max_{t\leq T}
\bigl|\SSS_t-\SSS_t^{(\infty,K)}\bigr|\geq\vep\Bigr) \leq CT^2 \vep
^{-1} K^{1-3\alpha}+\vep^{-1} \sqrt{CT K^{1-3\alpha}},
\end{equation}
so that again we can take $K=K(\vep)$ so large that
$\prob(\max_{t\leq T} |\SSS_t-\SSS_t^{(\infty,K)}|\geq\vep
)\leq\vep$.
This proves the claim.
\end{pf*}

\subsection{\texorpdfstring{Proof of Theorem \protect\ref{thm-WC-C1}}{Proof of Theorem 2.1}}
\label{sec-pf-thm-WC-C1}

In this section, we give a proof of Theorem~\ref{thm-WC-C1}.
We start by looking at the first hitting time of zero of the process
$l\mapsto Z_l$, and use the fact that by~(\ref{V1-def}), $V(1)=\inf
\{l\dvtx Z_l=0\}$,
where we recall that $V(1)$ denotes the number of
vertex checks performed in exploring the cluster of vertex 1.
Recall further that $\cluster(1)$ denotes the
cluster of vertex $1$, $|\cluster(1)|$ the number
of vertices in it, and $\weight(1)=\sum_{j\in\cluster(1)}w_j$
its weight.

The proof proceeds as follows. We shall first use
Theorem~\ref{thm-weak-lim-S} and Lem\-ma~\ref{lem-approxSZ} to prove that
$V(1)n^{-\rho}$ converges in distribution to $H_{\SSS}(0)$,
where $H_{\SSS}(0)$ denotes the first hitting time of 0 of the process
$(\SSS_t)_{t\geq0}$; see Corollary~\ref{cor-conv-hit-times}
below. We then prove that
$V(1)n^{-\rho}$, $|\cluster(1)|n^{-\rho}$ and $\weight(1)n^{-\rho}$
have identical scaling limits, by looking at the contribution due to the
second term in~(\ref{cluster1-comp}) for $|\cluster(1)|n^{-\rho}$,
and a similar computation for $\weight(1)n^{-\rho}$;
see Lemma~\ref{lem-multiple-hits} below.
We then complete the proof of Theorem~\ref{thm-WC-C1},
both for $|\cluster(1)|n^{-\rho}$ and for $\weight(1)n^{-\rho}$.
Finally, in Proposition~\ref{prop-wc-function},
we state and prove an auxiliary result concerning joint
convergence of $|\cluster(1)|n^{-\rho}$ and the indicators $\indic
{q\in\cluster(1)}$
for all $q$. This result is useful
in the proofs of Theorems~\ref{thm-WC-3,4} and~\ref{thm-prop-3,4}
and plays a crucial role in the proof of Theorem~\ref{thm-WC-C[K]} in the
next section, where we investigate the scaling limit of several
clusters simultaneously.

By Theorem~\ref{thm-weak-lim-S} and Lemma~\ref{lem-approxSZ}, the process
$(\ZZ_t^{(n)})_{t\geq0}$, where
$\ZZ_t^{(n)}=n^{-\alpha} Z_{tn^{\rho}}$, converges in
distribution to the process $(\SSS_t)_{t\geq0}$. By~(\ref{SS-conv}),
the same applies to $(\SSS_t^{(n)})_{t\geq0}$. Note that
%
%
\begin{equation} \label{HTT-V1}
n^{-\rho} V(1) =\min\bigl\{t\dvtx\ZZ
_t^{(n)}=0\bigr\}\equiv H^{(n)}(0).
\end{equation}
%
We next prove convergence in distribution of
$n^{-\rho} V(1)$:
\begin{Corollary}[(Convergence of hitting times)]
\label{cor-conv-hit-times}
As $n\rightarrow\infty$,
%
%
\begin{equation} n^{-\rho} V(1) \convd H_{\SSS}(0),
\end{equation}
where
%
%
\begin{equation} H_{\SSS}(x)=\inf\{t\dvtx\SSS_t\leq x\}
\end{equation}
is the first hitting time of level $x$ of $(\SSS_t)_{t\geq0}$.
\end{Corollary}
\begin{pf}
Since the process $(\SSS_t)_{t\geq0}$ has only positive jumps
(\cite{JacShi03}, Proposition~2.11 in Chapter 6) implies that the
hitting time of zero is a continuous function a.s. under the probability
measure of the limiting process on the space of c\`adl\`ag functions
equipped with the $J_1$-Skorokhod topology.
\end{pf}

\begin{Lemma}[($\SSS_t$ has a density)]
\label{lem-dens-SSt}
For all $t>0$, $\SSS_t$ has a density. As a result, the distribution of
$H_{\SSS}(0)$ has no atoms.
\end{Lemma}
\begin{pf}
We note that $\SSS_t$ has a density if and only if $\SSS'_t$ has, where
%
%
\begin{equation} \SSS'_t=\sum_{j=2}^{\infty} j^{-\alpha} [\II
_j'(t)-tj^{-\alpha}],
\end{equation}
and $(\II_j'(t))_{j\geq2}$ are independent indicator processes with
rate $j^{-\alpha}$.
This, in turn, follows when the characteristic function of $\SSS'_t$ is
integrable; see, for example,~\cite{GriSti01}, page 189.

The characteristic function of $\SSS'_t$ is given by
%
%
\begin{equation} \hat{f}_{\SSS'_t}(\varthetas)=\expec[{\mathrm
e}^{\i\varthetas\SSS'_t}] =\prod_{j=2}^\infty{\mathrm
e}^{-j^{-2\alpha}\i\varthetas}\bigl(1+({\mathrm e}^{-j^{-\alpha}\i
\varthetas}-1){\mathrm e}^{-j^{-\alpha}t}\bigr).
\end{equation}
Thus, for every $j_\varthetas\geq2$,
%
%
\begin{equation} |\hat{f}_{\SSS'_t}(\varthetas)| \leq\prod
_{j\geq j_\varthetas}^\infty|1+({\mathrm e}^{-j^{-\alpha}\i
\varthetas}-1){\mathrm e}^{-j^{-\alpha}t}|.
\end{equation}
Next, note that
\begin{eqnarray*}
&&
|1+({\mathrm e}^{-j^{-\alpha}\i\varthetas
}-1){\mathrm e}^{-j^{-\alpha}t}|^2 \\
&&\qquad= {\mathrm e}^{-2j^{-\alpha
}t} \sin(j^{-\alpha}\varthetas)^2 +\bigl(1-{\mathrm e}^{-j^{-\alpha
}t}+\cos(j^{-\alpha}\varthetas){\mathrm e}^{-j^{-\alpha}t}
\bigr)^2\\
&&\qquad= 1-2(1-{\mathrm e}^{-j^{-\alpha}t}){\mathrm
e}^{-j^{-\alpha}t}[1-\cos(j^{-\alpha}\varthetas)]\\
&&\qquad\leq
{\mathrm e}^{-2(1-{\mathrm e}^{-j^{-\alpha}t}){\mathrm
e}^{-j^{-\alpha}t}[1-\cos(j^{-\alpha}\varthetas)]} \\
&&\qquad\leq{\mathrm
e}^{-j^{-\alpha}t[1-\cos(j^{-\alpha}\varthetas)]},
\end{eqnarray*}
so that
%
%
\begin{equation} |\hat{f}_{\SSS'_t}(\varthetas)| \leq{\mathrm
e}^{-t\sum_{j\geq j_{\varthetas}} j^{-\alpha}t[1-\cos(j^{-\alpha
}\varthetas)]} \equiv{\mathrm e}^{-t \Phi(\varthetas)}.
\end{equation}
We choose
%
%
\begin{equation} j_{\varthetas}=\max\{j\geq2\dvtx b\varthetas
j^{-\alpha}\geq\pi/2\},
\end{equation}
so that
%
%
\begin{equation} \label{j-theta-form} j_{\varthetas}=\lfloor
(2b\varthetas/\pi)^{1/\alpha}\rfloor\vee2 =\lfloor(2b\varthetas
/\pi)^{\tau-1}\rfloor\vee2.
\end{equation}
Then we bound
%
%
\begin{equation} \Phi(\varthetas) \geq\sum_{j=j_{\varthetas
}}^{\infty} \frac{a}{j^\alpha}[1-\cos(b\varthetas j^{-\alpha})].
\end{equation}
Next, we use that
%
%
\begin{equation} 1-\cos(x)\geq\frac{2}{\pi}x^2,\qquad x\in\biggl[-\frac
{1}{2}\pi,\frac{1}{2}\pi\biggr],
\end{equation}
to arrive at
%
%
\begin{equation} \Phi(\varthetas) \geq c \varthetas^2 \sum
_{j=j_{\varthetas}}^{\infty} j^{-3\alpha},
\end{equation}
where $c>0$ denotes a positive constant appearing in
lower bounds that possibly
changes from line to line. We arrive at the fact that
%
%
\begin{equation} \Phi(\varthetas) \geq c \varthetas^2 j_\varthetas
^{1-3\alpha}\geq c \varthetas^2\vee\varthetas^{\tau-2},
\end{equation}
so that $|\hat{f}_{\SSS'_t}(\varthetas)|$ is integrable.
To prove that $H_{\SSS}(0)$ has no atoms,
note that when $\prob(H_{\SSS}(0)=u)>0$ for some $u\geq0$;
then, in particular, $\prob(\SSS_u=0)>0$, which contradicts the fact that
$\SSS_u$ has a density.
\end{pf}
%

We proceed by showing that the scaling limits of the number of vertex
checks of a
cluster and the cluster size are identical. For this, we shall make use
of the following lemma:
\begin{Lemma}[(Number of multiple hits is small)]
\label{lem-multiple-hits}
As $n\rightarrow\infty$, for any $m\geq1$,
%
%
\begin{equation} \label{expec-mult-hits} \expec\Biggl[\sum
_{j=2}^m[1-J_j]\Biggr]\leq\frac{mw_1}{\ell_n}+\frac{m(m-1)\nu
_n}{2\ell_n}.
\end{equation}
Consequently, there exists $t_n\rightarrow\infty$, such that
%
%
\begin{equation} \label{conv-mult-hits} n^{-\rho}\sum
_{j=2}^{t_nn^{\rho}} [1-J_j]\convp0.
\end{equation}
\end{Lemma}
\begin{pf}
We note that $J_j=0$ precisely when $M_l=1$ or when
there exists an $l<j$ and $i\in[n]$ such that $M_l=M_j=i$.
By independence and~(\ref{dist-mark}),
%
%
\begin{equation} \prob(M_l=M_j=i)=\prob(M_l=i)\prob
(M_j=i)=w_i^2/\ell_n^2.
\end{equation}
Therefore,
%
%
\begin{equation} \expec[1-J_j]\leq\frac{w_1}{\ell_n} + \sum
_{l=2}^{j-1} \sum_{i=2}^n \frac{w_i^2}{\ell_n^2} \leq\frac
{w_1}{\ell_n} +(j-1) \frac{\nu_n}{\ell_n}.
\end{equation}
%
Summing the above inequality over $2\leq j\leq m$ proves the claim in
(\ref{expec-mult-hits}).

For~(\ref{conv-mult-hits}), we use the Markov inequality to bound
\[
\prob\Biggl(n^{-\rho}\sum_{j=2}^{t_nn^{\rho}}
[1-J_j]\geq\vep_n\Biggr) \leq\vep_n^{-1} n^{-\rho}\expec
\Biggl[\sum_{j=2}^{t_nn^{\rho}} [1-J_j]\Biggr]\leq\frac{t_nw_1}{\vep_n
\ell_n}+\frac{t_n^2 n^{\rho} \nu_n}{2\ell_n\vep_n} =o(1),
\]
whenever $t_n^2 n^{-\alpha}/\vep_n=o(1)$.
Choosing, for example, $t_n=\log{n}$ and $\vep_n=1/\log{n}$ does the trick.
\end{pf}

Now we are ready to complete the proof of Theorem~\ref{thm-WC-C1}.
\begin{pf*}{Proof of Theorem~\ref{thm-WC-C1}}
By Corollary~\ref{cor-conv-hit-times},
$n^{-\rho} V(1)\convd H_{\SSS}(0)$. In particular, this implies
that $|\cluster(1)|\leq V(1)\leq n^{\rho}t_n$
for any $t_n\rightarrow\infty$. Therefore, by~(\ref{cluster1-comp}),
and \textit{whp},
%
%
\begin{equation} \label{sandwich-cluster1} n^{-\rho}V(1) -n^{-\rho
}\sum_{j=2}^{t_nn^{\rho}} [1-J_j] \leq n^{-\rho}|\cluster(1)|\leq
n^{-\rho}V(1).
\end{equation}
Now, by Lemma~\ref{lem-multiple-hits}, the difference between the
left-hand and right-hand sides of
(\ref{sandwich-cluster1}) converges to zero in probability, so that also
%
%
\begin{equation} n^{-\rho} |\cluster(1)|\convd H_{\SSS}(0).
\end{equation}
This completes the proof of Theorem~\ref{thm-WC-C1} and identifies
$H_1(0)=H_{\SSS}(0)$. In the same vein,
%
%
\begin{equation} \weight(1)=\sum_{i\in\cluster(1)} w_i =\sum
_{j=1}^{V(1)} w_{M_j} J_j.
\end{equation}
Now, by~(\ref{Sl-recurs}), for any $l\geq1$,
%
%
\begin{equation} \sum_{j=1}^{l} w_{M_j} J_j =S_l+l.
\end{equation}
As a result,
%
%
\begin{equation} \weight(1)=V(1)+S_{V(1)},
\end{equation}
so that
%
%
\begin{equation} n^{-\rho} \weight(1)=n^{-\rho}V(1)+n^{-\rho
}S_{V(1)}.
\end{equation}
Finally, $n^{-\rho}V(1)\convd H_{\SSS}(0)$, and,
since $\alpha<\rho$,
$n^{-\rho}|S_{V(1)}|=o(1)n^{-\alpha}|S_{V(1)}|\convp0$.
This proves that $n^{-\rho} \weight(1)\convd H_{\SSS}(0)$ as well.
\end{pf*}

In the next section, where we study the joint convergence of various
clusters simultaneously,
we shall also need the following joint convergence result:

\begin{Proposition}[(Weak convergence of functionals)]
\label{prop-wc-function}
As $n\rightarrow\infty$,
%
%
\begin{equation} \label{extend-WC-C1}
\bigl(n^{-\rho}\bigl|\cluster
(1)\bigr|, \bigl(\indic{q\in\cluster(1)}\bigr)_{q\geq1}\bigr) \convd
(H_{\SSS}(0), (\II_q(H_{\SSS}(0)))_{q\geq1})
\end{equation}
in the product topology, where $\II_q(H_{\SSS}(0))$ denotes
the indicator that $\II_q(t)=1$ at the hitting time of 0
of $(\SSS_t)_{t\geq0}$. Moreover, \textup{(i)}
the random variable $H_{\SSS}(0)$ is nondegenerate; and
\textup{(ii)} the indicators
$(\II_q(H_{\SSS}(0)))_{q\geq2}$ are nontrivial in the sense
that they take
the values 0 and 1 each with positive probability.
\end{Proposition}

We note that, while the indicator processes $(\II_q(t))_{t\geq0}$ are
independent
for different $q$, the random variables $(\II_q(H_{\SSS
}(0)))_{q\geq1}$ are \textit{not
independent} since $H_{\SSS}(0)$, the hitting time of 0 of the
process $(\SSS_t)_{t\geq0}$, depends
sensitively on all of the indicator processes.
%
\begin{pf}
We shall use a randomization trick.
Indeed, let $(N_j^{(n)}(t))_{t\geq0}$ be a sequence of
independent Poisson
processes with rate $w_j/\ell_n$. Let
%
%
\begin{equation} \label{Tj-def} T_j=\inf\{t\dvtx N(t)=j\}\qquad
\mbox{where}\qquad N(t)=\sum_{j\in[n]} N_j^{(n)}(t).\vadjust{\goodbreak}
\end{equation}
Then $t\mapsto N(t)$ is a rate 1 Poisson process, and
we have that [recall~(\ref{Sl-rewrite})]
%
%
\begin{equation} S_l=S'_{T_l},
\end{equation}
where the continuous-time process $(S_t')_{t\geq0}$ is defined by
%
%
\begin{eqnarray}
\label{Sprime-def}
S_t' &=& w_1(\lambda)-\frac
{w_1w_1(\lambda)N(t)}{\ell_n}+\sum_{i=2}^n w_i(\lambda)\biggl[\indic
{N_i^{(n)}(t)\geq1}-\frac{w_i N(t)}{\ell_n}\biggr]\nonumber\\[-8pt]\\[-8pt]
&&{} +\bigl(\nu
_n(\lambda)-1\bigr)N(t).\nonumber
\end{eqnarray}

By construction, the processes
$(\indic{N_q^{(n)}(n^{\rho}t)\geq1})_{t\geq0}$ are \textit
{independent}, and
are characterized by the birth times
%
%
\begin{equation} \label{birth-time-a} E_q^{(n)}=\inf\bigl\{t\dvtx
N_q^{(n)}(n^{\rho}t)\geq1\bigr\}.
\end{equation}
Again\vspace*{1pt} by construction, these birth times are independent
for different $q\geq2$, and $E_q^{(n)}$ has an exponential distribution
with parameter $n^{\rho}w_q/\ell_n$. The parameters of these
exponential random variables converge to
%
%
\begin{equation} \label{parameter-conv} n^{\rho}w_q/\ell
_n\rightarrow a q^{-\alpha},
\end{equation}
where $a=c_{ F}^{\alpha}/\expec[W]$, and which are the parameters
of the
limiting exponential random variables in terms of
which we can identify $\II_q(t)=\indic{N_q(t)\geq1}=\indic{
\Exp(a q^{-\alpha})\leq t}$; see~(\ref{SSt-def-abc-rep}). By the
convergence of the parameters, we can \textit{couple}
$E_q^{(n)}$ with $E_q=\Exp(a q^{-\alpha})$ in such a way that,
for every $q\geq2$ fixed,
%
%
\begin{equation} \label{coupling-exp} \prob\bigl(E_q^{(n)}\neq
E_q\bigr)=o(1).
\end{equation}
Indeed,~(\ref{coupling-exp}) follows by noting that,
by~(\ref{parameter-conv}), the density of $E_q^{(n)}$ converges
pointwise to that of
$E_q$, which, by~\cite{Thor00}, (7.3), implies that we can couple
$(E_q^{(n)})_{n\geq1}$ to $E_q$ in such a way that (\ref
{coupling-exp}) holds.

Equation~(\ref{coupling-exp}), jointly with the \textit{independence}
of $(E_q^{(n)})_{n\geq1}$ for different $q$'s, immediately implies that,
for each $K\geq1$,
%
%
\begin{equation} \label{coupling-indic} \prob\bigl(\indic{N_q^{
(n)}(n^{\rho}t)\geq1} =\II_q(t) \mbox{ } \forall t\geq0, q\in[K]
\bigr)=1-o(1),
\end{equation}
so that we have also, \textit{whp}, perfectly coupled the entire
processes
\[
\bigl(\indic{N_q^{(n)}(n^{\rho}t)\geq1}\bigr)_{t\geq0, q\in
[K]}\quad\mbox{and}\quad(\II_q(t))_{t\geq0, q\in[K]}.
\]
In particular, this implies that, for every $K\geq2$,
%
%
\begin{equation} \label{exact-coupling} \prob\bigl(\indic{N_q^{
(n)}(T_l)\geq1} =\II_q(T_l)\mbox{ }\forall l\geq1, q\in[K]\bigr)=1-o(1)
\end{equation}
and, by construction, $\indic{N_q^{(n)}(T_l)\geq1}=\II^{(n)}_q(l)$.

Applying\vspace*{1pt} the perfect coupling to $l=V(1)$, for which
$\indic{N_q^{(n)}(T_l)\geq1}=\indic{q\in\cluster(1)}$,
this provides a perfect coupling between $\indic{q\in\cluster(1)}$ and
$\II_q(T_{V(1)})$. We then note that
%
%
\begin{equation} n^{-\rho}|\cluster(1)|\convd H_{\SSS}(0)
\end{equation}
and, since $T_j$ is the birth time of the $j$th individual in a
rate 1 Poisson process,
%
%
\begin{equation} \sup_{t\leq u} |n^{-\rho}T_{tn^{\rho}}-t| \convp0,
\end{equation}
where, for noninteger $tn^{\rho}$, we recall the convention
below~(\ref{Ztn-def}).

Weak convergence of $(\indic{q\in\cluster(1)})_{q\geq1}$
in the product topology is equivalent to the weak convergence of
$(\indic{q\in\cluster(1)})_{q\in[m]}$ for any $m\geq1$; see
\cite{Kall02}, Theorem~4.29. Therefore,
together with the exact coupling in~(\ref{exact-coupling}),
this completes the proof of~(\ref{extend-WC-C1}), since
the processes $(\II_i(t))_{t\geq0}$ have a.s. no jump close to
$H_{\SSS}(0)$.

We continue to show the properties of the limiting variables.
The random variable $H_{\SSS}(0)$ is nondegenerate, since its
distribution does not have any atoms.
We shall next show that $\indic{q\in\cluster(1)}$ is nontrivial.
We shall show this only for $q=2$, the proof for $q>2$ being identical.
For this, we use the fact that
%
%
\begin{eqnarray}
\lim_{n\rightarrow\infty}\prob\bigl(2\in\cluster
(1)\bigr)&\geq&\prob\bigl(H_{\SSS}(0)\geq\vep, \II_2(\vep)=1\bigr)\nonumber\\[-8pt]\\[-8pt]
&\geq&
\prob\bigl(H_{\SSS}(0)\geq\vep\bigr)\prob\bigl(\II_2(\vep)=1\bigr)>0\nonumber
\end{eqnarray}
by the Fortuin--Kasteleyn--Ginibre (FKG) inequality (see
\cite{Grim99}, Theorem 2.4)
and the fact that both random variables $H_{\SSS}(0)$ and $\II
_2(\vep)$ are
monotone in the independent exponential random variables that
describe the first hit of $q$ for all $q\geq1$, so that
both $\{H_{\SSS}(0)\geq\vep\}$ and $\II_2(\vep)=1$ are
increasing events.

Further,
%
%
\begin{eqnarray}
\lim_{n\rightarrow\infty}\prob\bigl(2\notin\cluster
(1)\bigr)&\geq&\prob\bigl(H_{\SSS}(0)\leq K, \II_2(K)=0\bigr) \nonumber\\[-8pt]\\[-8pt]
&\geq&\prob
\bigl(H_{\SSS}(0)\leq K\bigr)\prob\bigl(\II_2(K)=0\bigr),\nonumber
\end{eqnarray}
again by FKG, now using that both
$\{H_{\SSS}(0)\leq K\}$ and $\{\II_2(K)=0\}$ are decreasing events.
Thus
%
%
\begin{equation} \lim_{n\rightarrow\infty}\prob\bigl(2\notin\cluster
(1)\bigr)\geq\prob\bigl(H_{\SSS}(0)\leq K\bigr)\prob\bigl(\II_2(K)=0\bigr)>0,
\end{equation}
which proves the claim.
\end{pf}
\begin{Remark}[(Convergence in the uniform topology)]
\label{rem-conv-unif}
In fact, by the proof of Proposition~\ref{prop-wc-function},
we even obtain that the weak convergence in Theorem~\ref{thm-WC-C1}
holds in the \textit{uniform topology}.
Indeed, the\vspace*{1pt} coupling obtained in the proof of Proposition \ref
{prop-wc-function}
[see in particular~(\ref{coupling-indic})] shows that
we can couple $(\SSS_t^{(n,K)})_{t\geq0}$
and $(\SSS_t^{(\infty,K)})_{t\geq0}$ such that these processes are
\textit{whp} equal \textit{for all} $t\geq0$. By
(\ref{fin-approx-Zn}) in Corollary~\ref{cor-fin-sum},
$(\ZZ_t^{(n)})_{t\geq0}$ is close to
$(\SSS_t^{(n,K)})_{t\geq0}$ in the uniform topology on $[0,T]$,
while~(\ref{max-ineq-res-rep}) shows that
$(\SSS_t^{(\infty,K)})_{t\geq0}$ is uniformly close to
$(\SSS_t)_{t\geq0}$. This proves the convergence in the uniform topology.
\end{Remark}
\begin{Remark}[(Convergence of cluster size of vertex $i$)]
\label{rem-conv-cluster-i}
We next remark on the scaling limits of $|\Ccal(i)|$ and $|\Ccal
_{\leq}(i)|$.
As in~(\ref{SSt-def-abc}), define
%
%
\begin{equation} \label{SSt-def-abc-i} \SSS_t^{(i)}=b-abt
i^{-\alpha}+ct+\sum_{j=1\dvtx j\neq i}^{\infty} b j^{-\alpha}
[\II_j(t)- at j^{-\alpha}],
\end{equation}
so that $(\SSS_t)_{t\geq0}=(\SSS_t^{(i)})_{t\geq0}$. Define
%
%
\begin{equation} H^{(i)}(0)=\inf\bigl\{t\dvtx\SSS_t^{(i)}=0\bigr\}
\end{equation}
to be the first hitting time of zero of the process $(\SSS_t^{
(i)})_{t\geq0}$.
Then, in an identical way as in the proof of Proposition \ref
{prop-wc-function}, it follows that,
as $n\rightarrow\infty$,
%
%
\begin{equation} \label{extend-WC-Ci}
\bigl(n^{-\rho}|\cluster
(i)|, \bigl(\indic{q\in\cluster(i)}\bigr)_{q\geq1}\bigr) \convd
\bigl(H^{(i)}(0), \bigl(\II_q\bigl(H^{(i)}(0)\bigr)\bigr)_{q\geq1}\bigr)
\end{equation}
in the product topology. As a result,
%
%
\begin{equation} n^{-\rho}|\cluster_{\leq}(i)|\convd H^{
(i)}(0)\prod_{j=1}^{i-1} \bigl(1-\II_q\bigl(H^{(i)}(0)\bigr)\bigr).
\end{equation}
\end{Remark}

%

\section{Convergence of multiple clusters}
\label{sec-mult-clusters}

In this section, we extend the analysis of one cluster in Section \ref
{sec-WC-C1}
to multiple clusters.
This sets the stage for the proof of Theorem~\ref{thm-WC-3,4}, which is completed in the next section.
The main result is as follows:
\begin{Theorem}[(Weak convergence of clusters of first vertices)]
\label{thm-WC-C[K]}
Fix the Norros--Reittu random graph with weights
$\bfw(\lambda)$ defined in~(\ref{w-lambda-def}).
Assume that $\nu=1$ and that~(\ref{F-bound-tau3,4b}) holds. Then,
for all $\lambda\in\Rbold$,
%
%
\begin{equation} (n^{-\rho}|\cluster_{\leq}(i)|)
_{i\geq1} \convd(H_i(0))_{i\geq1}
\end{equation}
for some nondegenerate limit $(H_i(0))_{i\geq1}$.
\end{Theorem}

In the remainder of this section, we shall prove Theorem~\ref{thm-WC-C[K]}
and use it to complete the proof of Theorem~\ref{thm-WC-3,4}.
We let $I_1^{(n)}=1$, and let
%
%
\begin{equation} I_2^{(n)}=\min[n]\setminus\cluster(1)\vadjust{\goodbreak}
\end{equation}
be the minimal element that is not part of $\cluster(1)$, where,
for a set of indices $A\subseteq[n]$, we let $\min A$ denote the
minimal element of $A$.
To extend the above definitions further, we define, recursively,
%
%
\begin{equation} \label{Di-def} \mathcal{D}_{i}^{(n)}=\cluster
_{\leq}\bigl(I_i^{(n)}\bigr)\quad \mbox{and} \quad\mathcal
{D}_{\leq i}^{(n)}=\bigcup_{j\leq i} \mathcal{D}_{j}^{
(n)}.
\end{equation}
Then we define $I_{i+1}^{(n)}$ by
%
%
\begin{equation} I_{i+1}^{(n)}=\min[n]\setminus\mathcal
{D}_{\leq i}^{(n)},
\end{equation}
which is the vertex with the smallest index of which we have not yet
explored its cluster.

Obviously,\vspace*{1pt} $|\cluster_{\leq}(i)|=0$ unless $i=I_j^{(n)}$
for some $j$. This
prompts us to investigate the\vspace*{1pt} weak convergence of $n^{-\rho}|\mathcal
{D}_{i}^{(n)}|$.
This will be done by induction on $i$. The \textit{induction hypothesis}
is that
%
%
\begin{equation} \label{IH-form}
\bigl(n^{-\rho}\bigl|\mathcal
{D}_{j}^{(n)}\bigr|, \bigl(\indic{q\in\mathcal{D}_{\leq j}^{
(n)}}\bigr)_{q\geq1}\bigr)_{j\in[i]} \convd\bigl(H_{j}(0), \bigl(\indic
{q\in\mathcal{D}_{\leq j}}\bigr)_{q\geq1}\bigr)_{j\in[i]}
\end{equation}
in the product topology, for some limiting random variables.
Part of the induction hypothesis is that these limiting random variables
satisfy the following facts: (1) the limiting random variables
$(H_{j}(0))_{j\in[i]}$
are \textit{nondegenerate}, in the sense that the essential support of
the random vector
$(H_{j}(0))_{j\in[i]}$ is $i$-dimensional, and (2) the random indicators
$(\indic{q\in\mathcal{D}_{\leq j}})_{j\in[i], q>i}$ are all
\textit{nontrivial}, in the sense that they take the values zero and
one, each with positive probability. By construction,
$\indic{q\in\mathcal{D}_{\leq j}^{(n)}}=1$ for $q\leq i$,
so the restriction
to $q>i$ in condition (2) is the most we can hope for.

We shall start by initializing the induction hypothesis for $j=1$,
which follows from
Proposition~\ref{prop-wc-function}, as we show now. Indeed, we have that
$\mathcal{D}_{1}^{(n)}=\mathcal{D}_{\leq1}^{
(n)}=\cluster(1)$, so that~(\ref{IH-form})
is identical to the statement in Proposition~\ref{prop-wc-function}.

We next advance the induction hypothesis by verifying that
(\ref{IH-form}) also holds for $j=i+1$. We first\vspace*{1pt} intuitively explain
our approach. The random variable $H_{i+1}(0)$ shall be the weak limit of
$n^{-\rho}|\mathcal{D}_{i+1}^{(n)}|$. We shall show that
$H_{i+1}(0)$ is the hitting time of zero of a process similar to $(\SSS
_t)_{t\geq0}$
in Section~\ref{sec-WC-C1}. We now start by explaining how this
process arises.

Assume that the induction hypothesis~(\ref{IH-form}) holds for $i$.
By~(\ref{IH-form}), the index set $\mathcal{D}_{\leq i}$ is the
(random) set of indices for which
%
%
\begin{equation} \bigl(\indic{q\in\mathcal{D}_{\leq i}^{
(n)}}\bigr)_{q\geq1}\convd\bigl(\indic{q\in\mathcal{D}_{\leq
i}}\bigr)_{q\geq1}.
\end{equation}
Then, we note that,
by~(\ref{IH-form}), we have that
%
%
\begin{equation} \label{conv-dist-index}\qquad
I_{i+1}^{(n)}\equiv\min
\bigl\{q\dvtx\indic{q\in\mathcal{D}_{\leq i}^{(n)}}=0\bigr\}
\quad\convd\quad I_{i+1}\equiv\min\bigl\{q\dvtx\indic{q\in\mathcal{D}_{
\leq i}}=0\bigr\},
\end{equation}
and we see that $I_{i+1}^{(n)}$ and $I_{i+1}$ are
\textit{deterministic} functions of the sets $\mathcal{D}_{\leq
i}^{(n)}$
and $\mathcal{D}_{\leq i}$, respectively. The random variable
$I_{i+1}$ is \textit{finite},
since, for $K,Q\geq1$ large,
%
%
\begin{equation} \label{split-I-large} \prob\bigl(I_{i+1}^{(n)}>K\bigr)
\leq\prob\bigl(\bigl|\mathcal{D}_{\leq j}^{(n)}\bigr|\geq Qn^{\rho}\bigr)
+\prob\bigl(I_{i+1}^{(n)}\geq K, \bigl|\mathcal{D}_{\leq j}^{
(n)}\bigr|< Qn^{\rho}\bigr).
\end{equation}
The first probability converges, by~(\ref{IH-form}) and the
continuous-mapping theorem, to
$\prob(H_1(0)+\cdots+H_i(0)\geq Q)$, which is small for $Q\geq1$
large. For the second probability in~(\ref{split-I-large}),
and for $i\leq K/2$, we can bound
%
%
\begin{eqnarray} \prob\bigl(I_{i+1}^{(n)}>K, \bigl|\mathcal{D}_{\leq
i}^{(n)}\bigr|< Qn^{\rho}\bigr) &\leq& \prob(\mbox{vertex } K \mbox{
drawn in } Qn^{\rho} \mbox{ vertex checks})\nonumber\\
&\leq& \prob(\exists
l\leq Qn^{\rho}\dvtx M_l=K) \leq\sum_{l=1}^{Qn^{\rho}} \frac
{w_K}{\ell_n} \\
&\leq& C Q K^{-\alpha},\nonumber
\end{eqnarray}
which converges to zero as $K\to\infty$
when $Q=K^{\beta}$ with $\beta<\alpha$.
As a result, we have that
%
%
\begin{equation} \prob(I_{i+1}>K)=\lim_{n\rightarrow\infty} \prob
\bigl(I_{i+1}^{(n)}>K\bigr)
\end{equation}
is small for $K$ large.

We conclude that, from the induction hypothesis in~(\ref{IH-form}),
we obtain the \textit{joint} convergence
%
%
\begin{eqnarray} \label{IH-conseq}
&&\bigl(\bigl(n^{-\rho}\bigl|\mathcal
{D}_{i}^{(n)}\bigr|, \bigl(\indic{q\in\mathcal{D}_{\leq j}^{
(n)}}\bigr)_{q\geq1}\bigr)_{j\in[i]},I_{i+1}^{(n)}\bigr)\nonumber\\[-8pt]\\[-8pt]
&&\qquad \convd
\bigl(\bigl(H_j(0), \bigl(\indic{q\in\mathcal{D}_{\leq j}}\bigr)_{q\geq1}
\bigr)_{j\in[i]}, I_{i+1}\bigr).\nonumber
\end{eqnarray}
We now start exploring the cluster of $I_{i+1}^{(n)}$,
and we need to show that this cluster size, as well as the indices in it,
converge in distribution. More precisely, the joint
convergence in~(\ref{IH-form}) for $i+1$ (and thus
the advancement of the induction hypothesis)
follows when we prove that, conditionally on
$\mathcal{D}_{\leq i}^{(n)}$,
%
%
\begin{eqnarray} \label{induct-aim}
&&\bigl(n^{-\rho}\bigl|\mathcal
{D}_{i+1}^{(n)}\bigr|, I_{i+1}^{(n)}, \bigl(\indic{q\in\mathcal
{D}_{\leq i+1}^{(n)}}\bigr)_{q\geq1}\bigr) \nonumber\\[-8pt]\\[-8pt]
&&\qquad\convd
\bigl(H_{i+1}(0), I_{i+1}, \bigl(\indic{q\in\mathcal{D}_{\leq
i+1}}\bigr)_{q\geq1}\bigr).\nonumber
\end{eqnarray}
To prove~(\ref{induct-aim}), we follow the approach in Section
\ref{sec-WC-C1} as closely as possible. A~crucial observation is that
after the exploration of $\mathcal{D}_{\leq i}^{(n)}$ and
conditionally on it, the
remaining graph is again a rank-1 inhomogeneous random graph, with
(a) vertex set $[n]\setminus\mathcal{D}_{\leq i}^{(n)}$, and
(b) edge\vspace*{1pt} probabilities, for $u,v\in[n]\setminus\mathcal{D}_{
\leq i}^{(n)}$, given by
$p_{uv}=1-{\mathrm e}^{-w_uw_v/\ell_n}$.

We now extend the exploration process of clusters described in
Section~\ref{sec-cluster-BP} to the setting above.
As in\vspace*{1pt} Section~\ref{sec-WC-C1}, we set
$Z_0(i)=1$ and let $Z_1(i)$ denote the number of neighbors of
the vertex $I_{i+1}^{(n)}$ outside $\mathcal{D}_{\leq
i}^{(n)}$, that is,
%
%
\begin{equation} Z_1(i)=\sum_{j\notin\mathcal{D}_{\leq
i}^{(n)}} \Poi\bigl(w_{I_{i+1}^{(n)}}(\lambda)w_j/\ell_n\bigr)
=\Poi\bigl(w_{I_{i+1}^{(n)}}(\lambda)\ell_n(i)/\ell_n\bigr),
\end{equation}
where we let
%
%
\begin{equation} \ell_n(i)=\sum_{j\notin\mathcal{D}_{\leq
i}^{(n)}} w_j
\end{equation}
be the total weight of vertices \textit{outside} $\mathcal{D}_{
\leq i}^{(n)}$.
For $l\geq2$, $(Z_l(i))_{l\geq1}$ satisfies the recursion relation
%
%
\begin{equation} Z_l(i)=Z_{l-1}(i)+X_l(i)-1,
\end{equation}
where $X_l(i)$ denotes the number of potential neighbors outside of
$\mathcal{D}_{\leq i}^{(n)}$ of the $l$th vertex
which is explored. As explained in more detail in
Section~\ref{sec-WC-C1},
the distribution of $X_l(i)$ (for $2\leq l\leq n$) is equal to $\Poi
(w_{M_{l}(i)}\ell_n(i)/\ell_n)J_{l}(i)$,
where now the marks $(M_{l}(i))_{l=1}^{\infty}$ are i.i.d. random
variables with
distribution $M(i)$ given by
%
%
\begin{equation} \label{dist-mark-i} \prob\bigl(M(i)=m\bigr)=
w_m/\ell_n(i),\qquad
m\in[n]\setminus\mathcal{D}_{\leq i}^{(n)},
\end{equation}
and
%
%
\begin{equation} J_{l}(i)=\indic{M_{l}(i)\notin\{I_{i+1}^{
(n)}\}\cup\{M_{2}(i), \ldots, M_{l-1}(i)\}}
\end{equation}
is the indicator that the mark $M_{l}(i)$ has not been found up to time
$l$ and is not equal to vertex $I_{i+1}^{(n)}$.


Then the number of vertex checks $V(I_{i+1}^{(n)})$ in the
exploration of
$\mathcal{D}_{i+1}^{(n)}=\cluster_{\leq}(I_{i+1}^{
(n)})$ equals
%
%
\begin{equation} \label{clusterIi+1-comp}
V\bigl(I_{i+1}^{(n)}\bigr)=\inf
\{l\dvtx Z_l(i)=0\}
\end{equation}
and
%
%
\begin{equation} \bigl(\indic{q\in\mathcal{D}_{i+1}^{(n)}}
\bigr)_{q\neq I_{i+1}^{(n)} =(\indic{\exists l\leq|\mathcal
{D}_{i+1}^{(n)}|\dvtx M_l(i)=q})_{q\neq I_{i+1}^{
(n)}}},
\end{equation}
while $\indic{I_{i+1}^{(n)}\in\mathcal{D}_{i+1}^{(n)}}=1$.
We again note that
%
%
\begin{equation} \bigl|\mathcal{D}_{i+1}^{(n)}\bigr|=\bigl|\cluster_{\leq
}\bigl(I_{i+1}^{(n)}\bigr)\bigr|\leq V\bigl(I_{i+1}^{(n)}\bigr),
\end{equation}
while
%
%
\begin{equation} n^{-\rho} \bigl[V\bigl(I_{i+1}^{(n)}\bigr)-\bigl|\mathcal
{D}_{i+1}^{(n)}\bigr|\bigr] \convp0,
\end{equation}
which can be proved along the lines of the proof of Lemma \ref
{lem-multiple-hits}.
This gives us a convenient description of all the random variables
needed to advance
the induction hypothesis.

In order to prove the weak convergence of $n^{-\rho}V(I_{i+1}^{(n)})$,
we again investigate the scaling limit of the process $(Z_l(i))_{l\geq
0}$. For this,
we define $S_0(i)=1$, $S_1(i)=w_{I_{i+1}^{(n)}}(\lambda)\ell
_n(i)/\ell_n$ and, for $l\geq2$,
%
%
\begin{equation} \label{Sli-def}
S_l(i)=S_{l-1}(i)+w_{M_{l}(i)}(\lambda)J_{l}(i)\ell_n(i)/\ell_n -1.
\end{equation}
Then, as in Lemma~\ref{lem-approxSZ}, it is easy to show that,
conditionally on $\mathcal{D}_{\leq i}^{(n)}$, the processes
$(S_l(i))_{l\geq0}$ and
$(Z_l(i))_{l\geq0}$ are uniformly close. Denote by $\mathcal{B}_{
i}^{(n)}=\mathcal{D}_{\leq i}^{(n)}\cup\{I_{
i+1}^{(n)}\}$ the union of all vertices explored in the first $i$
clusters and the minimal element not in the first $i$ clusters. We rewrite
%
%
\begin{eqnarray} 
S_l(i) &=& w_{I_{i+1}^{
(n)}}(\lambda)\frac{\ell_n(i)}{\ell_n} +\sum_{j=2}^l
w_{M_{j}(i)}\frac{\ell_n(i)}{\ell_n}J_{j}(i)
-(l-1)\nonumber\\[-8pt]\\[-8pt]
&=& w_{I_{i+1}^{(n)}}(\lambda)\frac{\ell_n(i)}{\ell_n} +\sum
_{q\in[n]\setminus\mathcal{B}_{ i}^{(n)}} w_q(\lambda
)\frac{\ell_n(i)}{\ell_n} \II_q^{(n)}(l;i)-(l-1),\nonumber
\end{eqnarray}
where
%
%
\begin{equation} \II_q^{(n)}(l;i)=\indic{\exists j\leq l\dvtx
M_j(i)=q}.
\end{equation}
We further rewrite the above as
%
%
\begin{eqnarray} \label{Sli-rewr} S_l(i)&=&w_{I_{i+1}^{
(n)}}(\lambda)\frac{\ell_n(i)}{\ell_n} +\sum_{q\in[n]\setminus
\mathcal{B}_{ i}^{(n)}} w_q(\lambda)\frac{\ell_n(i)}{\ell
_n} \biggl(\II_q^{(n)}(l;i)-\frac{lw_q}{\ell_n(i)}\biggr)\nonumber\\[-8pt]\\[-8pt]
&&{} +l
\biggl(\sum_{q\in[n]\setminus\mathcal{B}_{ i}^{(n)}} \frac
{w_q(\lambda) w_q}{\ell_n}-1\biggr)+1.\nonumber
\end{eqnarray}
We note that we can rewrite the last sum, using (\ref
{tildenun-asymp}), as
%
%
\begin{eqnarray} (1+\lambda n^{-\eta})\sum_{q\in[n]\setminus
\mathcal{B}_{ i}^{(n)}} \frac{w_q^2}{\ell_n}-1 &=& \bigl(\nu
_n(\lambda)-1\bigr)-(1+\lambda n^{-\eta})\sum_{q\in\mathcal{B}_{
i}^{(n)}} \frac{w_q^2}{\ell_n}\nonumber\\[-8pt]\\[-8pt]
&=& \theta n^{-\eta} -\sum
_{q\in\mathcal{B}_{ i}^{(n)}} \frac{w_q^2}{\ell_n}
+o(n^{-\eta}).\nonumber
\end{eqnarray}
In turn, the sum can be approximated by
%
%
\begin{equation} \sum_{q\in\mathcal{B}_{ i}^{(n)}} \frac
{w_q^2}{\ell_n} =d n^{-\eta} \sum_{q\in\mathcal{B}_{ i}^{
(n)}} q^{-2\alpha}\bigl(1+o_{\prob}(1)\bigr),
\end{equation}
where $d=c_{ F}^{2\alpha}/\expec[W]$. Denoting
%
%
\begin{equation} D_i^{(n)}=d \sum_{q\in\mathcal{B}_{
i}^{(n)}} q^{-2\alpha},
\end{equation}
we therefore have that
%
%
\begin{eqnarray} \label{Sli-rewr-rep2}
S_l(i)&=&w_{I_{i+1}^{
(n)}}\frac{\ell_n(i)}{\ell_n} +\sum_{q\in[n]\setminus\mathcal
{B}_{ i}^{(n)}} w_q \frac{\ell_n(i)}{\ell_n}\biggl(\II
_q^{(n)}(l;i)-\frac{lw_q} {\ell_n(i)}\biggr)\nonumber\\[-8pt]\\[-8pt]
&&{}+l\bigl(\theta-D_i^{
(n)}\bigr)n^{-\eta} + o_{\prob}(l n^{-\eta}).\nonumber
\end{eqnarray}
We conclude that we arrive at a similar process
as when exploring $\cluster(1)$, apart from the fact that:
(i) fewer vertices are allowed to participate,
(ii) a negative drift $-D_i^{(n)}$ is introduced and
(iii) a factor $\frac{\ell_n(i)}{\ell_n}=1+o_{\prob}(1)$ is introduced.

We proceed by investigating the convergence of $D_i^{(n)}$:
\begin{Lemma}[(Weak convergence of random drift)]
\label{lem-wc-drift}
As $n\rightarrow\infty$, and assuming~(\ref{IH-conseq}),
%
%
\begin{equation} D_i^{(n)}\convd D_i\equiv\sum_{q\in\mathcal
{D}_{\leq i}\cup\{I_{i+1}\}} q^{-2\alpha},
\end{equation}
where $(\mathcal{D}_{\leq i},I_{i+1})$ is the weak limit of
$(\mathcal{D}_{\leq i}^{(n)},I_{i+1}^{(n)})$ given in
(\ref{IH-conseq}).
\end{Lemma}
\begin{pf}
We start by bounding $\prob(q\in\mathcal{D}_{\leq i}^{
(n)})$, for
$q>0$ large. We shall first prove that the probability that
$|\mathcal{D}_{\leq i}^{(n)}|\leq n^{\rho} K$ is
$1-o(1)$ when\vspace*{1pt} $K>0$ grows large. Indeed, by~\cite{Hofs09a},
Theorem 1.2, we have that,
$|\Cmax|=\max_i|\cluster_{\leq i}|\leq\omega n^{\rho}$
with probability $1-o(1)$, as $\omega\rightarrow\infty$.
Thus, $|\mathcal{D}_{\leq i}^{(n)}|\leq n^{\rho} (i\omega
)=n^{\rho}K$,
with probability\vspace*{1pt} $1-o(1)$ as $K\rightarrow\infty$, when we take
$K=\omega i$. Denoting
%
%
\begin{equation} \mathcal{E}_{i,K}^{(n)}=\bigl\{\bigl|\mathcal{D}_{
\leq i}^{(n)}\bigr|\leq n^{\rho} K\bigr\},
\end{equation}
we have that
%
%
\begin{equation} \prob\bigl(\bigl\{q\in\mathcal{D}_{\leq i}^{
(n)}\setminus\bigl\{ I_{j}^{(n)}\bigr\}_{j=1}^i\bigr\}\cap\mathcal
{E}_{i,K}^{(n)}\bigr) \leq n^{\rho} K \frac{w_q}{\sum_{j>K n^{\rho
}} w_j},
\end{equation}
since, independently of the choices before, the probability of drawing
$q$ is at most
$w_q/\sum_{j>K n^{\rho}} w_j$. Now,
%
%
\begin{equation} \sum_{j>K n^{\rho}} w_j=\ell_n\bigl(1+o(1)\bigr)=\expec[W]
n\bigl(1+o(1)\bigr).
\end{equation}
Thus, for some $C>0$,
%
%
\begin{equation} \prob\bigl(\bigl\{q\in\mathcal{D}_{\leq i}^{
(n)}\setminus\bigl\{ I_{j}^{(n)}\bigr\}_{j=1}^i\bigr\}\cap\mathcal
{E}_{i,K}^{(n)}\bigr)\leq CK q^{-\alpha},\vadjust{\goodbreak}
\end{equation}
so that
%
%
\begin{equation} \label{expec-Di-tight} \expec\biggl[\sum_{q\in
\mathcal{B}_{ i}^{(n)}\dvtx q>Q} q^{-2\alpha}\indicwo
{\mathcal{E}_{i,K}^{(n)}}\biggr] \leq iQ^{-2\alpha}+CK
Q^{1-3\alpha},
\end{equation}
where\vspace*{1pt} the first contribution arises from the (at most $i$) values of
$q=I_{j}^{(n)}$ for $j\in[i+1]$
for which $I_{j}^{(n)}>Q$, and the second contribution from the
$q\notin\{I_{j}^{(n)}\}_{j\in[i+1]}$.

Equation~(\ref{expec-Di-tight}) implies that the weak convergence of
$D_i^{(n)}$ follows
from the weak convergence of
%
%
\begin{equation} \sum_{q\in\mathcal{B}_{ i}^{(n)}\dvtx
q\leq Q} q^{-2\alpha},
\end{equation}
which, in turn, follows from~(\ref{IH-conseq}) and the continuous
mapping theorem.
\end{pf}

Now we are ready to complete the proof of Theorem~\ref{thm-WC-C[K]}.
\begin{pf*}{Proof of Theorem~\ref{thm-WC-C[K]}}
We start by setting the stage for the weak convergence of processes
needed to
advance the induction hypothesis as formulated in~(\ref{induct-aim}).
Define
%
%
\begin{equation} \label{SStni-def}
\ZZ_t^{(n)}(i)=n^{-\alpha}
Z_{tn^{\rho}}(i),\qquad \SSS_t^{(n)}(i)=n^{-\rho} S_{tn^{\rho
}}(i)
\end{equation}
and
%
%
\begin{equation} \label{SSti-def} \quad\SSS_t(i)=bI_{i+1}^{-\alpha}
+\sum_{q\in\mathcal{D}_{\leq i}\cup\{I_{i+1}\}} aq^{-\alpha
}\bigl(\II_q(t)-bt q^{-\alpha}\bigr) +t(c-D_i).
\end{equation}
Then, using Lemma~\ref{lem-wc-drift},
the proof of Theorem~\ref{thm-WC-C1} can easily be adapted to prove that
$n^{-\rho}|\mathcal{D}_{i+1}^{(n)}|\convd H_{i+1}(0)$,
where $H_{i+1}(0)$ is the hitting time of 0 of $(\SSS_t(i))_{t\geq0}$,
and where $a,b,c$ are given by $a=c_{ F}^{\alpha}/\expec[W]$,
$b=c_{ F}^{\alpha}$ and $c=\theta$.

Indeed, in more detail, we shall work \textit{conditionally}
on $\mathcal{D}_{\leq i}^{(n)}$.
The proof of Theorem~\ref{thm-WC-C1} reveals that the main contribution
to\vspace*{1pt} $(\SSS_t(i))_{t\geq0}$ and $(\SSS_t^{(n)}(i))_{t\geq0}$
arises from the vertices $q\in[K]$. Now, since
$(\indic{a\in\mathcal{D}_{\leq i}^{(n)}})_{a\in[K]}$ is
a sequence of \textit{discrete} random variables taking a finite number
of outcomes and that converge in distribution, we have that its
probability mass function converges pointwise.
By~\cite{Thor00}, (6.3) on page 16, this implies that we can \textit{couple}
$(\indic{a\in\mathcal{D}_{\leq i+1}^{(n)}})_{a\in[K]}$
to $(\indic{a\in\mathcal{D}_{\leq i+1}})_{a\in[K]}$ in such a way
that
%
%
\begin{equation} \prob\bigl(\bigl(\indic{a\in\mathcal{D}_{\leq
i+1}^{(n)}}\bigr)_{a\in[K]} \neq\bigl(\indic{a\in\mathcal{D}_{
\leq i+1}}\bigr)_{a\in[K]}\bigr)=o(1).
\end{equation}
Therefore, \textit{whp}, there is a perfect coupling between the elements
of $\mathcal{D}_{\leq i+1}^{(n)}\cap[K]$ and $\mathcal
{D}_{\leq i+1}\cap[K]$.
When this is the case, we can basically think of the set of summands in
(\ref{Sli-rewr}) as being deterministic and follow the proof of
Theorem~\ref{thm-WC-C1} verbatim.

Further, the proof of Proposition~\ref{prop-wc-function} can be
adapted to prove
the joint convergence of
%
%
\begin{equation}\quad \bigl(n^{-\rho}\bigl|\mathcal{D}_{i+1}^{(n)}\bigr|,
\bigl(\indic{q\in\mathcal{D}_{i+1}^{(n)}}\bigr)_{q\geq1}\bigr) \convd
(H_{i+1}(0), (\II_q(H_{i+1}(0)))_{q\geq1}).
\end{equation}
Together with the induction hypothesis, this proves that~(\ref{IH-form})
also holds for all $j\leq i+1$, and, thus, we have advanced the
induction hypothesis.
This, in particular, proves Theorem~\ref{thm-WC-C[K]}. The proof for
cluster weights
follows in an identical way as the convergence proof of $n^{-\rho
}\weight(1)$
in the proof of Theorem~\ref{thm-WC-C1}.
\end{pf*}


\section{\texorpdfstring{Proofs of Theorems \protect\ref{thm-WC-3,4}, \protect\ref{thm-prop-3,4} and \protect\ref{prop-max-clusters-high-weight}}
{Proofs of Theorems 1.1, 1.5 and 1.6}}
\label{sec-pf-prop-max-clusters-high-weight}
In this section, we prove Theorems~\ref{thm-WC-3,4},~\ref{thm-prop-3,4} and
\ref{prop-max-clusters-high-weight} using the results in
Theorems~\ref{thm-WC-C1} and~\ref{thm-WC-C[K]}, as well as
Proposition~\ref{prop-wc-function}.
We start with a proof of Theorem \ref
{prop-max-clusters-high-weight}, followed by those of Theorems~\ref{thm-prop-3,4}
and~\ref{thm-WC-3,4}. Note that, combining parts (a) and (b)
in Theorem~\ref{prop-max-clusters-high-weight},
we obtain that, with high probability as $K$ becomes large,
the largest $m$ clusters are all
among the first $(|\Ccal_{\leq}(i)|)_{i\in[K]}$. This explains
why we start
the cluster exploration from the vertices with the highest weights.

\begin{pf*}{Proof of Theorem~\ref{prop-max-clusters-high-weight}}
(a) For $\max_{i\geq K} |\Ccal_{\leq}(i)|\geq\vep n^{\rho}$
to occur,
we must have that there exists a cluster using the vertices in
$[n]\setminus[K]$ such that (1) $|\Ccal_{\leq}(i)|\geq\vep
n^{\rho}$, and (2) the cluster $\Ccal_{\leq}(i)$ is not
connected to any of
the vertices in $[K]$.

By construction, the graph restricted to the vertices in
$[n]\setminus[K]$ is again a Norros--Reittu model, with edge
probabilities $p_{ij}=1-{\mathrm e}^{-w_iw_j/\ell_n}$, for all
$i,j\in[n]\setminus[K]$. However, no vertex in $[n]\setminus [K]$ found
to be in the cluster $\Ccal(i)$ is allowed to have an edge to any of
the vertices in $[K]$. We shall now bound this probability, making use
of the results in~\cite{Hofs09a}.

With
%
%
\begin{equation} Z_{\geq k}^{[K]} =\sum_{v=1}^n \indic
{|\cluster(v)|\geq k, \cluster(v)\cap[K]=\varnothing},
\end{equation}
we have
%
%
\begin{eqnarray} \prob\Bigl(\max_{i\geq K} |\Ccal_{\leq
}(i)|\geq k\Bigr) &=& \prob\bigl(Z_{\geq k}^{[K]}\geq k\bigr) \leq\frac
{\expec[Z_{\geq k}^{[K]}]}{k} \nonumber\\[-8pt]\\[-8pt]
&=&\frac{1}{k}\sum_{v=K+1}^n
\prob\bigl(|\cluster(v)|\geq k,
\cluster(v)\cap[K]=\varnothing\bigr).\nonumber
\end{eqnarray}
Denote by $\cluster^{[K]}(v)$ the cluster of $v$ restricted to
the vertices
$[n]\setminus[K]$. Then, due to the independence of disjoint sets of edges,
and the fact that $\cluster(v)\cap[K]=\varnothing$ only depends
on edges between $[K]$ and $[n]\setminus[K]$, while $|\cluster^{
[K]}(v)|\geq k$
depends only on edges between pairs of vertices in $[n]\setminus[K]$,
we obtain
%
%
\begin{eqnarray}
&&
\prob\bigl(|\cluster(v)|\geq k, \cluster(v)\cap
[K]=\varnothing\bigr) \nonumber\\
&&\qquad= \expec\bigl[\prob\bigl(\cluster(v)\cap
[K]=\varnothing\mid\cluster^{[K]}(v)\bigr) \indic{|\cluster^{
[K]}(v)|\geq k}\bigr]\\
&&\qquad= \expec\bigl[{\mathrm e}^{-\mathcal
{W}_{[K]} \mathcal{W}^{[K]}(v)/\ell_n} \indic{|\cluster
^{[K]}(v)|\geq k}\bigr],\nonumber
\end{eqnarray}
where, similarly to~(\ref{weight-cluster-def}), we define
%
%
\begin{equation} \mathcal{W}^{[K]}(v)=\sum_{a\in\cluster^{
[K]}(v)} w_a \quad\mbox{and}\quad \mathcal{W}_{[K]}=\sum
_{j=1}^K w_j.
\end{equation}
We split depending on whether $\mathcal{W}^{[K]}(v)\geq k/2$ or
not, to obtain
%
%
\begin{eqnarray}\quad
\label{first-term}
\prob\Bigl(\max_{i\geq K} |\Ccal_{\leq
}(i)|\geq k\Bigr) &\leq& \frac{1}{k}\sum_{v=K+1}^n{\mathrm
e}^{-\mathcal{W}_{[K]} k/(2\ell_n)} \prob\bigl(\bigl|\cluster^{
[K]}(v)\bigr|\geq k\bigr)\\
\label{sec-term}
&&{} +\frac{1}{k}\sum
_{v=K+1}^n\prob\bigl(\bigl|\cluster^{[K]}(v)\bigr|\geq k, \mathcal{W}^{
[K]}(v)\leq k/2\bigr).
\end{eqnarray}
For the first term we compute that, for some $C>0$,
%
%
\begin{equation} \mathcal{W}_{[K]}\geq c_{ F}\sum_{j=1}^K
(n/j)^{\alpha}\bigl(1+o(1)\bigr) \geq C n^{\alpha} K^{\rho}.
\end{equation}
Thus, when $k=k_n=\vep n^{\rho}$, we obtain,
for some $u>0$, and using $\alpha+\rho=1$ [see~(\ref{defs-al-rho-eta})],
%
%
\begin{eqnarray}
&&
\frac{1}{k_n}{\mathrm e}^{-\mathcal{W}_{[K]}
k_n/(2\ell_n)} \sum_{v\in[n]}\prob\bigl(\bigl|\cluster^{[K]}(v)\bigr|\geq
k_n\bigr) \nonumber\\
&&\qquad\leq \frac{1}{k_n}{\mathrm e}^{-u\vep K^{\rho}} \sum_{v\in
[n]}\prob\bigl(\bigl|\cluster^{[K]}(v)\bigr|\geq
k_n\bigr)\nonumber\\[-8pt]\\[-8pt]
&&\qquad\leq{\mathrm e}^{-u\vep
K^{\rho}} \frac{1}{k_n}\sum_{v\in[n]}\prob\bigl(|\cluster(v)|\geq
k_n\bigr)\nonumber\\
&&\qquad= {\mathrm e}^{-u\vep K^{\rho}} \frac{n}{k_n} \prob
\bigl(|\cluster(V)|\geq k_n\bigr),\nonumber
\end{eqnarray}
where $V\in[n]$ is a vertex chosen uniformly at random from
$[n]$. By~\cite{Hofs09a}, Proposition 2.4(a), there exists a
constant $a_1<\infty$ such that
%
%
\begin{eqnarray}
\prob\bigl(|\cluster(V)|\geq k_n\bigr)&\leq& a_1
\bigl(k_n^{-1/(\tau-2)}+\bigl(\vep_n\vee n^{-(\tau-3)/(\tau-1)}
\bigr)^{1/(\tau-3)}\bigr) \nonumber\\[-8pt]\\[-8pt]
&\leq& a_1\bigl(k_n^{-1/(\tau-2)}+n^{-1/(\tau
-1)}\bigr),\nonumber
\end{eqnarray}
so that, for $k=k_n=\vep n^{\rho}$ with $\vep<1$ and with $a_1'=2a_1$,
%
%
\begin{equation}
\frac{n}{k_n}\prob\bigl(|\cluster(V)|\geq k_n\bigr)\leq
a_1'\vep^{-(\tau-1)/(\tau-2)} n^{-\rho}.
\end{equation}
Therefore, the term in~(\ref{first-term}) is bounded by
%
%
\begin{equation} {\mathrm e}^{-a \vep K^{\rho}} a_1'\vep^{-(\tau
-1)/(\tau-2)}.
\end{equation}
When we pick $K=K(\vep)$ sufficiently large, we can make this as small
as we wish.

We continue with the term in~(\ref{sec-term}), for which we use a
large deviation argument.
We formulate this result in the following lemma:
\begin{Lemma}[(Large deviations for cluster weights)]
\label{lem-LD-CW}
For every $k=o(n)$ and $K=o(n)$, there exists a $J>0$ such that
%
%
\begin{equation} \prob\bigl(\exists v\dvtx\bigl|\Ccal^{[K]}(v)\bigr|\geq
k, \weight^{[K]}(v)\leq k/2\bigr)\leq n{\mathrm e}^{-J k}.
\end{equation}
\end{Lemma}
\begin{pf}
When\vspace*{1pt} $|\cluster^{[K]}(v)|\geq k$, then $\weight^{[K]}(v)$
is stochastically
bounded from below by the sum $\sum_{i=1}^k w_{v(i)}$, where
$(v(i))_{i=1}^k$ is
the sized-biased ordering of $[n]$, that is,
for every $j\notin(v(s))_{s\in[i-1]}$,
%
%
\begin{equation} \label{size-biased-reor} \prob\bigl(v(i)=j\mid
(v(s))_{s\in[i-1]}\bigr)=\frac{w_j}{\sum_{l\notin(v(s))_{s\in[i-1]}}
w_l}.
\end{equation}
See~\cite{BhaHofLee09a}, Section 2, Lemma 2.1, for more details about the
size-biased reordering. Indeed, each time we draw a random mark and,
conditionally on
this mark not being one that has been found earlier as well as on all
the marks found so far, it will be equal to $j$ with
the probability in~(\ref{size-biased-reor}). When $|\cluster^{
[K]}(v)|\geq k$,
we must draw a vertex that we have not seen yet, a total of at least
$k$ times.

We apply the size-biased reordering to the vertex set $[n]\setminus[K]$.
Then, for each $i$ and conditionally on $(v(s))_{s\in[i-1]}$,
the random variable $w_{v(i)}$ is stochastically bounded
from above by the random variable $W_i'$ with distribution
%
%
\begin{equation} \prob(W_i'=w_j)=\frac{w_j}{\ell_n-\sum
_{s=1}^{i-1+K} w_s},\qquad j\in[n]\setminus[i-1+K],
\end{equation}
that is, we have removed the vertices with the largest $i-1+K$ weights.
As a result, the random variables $(W_i')_{i\geq1}$
are \textit{independent}. Now take $\kappa>0$ very small,
and note that, whenever $k-1+K\leq\kappa n$ and for
every $i\leq k$, $W_i'$
is stochastically bounded
from above by a random variable $W_i^{(n)}(\kappa)$ with distribution
%
%
\begin{equation} \prob\bigl(W_i^{(n)}(\kappa)=w_j\bigr)=\frac{w_j}{\ell
_n-\sum_{s=1}^{\kappa n} w_s},\qquad j\in[n]\setminus[\kappa n],
\end{equation}
where the random variables $(W_i^{(n)}(\kappa))_{i=1}^k$ are i.i.d.
Now take $\kappa>0$ so small that
%
%
\begin{equation} \expec\bigl[W_i^{(n)}(\kappa)\bigr] =\sum_{j=\kappa n}^n
\frac{w_j^2}{\ell_n-\sum_{s=1}^{\kappa n} w_s} \geq3/4.
\end{equation}
Then,
%
%
\begin{eqnarray}
&&\prob\bigl(\exists v\dvtx\bigl|\Ccal^{[K]}(v)\bigr|\geq
k, \weight^{[K]}(v)\leq k/2\bigr)\nonumber\\
&&\qquad\leq \sum_{v=K+1}^n\prob
\bigl(\bigl|\cluster^{[K]}(v)\bigr|\geq k, \weight^{[K]}(v)\leq k/2\bigr)\\
&&\qquad\leq n\prob\Biggl(\sum_{i=1}^k W_i^{(n)}(\kappa)\leq k/2\Biggr).\nonumber
\end{eqnarray}
Intuitively, since $\expec[W_i^{(n)}(\kappa)]\approx\nu
_n\approx\nu=1$ for $\kappa>0$ small,
the Chernoff bound proves that
$\prob(\sum_{i=1}^k W_i^{(n)}(\kappa)\leq k/2)$ is
exponentially small in $k$, so that
the term in~(\ref{sec-term}) is exponentially small. We now make this
intuition precise.

By the Chernoff bound, for each $\varthetas\geq0$,
and by the fact that $(W_i^{(n)}(\kappa))_{i\in[k]}$ are i.i.d.
random variables,
we have
%
%
\begin{equation} \label{chernoff-bd} \prob\Biggl(\sum_{i=1}^k
W_i^{(n)}(\kappa)\leq k/2\Biggr) \leq{\mathrm e}^{\varthetas k/2}
\expec\bigl[{\mathrm e}^{-\varthetas\sum_{i=1}^k W_i^{
(n)}(\kappa)}\bigr] = ({\mathrm e}^{\varthetas/2}\phi_{n,\kappa
}(\varthetas))^{k},\hspace*{-35pt}
\end{equation}
where
%
%
\begin{equation} \phi_{n,\kappa}(\varthetas)=\expec\bigl[{\mathrm
e}^{-\varthetas W_1^{(n)}(\kappa)}\bigr]
\end{equation}
denotes the Laplace transform of $W_1(\kappa)$. By~(\ref{chernoff-bd}),
it suffices to prove that there exists a $\varthetas>0$ such that,
uniformly in $n$ sufficiently large, $\varthetas/2+\log\phi
_{n,\kappa}(\varthetas)<0$.
This is what we shall show now. By dominated convergence,
for each fixed $\varthetas>0$,
%
%
\begin{equation} \log\phi_{n,\kappa}(\varthetas)\rightarrow\log
\phi_{\kappa}(\varthetas) =\log\expec\bigl[{\mathrm e}^{-\varthetas
W(\kappa)}\bigr],
\end{equation}
where
%
%
\begin{equation} \prob\bigl(W(\kappa)\leq x\bigr)=\expec\bigl[[1-F]^{-1}(U)\mid
U\geq\kappa\bigr],
\end{equation}
and $U$ is a uniform random variable on $[0,1]$. As a result, the
distribution of $U$ \textit{conditionally on $U\geq\kappa$}
is uniform on $[\kappa,1]$. Let $U_\kappa$ denote a uniform random
variable on $[\kappa,1]$,
so that $W(\kappa)\stackrel{d}{=} [1-F]^{-1}(U_\kappa)$. Then,
$W(\kappa)$ has mean
$\expec[W(\kappa)]\geq3/4$ and bounded variance $\sigma_\kappa^2$
(since $W(\kappa)\leq[1-F]^{-1}(\kappa)<\infty$ a.s.).
Therefore, a~Taylor expansion yields that, for \textit{fixed} $\kappa>0$,
%
%
\begin{equation} \log\phi_{\kappa}(\varthetas) \leq-3\varthetas
/4+\sigma_\kappa^2 \varthetas^2+o(\varthetas^2).
\end{equation}
Now, fix a $\varthetas>0$ so small that
%
%
\begin{equation} \varthetas/2-3\varthetas/4+\sigma_\kappa
^2\varthetas^2\leq-\varthetas/6,
\end{equation}
and then $N$ so large that, for all $n\geq N$,
%
%
\begin{equation} \log\phi_{n,\kappa}(\varthetas)\leq\log\phi
_{\kappa}(\varthetas)+\varthetas/12.
\end{equation}
Then, indeed, for $n\geq N$, since $\varthetas>0$,
%
%
\begin{equation} \varthetas/2+\log\phi_{n,\kappa}(\varthetas)\leq
-\varthetas/6+\varthetas/12 =-\varthetas/12<0,
\end{equation}
so that
%
%
\begin{equation} {\mathrm e}^{\varthetas/2}\phi_{n,\kappa
}(\varthetas)\leq{\mathrm e}^{-\varthetas/12},
\end{equation}
which, in turn, implies that
%
%
\begin{equation} \sum_{v=1}^n\prob\bigl(\bigl|\cluster^{[K]}(v)\bigr|\geq k,
\mathcal{W}^{[K]}(v)\leq k/2\bigr) \leq n {\mathrm e}^{-k\varthetas
/12}.
\end{equation}
When $n\rightarrow\infty$, this proves the claim for $J=\varthetas/12$.
\end{pf}

To prove Theorem~\ref{prop-max-clusters-high-weight}(a), we apply
Lemma~\ref{lem-LD-CW} to the term in~(\ref{sec-term}), which
is then bounded by ${\mathrm e}^{-\Theta(\vep n^{\rho})}$ when we
take $k=\vep n^{\rho}$.

(b) We denote by
%
%
\begin{equation} \label{Zgeqk-def} Z_{\geq k}=\sum_{v=1}^n
\indic{|\cluster(v)|\geq k}
\end{equation}
the number of vertices that are contained in connected components
of size at least~$k$. In~\cite{Hofs09a}, the random variable $Z_{
\geq k}$ has been used
in a crucial way to prove probabilistic bounds on $|\Cmax|$. We now
slightly extend these results.

We shall prove that, for all $\vep>0$ sufficiently small,
there exist constants $b_2,C$ such that
%
%
\begin{equation} \label{aim-pop-b} \prob\bigl(Z_{\geq\vep
n^{\rho}}\leq b_2 n^{\rho}\vep^{-1/(\tau-2)}\bigr) \leq C \vep
^{2/(\tau-2)}.
\end{equation}
We first note that it suffices to prove~(\ref{aim-pop-b}) when $\nu
_n\leq1-Kn^{-\eta}$.
Indeed, the random variable $Z_{\geq\vep n^{\rho}}$
is increasing in the edge occupation statuses, and, therefore, we may
take $\lambda<0$ so that $-\lambda>K$ to achieve the claim.

We shall use a second moment method. By~\cite{Hofs09a}, Proposition 2.4(b),
there exists $a_2=a_2(K)$ such that
%
%
\begin{equation} \expec[Z_{\geq\vep n^{\rho}}] \geq n \prob
\bigl(|\cluster(V)|\geq\vep n^{\rho}\bigr)\geq a_2 n^{\rho}\vep^{-1/(\tau
-2)},
\end{equation}
where $V$ is chosen uniformly from $[n]$. Therefore,
when we take $b_2=a_2/2$,
%
%
\begin{equation} \prob\bigl(Z_{\geq\vep n^{\rho}}\leq b_2
n^{\rho}\vep^{-1/(\tau-2)}\bigr) \leq\prob(Z_{\geq
\vep n^{\rho}}\leq\expec[Z_{\geq\vep n^{\rho}}]/2).
\end{equation}
We take $\vep>0$ small, and bound, by the Chebychev inequality,
%
%
\begin{equation} \prob(Z_{\geq\vep n^{\rho}}\leq\expec
[Z_{\geq\vep n^{\rho}}]/2) \leq\frac{4\Var(Z_{\geq
\vep n^{\rho}})} {\expec[Z_{\geq\vep n^{\rho}}]^2}.
\end{equation}
By~\cite{Hofs09a}, Proposition 2.2, and~\cite{Hofs09a}, Proposition 2.5 and
(2.22),
uniformly in $k\geq1$,
%
%
\begin{equation} \Var(Z_{\geq k})\leq n\expec[|\cluster(V)|]
\leq n^{1+\eta}=n^{\rho}.
\end{equation}
As a result, we obtain
%
%
\begin{equation} \prob(Z_{\geq\vep n^{\rho}}\leq\expec
[Z_{\geq\vep n^{\rho}}]/2) \leq\frac{4 n^{2\rho}}
{a_2^2\vep^{-2/(\tau-2)}n^{2\rho}} =C \vep^{2/(\tau-2)},
\end{equation}
which is small when $\vep>0$ is small. We conclude that, with probability
at least $1-o_{\vep}(1)$, where $o_\vep(1)$ denotes a function that
is $o(1)$
uniformly in $n$ as $\vep\downarrow0$,
%
%
\begin{equation} Z_{\geq\vep n^{\rho}} \geq\expec[Z_{
\geq\vep n^{\rho}}]/2 \geq\frac{a_2}{2}\vep^{-1/(\tau-2)} n^{\rho
}.
\end{equation}
Since, by~\cite{Hofs09a}, Theorem 1.2, $|\Cmax|\leq\vep^{-1/2}
n^{\rho}$ with probability
at least $1-o_{\vep}(1)$,
there are, again with probability
at least $1-o_{\vep}(1)$, at least
%
%
\begin{equation} \frac{a_2}{2}\vep^{-1/(\tau-2)} n^{\rho}/ (\vep
^{-1/2} n^{\rho})=C \vep^{1/2-1/(\tau-2)}
\end{equation}
clusters of size at least $\vep n^{\rho}$. Since
$1/2-1/(\tau-2)<0$, the number of clusters of size
at least $\vep n^{\rho}$ tends to infinity
when $\vep\downarrow0$. By part (a),
\textit{whp} for $K\geq1$ large, these clusters
will be part of $(|\Ccal_{\leq}(i)|)_{i\in[K]}$
when $K=K(\vep)\geq1$ is sufficiently large.
\end{pf*}

We now complete the proof of Theorem~\ref{thm-prop-3,4}.
\begin{pf*}{Proof of Theorem~\ref{thm-prop-3,4}}
We use Proposition~\ref{prop-wc-function} and note that the limiting
variables are all nontrivial (i.e., they are equal to
0 or 1 each with positive probability). This proves~(\ref{lim-qij}).
The proof of
(\ref{lim-qi}) is similar, noting that $|\cluster_{\leq}(i)|$
equals $|\Cmax|$ with strictly positive probability.
\end{pf*}

We finally use Theorem~\ref{prop-max-clusters-high-weight} to complete the proof of Theorem~\ref{thm-WC-3,4}:
%
\begin{pf*}{Proof of Theorem~\ref{thm-WC-3,4}}
Weak convergence of $(|\cluster_{(i)}|n^{-\rho})_{i\geq1}$
in the product topology is equivalent to the weak\vadjust{\goodbreak} convergence of
$(|\cluster_{(i)}|n^{-\rho})_{i\in[m]}$ for any $m\geq1$;
see~\cite{Kall02}, Theorem 4.29. In turn, by Theorem \ref
{prop-max-clusters-high-weight},
this follows from the convergence in distribution of
$(|\cluster_{\leq}(i)|n^{-\rho})_{i\in[m]}$
for all $m$. The latter follows from Theorem~\ref{thm-WC-C[K]}.
Since, \textit{whp} for large $K$, again by Theorem \ref
{prop-max-clusters-high-weight},
$(|\cluster_{(i)}|n^{-\rho})_{i\in[m]}$ is equal to
the largest $m$ components of $(|\cluster_{\leq}(i)|n^{-\rho
})_{i\in[K]}$,
we have identified
%
%
\begin{equation} \label{gamm-ident} (\gamma_i(\lambda))_{i\geq
1}\stackrel{d}{=} \bigl(H_{(i)}(0)\bigr)_{i\geq1},
\end{equation}
where $(H_{(i)}(0))_{i\geq1}$ is $(H_i(0))_{i\geq1}$ ordered in
size. This completes the proof
of Theorem~\ref{thm-WC-3,4} and identifies the limiting random variables.
\end{pf*}

\section{\texorpdfstring{Proof of Theorem \protect\ref{thm-sub-crit}}{Proof of Theorem 1.3}}
\label{sec-sub-crit}
In this section, we shall prove Theorem~\ref{thm-sub-crit} on the largest
subcritical clusters. We shall extend the result also to the ordered weights
of subcritical clusters as formulated in Theorem \ref
{thm-scal-lim-cluster-weights},
which shall be a crucial ingredient in the proof of
Theorem~\ref{thm-mult-coal}, which is given in Section~\ref{secmult-coal}
below.

We shall prove that Theorem~\ref{thm-sub-crit} holds for $\weight
_{(j)}$ as
well as for $|\Ccal_{(j)}|$. Indeed, it shall also follow from
the result that
\textit{whp},
$\weight_{(j)}=\sum_{i\in\Ccal_{(j)}}w_{i}$, that is, the
$j$th largest cluster weight is the
weight of the $j$th largest cluster, as claimed in
Theorem~\ref{thm-scal-lim-cluster-weights}.

To prove this scaling, we shall prove that, when the weights are
equal to $\bfwit(\lambda_n)$ as defined in~(\ref{w-lambda-def}),
and when $\lambda_n\rightarrow-\infty$,
%
%
\begin{equation} \label{aim-b1} |\lambda_n| n^{-\rho}|\Ccal
(j)|\convp c_j,\qquad |\lambda_n| n^{-\rho} \weight(j)\convp c_j,
\end{equation}
where we recall that
%
%
\begin{equation} c_j=c_{ F}^{\alpha} j^{-\alpha}=\lim
_{n\rightarrow\infty} n^{-\alpha} w_j.
\end{equation}
Since $j\mapsto c_j$ is strictly decreasing, this means that, \textit{whp},
$\Ccal(j)=\Ccal_{\leq}(j)$. Thus, this also implies that
\textit{whp}, $\Ccal(j)=\Ccal_{(j)}$ for all $j\leq m$.
Then~(\ref{aim-b1}) proves the result for the ordered
cluster sizes and weights.


Recall the definitions of $T$, $T(i)$ and their weights $w_T$ and $w_{T(i)}$
introduced in Section~\ref{sec-BPs}, where also their moments
are computed in Lemma~\ref{lem-mom-cluster-weight}. We make frequent
use of
these computations. The proof of Theorem~\ref{thm-sub-crit} consists
of four key steps,
which we shall prove one by one.

\subsection*{Asymptotics of mean cluster size and weight of
high-weight vertices}
In the following lemma we investigate the means of $|\Ccal(j)|$ and
$\weight(j)$:
\begin{Lemma}[(Mean cluster size and weights)]
\label{lem-mean-cluster-size-weight}
As $n\rightarrow\infty$, for every $j\in{\mathbb N}$ fixed,
and when $\lambda_n\rightarrow-\infty$ such that $\nu_n(\lambda
_n)\rightarrow1$,
%
%
\begin{eqnarray} \expec[|\Ccal(j)|]&=&\frac{w_j}{1-\nu_n(\lambda
_n)}\bigl(1+o(1)\bigr),\nonumber\\[-8pt]\\[-8pt]
\expec[\weight(j)]&=&\frac{w_j}{1-\nu_n(\lambda
_n)}\bigl(1+o(1)\bigr).\nonumber
\end{eqnarray}
\end{Lemma}
\begin{pf}
By the fact that $|\Ccal(j)|$ and $T(j)$ can be coupled so that
$|\Ccal(j)|\leq T(j)$ a.s., we obtain that
%
%
\begin{equation} \expec[|\Ccal(j)|]\leq\expec[T(j)]=\frac
{w_j}{1-\nu_n(\lambda_n)},
\end{equation}
the latter equality following from Lemma~\ref{lem-mom-cluster-weight}(c).
A similar upper bound follows for $\expec[\weight(j)]$ now using
Lemma~\ref{lem-mom-cluster-weight}(d).

For the lower bound, we rewrite
%
%
\begin{equation} \expec[|\Ccal(j)|]=\expec[T(j)]-\expec[T(j)-|\Ccal
(j)|].
\end{equation}
Now, for $a_n=n^{\rho}\gg\expec[T(j)]$, we bound
%
%
\begin{equation} \label{diff-means-split} \qquad\expec[T(j)-|\Ccal(j)|]
\leq\expec\bigl[T(j)\indic{T(j)>a_n}\bigr]+\expec\bigl[[T(j)-|\Ccal(j)|]\indic
{T(j)\leq a_n}\bigr].
\end{equation}
By Lemma~\ref{lem-mom-cluster-weight}(c), the first term in (\ref
{diff-means-split}) is bounded by
%
%
\begin{eqnarray} \label{Tj-large}
&&\expec\bigl[T(j)\indic
{T(j)>a_n}\bigr]\nonumber\\
&&\qquad\leq\frac{1}{a_n}\expec[T(j)^2]\nonumber\\[-8pt]\\[-8pt]
&&\qquad=\frac{1}{a_n}
\biggl(\biggl(1+\frac{w_j}{1-\nu_n(\lambda_n)}\biggr)^2+ \frac{w_j(1+\nu_n(\lambda
_n))}{(1-\nu_n(\lambda_n))^2}\nonumber\\
&&\qquad\quad\hspace*{69pt}{}+ \frac{w_j}{(1-\nu_n(\lambda
_n))^3}\frac{1}{\ell_n}\sum_{i\in[n]} w_i^3\biggr).\nonumber
\end{eqnarray}
The first two terms in~(\ref{Tj-large}) are $o(w_j/(1-\nu_n(\lambda_n)))$
since $\nu_n(\lambda_n)=1+n^{-\eta}\lambda_n+o(n^{-\eta}|\lambda
_n|)$ by~(\ref{tildenun-asymp}) and the fact that $\lambda
_n\rightarrow-\infty$,
so that
%
%
\begin{equation} w_j/\bigl(1-\nu_n(\lambda_n)\bigr)\leq cn^{\alpha+\eta
}|\lambda_n|^{-1}=cn^{\rho}|\lambda_n|^{-1}=o(n^{\rho})=o(a_n),
\end{equation}
since $\alpha+\eta=\rho$ [recall~(\ref{defs-al-rho-eta})].
The last term in~(\ref{Tj-large}) is bounded by
%
%
\begin{equation} \frac{w_j}{1-\nu_n(\lambda_n)} \frac{c n^{3\alpha
-1}}{a_n (1-\nu_n(\lambda_n))^2} =\frac{w_j}{1-\nu_n(\lambda_n)}
c n^{3\alpha-1-\rho-2\eta}|\lambda_n|^{-2}.
\end{equation}
By~(\ref{defs-al-rho-eta}), $3\alpha-1-\rho-2\eta=3(\tau-4)/(\tau
-1)<0$, so that
also this term is $o(w_j/(1-\nu_n(\lambda_n)))$.

For the second term in~(\ref{diff-means-split}), we note that
differences between $T(j)$ and $|\Ccal(j)|$ arise due to
vertices which have been used \textit{at least twice} in $T(j)$.
Indeed, as explained in more detail in Section~\ref{sec-BPs},
the law of $|\Ccal(j)|$ can be obtained from the branching process
by removing vertices (and their complete offspring) of which the mark
has already been used (see the description of the cluster exploration in
Section~\ref{sec-cluster-BP} and the relation to branching processes
described in Sections~\ref{sec-cluster-BP} and~\ref{sec-BPs}). Thus,
when we draw vertex $i$ twice, then the second time we must thin the
entire tree
that is rooted at this vertex with mark $i$. The expected number of\vadjust{\goodbreak}
vertices in the tree
equals $\expec[T(i)]$, so that we arrive at
%
%
\begin{eqnarray}
&&\expec\bigl[[T(j)-|\Ccal(j)|]\indic{T(j)\leq a_n}\bigr]\nonumber\\[-2pt]
&&\qquad\leq
\sum_{i\in[n]} \expec\bigl[[T(j)-|\Ccal(j)|]\indic{T(j)\leq a_n}
\indic{\mathrm{mark}\ i\ \mathrm{drawn}\ \mathrm{at}\ \mathrm{least}\ \mathrm{twice}}
\bigr]\\[-2pt]
&&\qquad\leq \sum
_{i\in[n]} \expec[T(i)] \sum_{s_1<s_2=1}^{a_n} \prob(\mbox
{mark $i$ drawn at times } s_1, s_2).\nonumber
\end{eqnarray}
Now, $i$ can only be chosen at time $s_1$ when $T(j)\geq s_1-1$,
which is independent of the event that the mark $i$ is chosen at times
$s_1, s_2$.
Therefore,
%
%
\begin{eqnarray} \label{conc-term-1}\quad \expec\bigl[[T(j)-|\Ccal(j)|]\indic
{T(j)\leq a_n}\bigr]&\leq&\sum_{i\in[n]} \expec[T(i)] \sum
_{s_1<s_2=1}^{a_n} \prob\bigl(T(j)\geq s_1-1\bigr) \frac{w_i^2}{\ell_n^2}\nonumber
\\[-2pt]
&\leq& a_n \sum_{s_1=1}^{a_n} \prob\bigl(T(j)\geq s_1-1\bigr)\sum_{i\in[n]}
\expec[T(i)]\frac{w_i^2}{\ell_n^2}\\[-2pt]
&\leq& a_n \expec[T(j)]\sum
_{i\in[n]} \expec[T(i)]\frac{w_i^2}{\ell_n^2}.\nonumber
\end{eqnarray}
This is $o(w_j/(1-\nu_n(\lambda_n)))$ when $\lambda_n\rightarrow
-\infty$, since
%
%
\begin{eqnarray} \label{conc-term-1b} a_n \sum_{i\in[n]} \expec
[T(i)]\frac{w_i^2}{\ell_n^2} &=& a_n \sum_{i\in[n]} \frac{w_i^3}{\ell
_n^2(1-\nu_n(\lambda_n))} \leq\frac{C}{|\lambda_n|} n^{\rho
-2+3\alpha+\eta} \nonumber\\[-9pt]\\[-9pt]
&=&\frac{C}{|\lambda_n|}=o(1).\nonumber
\end{eqnarray}
This completes the proof for $\expec[|\Ccal(j)|]$. The proof for
$w_{T(j)}$ is similar.
Indeed, we split
%
%
\begin{eqnarray} \label{weight-split-a}
\expec\bigl[w_{T(j)}-\weight(j)\bigr]
&\leq&\expec\bigl[w_{T(j)}\indic{T(j)>a_n}\bigr]\nonumber\\[-9pt]\\[-9pt]
&&{}+\expec\bigl[\bigl[w_{T(j)}-\weight
(j)\bigr]\indic{T(j)\leq a_n}\bigr].\nonumber
\end{eqnarray}
The first term is now bounded by
%
%
\begin{equation} \label{weight-split-b} \expec\bigl[w_{T(j)}\indic
{T(j)>a_n}\bigr]\leq\frac{1}{a_n}\expec\bigl[w_{T(j)}T(j)\bigr],
\end{equation}
which we can again bound using $\expec[w_{T(j)}T(j)]\leq\expec
[w_{T(j)}^2]+\expec[T(j)^2]$
together with Lemma~\ref{lem-mom-cluster-weight}(a) and (b).
Further,
%
%
\begin{eqnarray} \label{weight-split-c}
&&\expec\bigl[\bigl[w_{T(j)}-\weight
(j)\bigr]\indic{T(j)\leq a_n}\bigr] \nonumber\\[-2pt]
&&\qquad\leq \sum_{i\in[n]} \expec
\bigl[\bigl[w_{T(j)}-\weight(j)\bigr]
\indic{T(j)\leq a_n} \indic{\mathrm{mark}\ i\
\mathrm{drawn}\ \mathrm{at}\ \mathrm{least}\ \mathrm{twice}}\bigr]\nonumber\\[-8pt]\\[-8pt]
&&\qquad\leq\sum_{i\in[n]} \expec\bigl[w_{T(i)}\bigr]
\sum_{s_1,s_2=1}^{a_n} \prob(\mbox{mark $i$ drawn at times }
s_1, s_2)\nonumber\\
&&\qquad\leq a_n \expec[T(j)]\sum_{i\in[n]} \expec
\bigl[w_{T(i)}\bigr]\frac{w_i^2}{\ell_n^2} =a_n \frac{w_j\nu_n(\lambda
_n)}{(1-\nu_n(\lambda_n))^2}\sum_{i\in[n]}\frac{w_i^3}{\ell
_n^2},\nonumber
\end{eqnarray}
so that
%
%
\begin{equation} \label{weight-split-d} \expec[\weight(j)] \geq
\expec\bigl[w_{T(j)}\bigr] -\frac{1}{a_n}\expec\bigl[w_{T(j)}T(j)\bigr]-a_n \frac
{w_j}{(1-\nu_n(\lambda_n))^2}\sum_{i\in[n]}\frac{w_i^3}{\ell
_n^2}.\hspace*{-32pt}
\end{equation}
We bound $\expec[w_{T(j)}T(j)]\leq\expec[w_{T(j)}^2]+\expec[T(j)^2]$.
Now we can simply follow the argument for $\expec[|\cluster(j)|]$.
\end{pf}

\subsection*{Cluster size and weight of high weight vertices are concentrated}
We note that, by the stochastic domination and the fact that $\expec
[|\Ccal(j)|]=\frac{w_j}{1-\nu_n(\lambda_n)}(1+o(1))$,
we have
%
%
\begin{equation} \Var(|\Ccal(j)|) \leq\Var(T(j))+o(\expec
[T(j)]^2).
\end{equation}
By Lemma~\ref{lem-mom-cluster-weight}(a),
%
%
\begin{eqnarray} \Var(T(j))&=&\frac{w_j(1+\nu_n(\lambda
_n))}{1-\nu_n(\lambda_n)} +\frac{w_j}{(1-\nu_n(\lambda_n))^3}
\biggl(\frac{1}{\ell_n} \sum_{l\in[n]} w_l^3\biggr)\nonumber\\[-8pt]\\[-8pt]
&=&o\biggl(\frac
{w_j^2}{(1-\nu_n(\lambda_n))^2}\biggr),\nonumber
\end{eqnarray}
since $j$ is fixed and
%
%
\begin{equation} \biggl(\frac{1}{\ell_n} \sum_{l\in[n]} w_l^3
\biggr)\bigl(1-\nu_n(\lambda_n)\bigr)^{-1} =\frac{C}{|\lambda_n|}n^{3\alpha-1+\eta
}=o(n^{\alpha})=o(w_j).
\end{equation}
%
For $w_{T(j)}$ the argument is identical.
We conclude that, for $j$ fixed, $\Var(|\Ccal(j)|)=o(\expec
[|\Ccal(j)|]^2)$
and $\Var(\weight(j))=o(\expec[\weight(j)]^2)$,
so that
%
%
\begin{equation} \frac{|\Ccal(j)|}{\expec[|\Ccal(j)|]}\convp1,\qquad
\frac{\weight(j)}{\expec[\weight(j)]}\convp1,
\end{equation}
and then Lemma~\ref{lem-mean-cluster-size-weight} completes the proof of
(\ref{aim-b1}).

\subsection*{Cluster weight sums}
We start by proving a convenient result relating the cluster weights
$\weight(j)$ and
$\weight_{\leq}(j)$.
\begin{Lemma}[(Cluster weight properties)]
\label{lem-cluster-weight-prop}
\textup{(a)} For every integer $m\geq2$,
%
%
\begin{equation} \sum_{j\in[n]} \weight_{\leq}(j)^m =\sum
_{j\in[n]} w_j \weight(j)^{m-1}.
\end{equation}

\textup{(b)} For every $i,j\in[n]$,
%
%
\begin{equation} \expec\bigl[\weight(i)\weight(j)\indic{i\nc j}\bigr] \leq
\expec[\weight(i)]\expec[\weight(j)].
\end{equation}
\end{Lemma}
\begin{pf}
(a) We compute
%
%
\begin{eqnarray} \sum_{j\in[n]} \weight_{\leq}(j)^m &=& \sum
_{j\in[n]}\sum_{i_1,\ldots, i_{m}} \prod_{s=1}^m w_{i_s}\indic
{i_s\in\Ccal(i_1)\ \forall s=2, \ldots, m, \min\Ccal(i_1)=j}\nonumber\\
&=& \sum_{i_1,\ldots, i_{m}} \prod_{s=1}^m w_{i_s}\indic{i_s\in
\Ccal(i_1)\ \forall s=2, \ldots, m}\\
&=&\sum_{i_1\in[n]} w_{i_1}\weight
(i_1)^{m-1}.\nonumber
\end{eqnarray}

(b) We write out
%
%
\begin{eqnarray} \expec\bigl[\weight(i)\weight(j)\indic{i\nc j}\bigr] &=& \sum
_{k,l} w_{k}w_l \prob(i\conn k, j\conn l, i\nc j) \nonumber\\
&\leq&\sum_{k,l}
w_{k}w_l\prob(i\conn k)\prob(j\conn l)\\
&=& \expec[\weight(i)]\expec
[\weight(j)]\nonumber
\end{eqnarray}
by the BK-inequality; see~\cite{Grim99}, Section 2.3.
\end{pf}

\subsection*{Only high-weight vertices matter}
We start by proving that the probability that, for $K\geq1$, there
exists a
$j>K$ such that $\weight_{\leq}(j)\geq\vep n^{\rho}/|\lambda
_n|$ is small.
Since, for all $j\leq K$, we have that $|\lambda_n|n^{-\rho}\weight
(j)\convp c_j$,
we have that, for all $i\leq m$ and $m$ such that $c_m>\vep$, $\weight
(j)=\weight_{(j)}$.

Recall that $\weight^{[K]}(j)$ is the weight of the
cluster of $j$ in the random graph
only making use of the vertices in $[n]\setminus[K]$,
and let $\weight^{[K]}_{\leq}(j)=\weight^{[K]}(j)$
when $j$ is the minimal element in $\cluster^{[K]}(j)$.
If there exists a $j>K$ such that
$\weight^{[K]}(j)\geq\vep n^{\rho}/|\lambda_n|$, then
%
%
\begin{equation} \label{high-weights-late-conseq} \sum_{j>K} \weight
^{[K]}_{\leq}(j)^3 \geq\frac{\vep^3}{|\lambda_n|^3}
n^{3\rho}.
\end{equation}
Since
%
%
\begin{equation} \ell_n\geq\sum_{j>K} w_j,
\end{equation}
we see\vspace*{1pt} that this random graph is stochastically bounded by
the random graph having weights $\bfw^{[K]}$, where $w^{[K]}_j=0$ when
$j\leq K$ and $w^{[K]}_j=w_j$ otherwise. By the Markov inequality,
%
%
\begin{eqnarray} \label{Markov-weights}
&&\prob\bigl(\exists j>K\dvtx
\weight_{\leq}(j)\geq\vep n^{\rho}/|\lambda_n|\bigr)\nonumber\\[-8pt]\\[-8pt]
&&\qquad \leq
\sum_{j>K} \weight^{[K]}_{\leq}(j)^3 = \frac{|\lambda
_n|^3}{\vep^3} n^{-3\rho}\sum_{j>K} w_j \expec\bigl[\weight^{
[K]}(j)^2\bigr],\nonumber
\end{eqnarray}
where we have used Lemma~\ref{lem-cluster-weight-prop} for the
equality. We note that we can again stochastically dominate $|\cluster
^{[K]}(j)|$
by $T^{[K]}(j)$ and $\weight^{[K]}(j)$ by $w_{T^{[K]}(j)}$,
where now the offspring distribution is equal to $X^{[K]}_i=\XBP
_i\indic{M_i>K}$
(recall Section~\ref{sec-cluster-BP}).
Therefore, by Lemma~\ref{lem-mom-cluster-weight}(d), we obtain that
%
%
\begin{eqnarray} \label{sec-mom-weight}
\expec\bigl[\weight^{
[K]}(j)^2\bigr] &\leq&\expec\bigl[w_{T^{[K]}(j)}^2\bigr]\nonumber\\[-8pt]\\[-8pt]
&=&\biggl(\frac{w_j^{
[K]}}{1-\nu_n^{[K]}}\biggr)^2+ \frac{w_j^{[K]}}{(1-\nu
_n^{[K]})^3}\biggl(\frac{1}{\ell_n} \sum_{i\in[n]} \bigl(w_i^{
[K]}\bigr)^3\biggr),\nonumber
\end{eqnarray}
where
%
%
\begin{equation} w_j^{[K]}=w_j\indic{j>K},\qquad \nu^{
[K]}_n=\sum_{j\in[n]} \bigl(w_j^{[K]}\bigr)^2/\ell_n.
\end{equation}

It is not hard to see that, for each $K\geq1$ fixed, as $n\rightarrow
\infty$,
%
%
\begin{equation} \label{weight-tree-bd} \expec\bigl[w_{T^{
[K]}(j)}^2\bigr]\leq\biggl(\frac{w_j}{1-\nu_n(\lambda_n)}
\biggr)^{\!2}+\bigl(1+o(1)\bigr)\frac{w_j}{(1-\nu_n(\lambda_n))^3}
\Biggl(\frac{1}{\ell
_n} \sum_{i>K}^n w_i^3\Biggr).\hspace*{-35pt}
\end{equation}
Substitution of the bound~(\ref{sec-mom-weight})
in the right-hand side of~(\ref{Markov-weights})
and performing the sum over $j$ gives that
%
%
\begin{eqnarray} \label{aim-b2}
\sum_{j>K} w_j \expec\bigl[\weight^{
[K]}(j)^2\bigr] &\leq&\frac{1}{(1-\nu_n(\lambda_n))^2}\biggl(1+\frac
{1}{1-\nu_n(\lambda_n)}\biggr)\sum_{j>K} w_j^3\nonumber\\[-8pt]\\[-8pt]
&\leq& CK^{1-3\alpha}
(n^{\rho}/|\lambda_n|)^3,\nonumber
\end{eqnarray}
so that
%
%
\begin{eqnarray} \label{Markov-weights-rep}\quad
\prob\bigl(\exists j>K\dvtx
\weight_{\leq}(j)\geq\vep n^{\rho}/|\lambda_n|\bigr) &\leq&
|\lambda_n|^3\vep^{-3} n^{-3\rho}CK^{1-3\alpha} (n^{\rho
}/|\lambda_n|)^3 \nonumber\\[-8pt]\\[-8pt]
&=&CK^{1-3\alpha} \vep^{-3},\nonumber
\end{eqnarray}
which can be made arbitrarily small by taking $K=K(\vep)$ large.

We complete this section by proving that the probability
that there exists a $j>K$ such that $|\Ccal_{\leq}(j)|\geq\vep
n^{\rho}/|\lambda_n|$
is small. For this, we use Lemma~\ref{lem-LD-CW}, which proves that,
\textit{whp}, if $|\Ccal_{\leq}(j)|\geq\vep n^{\rho}/|\lambda_n|$,
then also $\weight_{\leq}(j)\geq\vep n^{\rho}/(2|\lambda_n|)$.
Thus, the result for cluster sizes follows from the proof for cluster weights.
This completes the proof of Theorem~\ref{thm-sub-crit}.

\section{\texorpdfstring{Proof of Theorem \protect\ref{thm-mult-coal}}{Proof of Theorem 1.2}}
\label{secmult-coal}
In this section, we prove Theorem~\ref{thm-mult-coal}.
We start by using~\cite{AldLim98}, Proposition 7, to show that
the random graph multiplicative coalescent converges
(recall Lemma~\ref{lem-mult-coal}).

\subsection*{Convergence of the random graph multiplicative coalescent
at fixed time}
We apply~\cite{AldLim98}, Proposition 7, which gives conditions
to show that, for fixed $\lambda\in{\mathbb R}$, the random sequence
$\mathbf{X}^{(n)}(|\lambda_n|+\lambda)$
converges in distribution to a random variable which has the
same distribution as the
$(0,\beta,\bfd)$-multiplicative coalescent at time $\lambda$ when
three conditions are satisfied about the
initial state $\bfitx^{(n)}
=\mathbf{X}^{(n)}(0)$.
To state these conditions, we define, for $r=2,3$,
with $\bfitx^{(n)}=(x_j^{(n)})_{j\geq1}$,
%
%
\begin{equation} \sigma_r\bigl(\bfitx^{(n)}\bigr)=\sum_j \bigl(x_j^{
(n)}\bigr)^r.
\end{equation}
Then, the conditions in~\cite{AldLim98}, Proposition 7, are that,
as $\lambda_n\rightarrow-\infty$:

\begin{longlist}[(a)]
\item[(a)]
%
%
\begin{equation} \label{cond-a} |\lambda_n|\bigl(|\lambda_n|\sigma
_2\bigl(\bfitx^{(n)}\bigr)-1\bigr)\convp-\beta;
\end{equation}
\item[(b)]
%
%
\begin{equation} \label{cond-b} \frac{x_j^{(n)}}{\sigma
_2(\bfitx^{(n)})}\convp d_j;
\end{equation}
\item[(c)]
%
%
\begin{equation} \label{cond-c} |\lambda_n|^3\sigma_3\bigl(\bfitx
^{(n)}\bigr)\convp\sum_{j=1}^\infty d_j^3.
\end{equation}
\end{longlist}
The conditions (a)--(c) above are not precisely what is in
\cite{AldLim98}, Proposition 7, and
we start by explaining how (a)--(c) imply the conditions for
\cite{AldLim98}, Proposition 7. Indeed, in
\cite{AldLim98}, Proposition 7, the condition in (a) is replaced by
$\sigma_2(\bfitx^{(n)})\rightarrow0$,
and the process
%
%
\begin{equation} \mathbf{X}^{(n)}\biggl(\frac{1}{\sigma
_2(\bfitx^{(n)})}+\lambda\biggr)
\end{equation}
is proved to converge to the realization
of a $(0,0,\bfd)$-multiplicative coalescent
at time $\lambda$. Under condition (a) (and the fact that $\lambda
_n\rightarrow-\infty$), (a) implies that $1/\sigma_2(\bfitx^{
(n)})=|\lambda_n|-\beta+o(1)$.
Since if $(\mathbf{X}(t))_{t}$ is a multiplicative coalescent
with parameters $(0,0,\bfd)$, then
$(\mathbf{X}(t-\beta))_{t}$ is a multiplicative coalescent
with parameters $(0,\beta,\bfd)$ (see~\cite{AldLim98}, (13)),
and using the continuity proved in~\cite{AldLim98}, Lemma 27,
this proves the fact that $\mathbf{X}^{(n)}(|\lambda
_n|+\lambda)$
converges in distribution to a random variable which has the
same distribution as a $(0,\beta,\bfd)$-multiplicative coalescent
at time $\lambda$.
Also, in~\cite{AldLim98}, Proposition 7, condition (c) is replaced by
the condition that
%
%
\begin{equation} \frac{\sigma_3(\bfitx^{(n)})}{\sigma_2(\bfitx
^{(n)})^3}\convp\sum_{j=1}^\infty d_j^3,
\end{equation}
which follows from a combination of (a) and (c). Further, in (a)--(c),
we work with \textit{convergence in probability} (as the initial state is
a random variable), while in~\cite{AldLim98}, Proposition 7, the initial
state is considered to be deterministic. This is a minor change.

In the remainder of this section, we shall show that conditions
(a)--(c) hold with
$\beta=-\zeta/\expec[W]$ and $d_j=c_j/\expec[W]$.


\subsection*{\texorpdfstring{Asymptotics of $\sigma_2(\bfitx^{(n)})$}
{Asymptotics of $\sigma_2(\bfitx^{(n)})$}}

In the following lemma, we state the properties of
$\sigma_2(\bfitx^{(n)})$ that we shall rely on.
In order to state the result, we recall that
%
%
\begin{equation} \sigma_2\bigl(\bfitx^{(n)}\bigr)=\sum_j \bigl(x_j^{
(n)}\bigr)^2,
\end{equation}
where $x_j^{(n)}=n^{-\rho} \weight_{(j)}$,
and where the vertex weights are now given by
%
%
\begin{equation} \bar w_j(0) =(1+\lambda_n\ell_nn^{-2\rho})w_j
=w_j(\lambda_n\ell_nn^{-2\rho+\eta})= w_j(\lambda_n\ell_n/n),
\end{equation}
since $2\rho-\eta=1$,
so that
%
%
\begin{equation} \label{bfw0-ident} \bar{\bfw}(0)=\bfw(\lambda
_n\ell_n/n)= \bfw(\expec[W]\lambda_n)\bigl(1+o(1)\bigr).
\end{equation}
Now,
%
%
\begin{equation} \sigma_2\bigl(\bfitx^{(n)}\bigr)=n^{-2\rho} \sum_{j\geq
1} \weight_{(j)}^2 =n^{-2\rho} \sum_{j\geq1} \weight_{
\leq}(j)^2
\end{equation}
and, thus, by Lemma~\ref{lem-cluster-weight-prop},
%
%
\begin{equation} \sigma_2\bigl(\bfitx^{(n)}\bigr)=n^{-2\rho}\sum_{j\in
[n]} \weight_{\leq}(j)^2 =n^{-2\rho}\sum_{j\in[n]} w_j
\weight(j).
\end{equation}
We continue by investigating the mean and variance of the above sum:
\begin{Lemma}[{[Mean and variance of $\sigma_2(\bfitx^{(n)})$]}]
\label{lem-mean-sigma2}
When the weights $\bfw(\lambda_n)$ satisfy that $\nu_n(\lambda
_n)<1-\lambda_n n^{-\eta}$, then:

\begin{longlist}
\item
%
%
\begin{equation} \expec\biggl[\sum_{i\in[n]} w_i \weight(i)
\biggr]=\frac{\ell_n}{1-\nu_n(\lambda_n)} +o(n^{2\rho}\lambda_n^{-2});
\end{equation}
\item
%
%
\begin{eqnarray} \label{Var-bd-clusters}\quad
&&\Var\biggl(\sum_{i\in
[n]} w_i \weight(i)\biggr) \nonumber\\[-8pt]\\[-8pt]
&&\qquad\leq\ell_n\expec[w_T^3]\leq C
\biggl(\expec[w_T]^4\frac{1}{\ell_n} \sum_{i\in[n]} w_i^4+\expec
[w_T]^2\expec[w_T^2]\frac{1}{\ell_n} \sum_{i\in[n]}
w_i^3\biggr).\nonumber
\end{eqnarray}
\end{longlist}
\end{Lemma}
\begin{pf}
(i) We bound
%
%
\begin{equation} \expec\biggl[\sum_{i\in[n]} w_i \weight(i)\biggr]
\leq\expec\biggl[\sum_{i\in[n]} w_i w_{ T(i)}\biggr]=\ell_n
\expec[w_{T}]=\frac{\ell_n \nu_n}{1-\nu_n}.
\end{equation}
For the lower bound, we make use of the bound alike in~(\ref{weight-split-d}),
%
%
\begin{eqnarray} \label{weight-split-c-rep}\qquad
\expec\bigl[w_{T(i)}-\weight
(i)\bigr] &\leq& \sum_{j\in[n]} \expec\bigl[\bigl[w_{T(i)}-\weight(i)\bigr]\indic
{\mathrm{mark}\ j\ \mathrm{drawn}\ \mathrm{at}\ \mathrm{least}\
\mathrm{twice}}\bigr]\nonumber\\[-8pt]\\[-8pt]
&\leq& \sum_{j\in[n]}
\expec\bigl[w_{T(j)}\bigr] \sum_{s_1<s_2} \prob(\mbox{mark $j$ drawn at
times } s_1, s_2)\nonumber.
\end{eqnarray}
Now, there are two contributions, depending on whether $s_2$ is in the
family tree of $s_1$ or not.
When it is not, then the events $\{\mbox{mark $j$ drawn at time } s_1\}
$ and $\{\mbox{mark $j$ drawn at time } s_2\}$
are completely independent, and we arrive at
%
%
\begin{eqnarray} \sum_{s_1<s_2} \prob(\mbox{mark $j$ drawn at
times } s_1, s_2) &=&\frac{w_j^2}{\ell_n^2} \sum_{s_1<s_2}\prob
\bigl(T(i)\geq s_1\bigr) \nonumber\\[-8pt]\\[-8pt]
&\leq&\frac{w_j^2}{\ell_n^2} \expec[T(i)^2].\nonumber
\end{eqnarray}
When $s_2$ is in the family tree of $s_1$, then we obtain the bound
%
%
\begin{equation}
\sum_{s_1<s_2} \prob(j\mbox{ chosen at times }
s_1, s_2) =\frac{w_j^2}{\ell_n^2} \sum_{s_1<s_2}\prob\bigl(T(i)\geq
s_1\bigr)\prob(s_2\in T_{s_1}),\hspace*{-35pt}
\end{equation}
where we denote the tree rooted at $s_1$ by $T_{s_1}$. Thus,
denoting by $|T_{s_1}|$ the number of elements in $T_{s_1}$,
%
%
\begin{equation} \sum_{s_2}\prob(s_2\in T_{s_1})\leq\expec
[|T_{s_1}|]=\expec[T(j)],
\end{equation}
and we arrive at a contribution of
%
%
\begin{eqnarray} \sum_{s_1<s_2} \prob(j\mbox{ chosen at times }
s_1, s_2) &\leq&\frac{w_j^2}{\ell_n^2}\sum_{s_1<s_2}\prob
\bigl(T(i)\geq s_1\bigr)\expec[T(j)] \nonumber\\[-8pt]\\[-8pt]
&=&\frac{w_j^2}{\ell_n^2}\expec
[T(i)]\expec[T(j)].\nonumber
\end{eqnarray}
Therefore,
%
%
\begin{eqnarray}
\expec\bigl[w_{T(i)}-\weight(i)\bigr]&\leq&\sum_{j\in[n]}
\expec\bigl[w_{T(j)}\bigr]\frac{w_j^2}{\ell_n^2} \bigl(\expec[T(i)^2]+\expec
[T(i)]\expec[T(j)]\bigr)\nonumber\\
&=&\expec[T(i)^2]\frac{1}{1-\nu_n(\lambda
_n)}\sum_{j\in[n]}\frac{w_j^3}{\ell_n^2} \\
&&{}+\expec[T(i)]\frac
{1}{(1-\nu_n(\lambda_n))^2}\sum_{j\in[n]} \frac{w_j^4}{\ell
_n^2}.\nonumber
\end{eqnarray}
Thus we obtain
\begin{eqnarray*}
\sum_{i\in[n]}w_i\expec
[\weight(i)] &\geq&\sum_{i\in[n]}w_i\expec\bigl[w_{T(i)}\bigr] -\sum_{i\in
[n]}w_i\expec[T(i)^2]\frac{1}{1-\nu_n(\lambda_n)} \sum_{j\in
[n]}\frac{w_j^3}{\ell_n^2} \\
&&{}-\sum_{i\in[n]}w_i^2\frac{1}{(1-\nu
_n(\lambda_n))^3}\sum_{j\in[n]} \frac{w_j^4}{\ell_n^2}.
\end{eqnarray*}
We bound
%
%
\begin{eqnarray}
&&\sum_{i\in[n]}w_i\expec[T(i)^2] \sum_{j\in
[n]}\frac{w_j^3}{\ell_n^2(1-\nu_n(\lambda_n))} \nonumber\\
&&\qquad\leq \frac{C}{\ell
_n^2(1-\nu_n(\lambda_n))^3} \biggl(\sum_{i\in[n]}w_i^3\biggr)^2
+C\sum_{i\in[n]}\frac{w_i^2}{(1-\nu_n(\lambda_n))^4} \sum_{j\in
[n]}\frac{w_j^3}{\ell_n^2}\\
&&\qquad\leq C |\lambda_n|^{-3} n^{3\eta
-2+6\alpha}+C |\lambda_n|^{-4} n^{4\eta-1+3\alpha}.\nonumber
\end{eqnarray}
Now, $3\eta-2+6\alpha=1<2\rho=(\tau-2)/(\tau-1)$, since $\tau>3$,
so that the first term is
$o(n^{2\rho}/|\lambda_n|^2)$. For the second term
$4\eta-1+3\alpha=2\rho+(\tau-4)/(\tau-1)<2\rho$, so this terms is also
$o(n^{2\rho}|\lambda_n|^{-2})$. Similarly,
%
%
\begin{eqnarray}
\sum_{i\in[n]}w_i^2\frac{1}{(1-\nu_n(\lambda
_n))^3}\sum_{j\in[n]} \frac{w_j^4}{\ell_n^2} &=&\nu_n(\lambda_n)
\frac{1}{(1-\nu_n(\lambda_n))^3}\sum_{j\in[n]} \frac{w_j^4}{\ell
_n} \nonumber\\[-8pt]\\[-8pt]
&\leq& C|\lambda_n|^{-3} n^{3\eta-1+4\alpha}.\nonumber
\end{eqnarray}
Again, $3\eta-1+4\alpha=2(\tau-3)/(\tau-1)<2\rho$, so also this
contribution is $o(n^{2\rho}|\lambda_n|^{-2})$.

(ii) We shall start by bounding the second moment. For this, we rewrite
%
%
\begin{equation} \expec\biggl[\biggl(\sum_{i\in[n]} w_i \weight
(i)\biggr)^2\biggr] =\sum_{i_1,i_2} w_{i_1} w_{i_2} \expec[\weight
(i_1)\weight(i_2)].
\end{equation}
Now we split
%
%
\begin{eqnarray}
\expec[\weight(i_1)\weight(i_2)] &=&\expec\bigl[\weight
(i_1)\weight(i_2)\indic{i_1\conn i_2}\bigr] \nonumber\\[-8pt]\\[-8pt]
&&{}+\expec\bigl[\weight(i_1)\weight
(i_2)\indic{i_1\nc i_2}\bigr].\nonumber
\end{eqnarray}
By Lemma~\ref{lem-cluster-weight-prop}(b), the second term is bounded
from above by
$\expec[\weight(i_1)]\expec[\weight(i_2)]$. Therefore, summing over
$i_1,i_2$, we obtain that
%
%
\begin{eqnarray}
\expec\biggl[\biggl(\sum_{i\in[n]} w_i \weight
(i)\biggr)^2\biggr] &\leq&\expec\biggl[\sum_{i\in[n]} w_i \weight
(i)\biggr]^2\nonumber\\[-8pt]\\[-8pt]
&&{}+ \sum_{i_1,i_2} w_{i_1} w_{i_2}\expec\bigl[\weight
(i_1)\weight(i_2)\indic{i_1\conn i_2}\bigr],\nonumber
\end{eqnarray}
so that
%
%
\begin{eqnarray} \Var\biggl(\sum_{i\in[n]} w_i\weight(i)\biggr)
&=&\expec\biggl[\biggl(\sum_{i\in[n]} w_i \weight(i)\biggr)^2\biggr]
-\expec\biggl[\sum_{i\in[n]} w_i \weight(i)\biggr]^2\nonumber\\
&\leq&\sum
_{i_1,i_2}w_{i_1}w_{i_2}\expec\bigl[\weight(i_1)\weight(i_2)\indic
{i_1\conn i_2}\bigr]\nonumber\\[-8pt]\\[-8pt]
&=&\sum_{i\in[n]}w_i\expec[\weight(i)^3]\nonumber\\
&\leq&
\sum_{i\in[n]}w_i\expec\bigl[w_{T(i)}^3\bigr]
=\ell_n \expec[w_{T}^3].\nonumber
\end{eqnarray}
The upper bound on $\expec[w_{T}^3]$ follows as in the proof of Lemma
\ref{lem-mom-cluster-weight}.
\end{pf}

\subsection*{Check of convergence conditions}
We conclude that we are left to prove that conditions (a), (b) and (c) in
(\ref{cond-a})--(\ref{cond-c}) hold.
We shall prove these conditions in the order (b), (c) and (a),
condition (a) being the
most difficult one.

Condition (b) follows from~(\ref{aim-b1}) and condition (a), as we
show now.
Substituting~(\ref{bfw0-ident}) into~(\ref{aim-b1}),
we obtain that
%
%
\begin{equation}\quad
x_j^{(0)}=n^{-\rho}\weight_{(j)} =\frac
{c_j}{\expec[W]|\lambda_n|}\bigl(1+o_{\prob}(1)\bigr)=\bigl(1+o_{\prob
}(1)\bigr)d_j/|\lambda_n|,
\end{equation}
where $d_j=c_j/\expec[W]$. Further, the first-order
asymptotics in condition (a) proves that
$|\lambda_n|\sigma_2(\bfitx^{(n)})\convp1$, so that
the factor $1/\sigma_2(\bfitx^{(n)})$ in condition (b) can be
replaced by
a multiplication by $|\lambda_n|$. We conclude that
$|\lambda_n|x_j^{(0)}\convp d_j$, where $d_j=c_j/\expec[W]$,
as required.

For condition (c), we apply similar ideas and
start with
%
%
\begin{equation} |\lambda_n|^3\sigma_3\bigl(\bfitx^{(n)}\bigr) =\sum
_{j\in[n]} (|\lambda_n| n^{-\rho} \weight_{\leq}(j)
)^3.
\end{equation}
The summands for $j>K$ can be bounded using Lemma \ref
{lem-cluster-weight-prop} by
%
%
\begin{equation} \sum_{j>K} (|\lambda_n| n^{-\rho} \weight
_{\leq}(j))^3 \leq(|\lambda_n| n^{-\rho}
)^{3}\sum_{j>K}w_j \weight^{[K]}(j)^2,
\end{equation}
which is small in probability by the Markov inequality and~(\ref{aim-b2}).
The summands for $j\leq K$ converge in probability by~(\ref{aim-b1}).
Thus condition (c) follows from~(\ref{aim-b1}) and~(\ref{aim-b2}).

We continue with condition (a), which is equivalent to the statement that
%
%
\begin{equation} \label{aim-a}
\sigma_2\bigl(\bfitx^{(n)}\bigr)=\frac
{1}{|\lambda_n|}-\frac{\beta} {\lambda_n^2}+o_{\prob
}(|\lambda_n|^{-2}),
\end{equation}
where $\beta=-\zeta/\expec[W]$.

We shall prove~(\ref{aim-a}) by a second moment method.
We first identify, by Lem\-ma~\ref{lem-cluster-weight-prop}(a),
%
%
\begin{equation} \sigma_2\bigl(\bfitx^{(n)}\bigr)=n^{-2\rho} \sum_{j\in
[n]} \weight_{\leq}(j)^2 =n^{-2\rho} \sum_{i\in[n]} w_i
\weight(i).
\end{equation}
Thus, in order to prove~(\ref{aim-a}), it suffices to show that
%
%
\begin{equation} \label{aim-a1} \expec\biggl[\sum_{i\in[n]} w_i \weight
(i)\biggr]=n^{2\rho} \biggl(|\lambda_n|^{-1}+\frac{\zeta}{\expec
[W]}|\lambda_n|^{-2} +o(|\lambda_n|^{-2})\biggr)
\end{equation}
and
%
%
\begin{equation} \label{aim-a2} \Var\biggl(\sum_{i\in[n]} w_i
\weight(i)\biggr)=o(n^{4\rho}|\lambda_n|^{-4}).
\end{equation}
Indeed, by~(\ref{aim-a1}), we have that, for $n$ sufficiently large,
%
%
\begin{eqnarray}
&&\prob\biggl(\biggl|\sigma_2\bigl(\bfitx^{
(n)}\bigr)-|\lambda_n|^{-1}-\frac{\zeta}{\expec[W]}|\lambda
_n|^{-2}\biggr|\geq\vep|\lambda_n|^{-2}\biggr) \nonumber\\[-8pt]\\[-8pt]
&&\qquad\leq\prob\bigl(
\bigl|\sigma_2\bigl(\bfitx^{(n)}\bigr)-\expec\bigl[\sigma_2\bigl(\bfitx^{
(n)}\bigr)\bigr]\bigr|\geq\vep|\lambda_n|^{-2}/2\bigr),\nonumber
\end{eqnarray}
which, by the Chebychev inequality is bounded by
%
%
\begin{eqnarray}\qquad
\prob\biggl(\biggl|\sigma_2\bigl(\bfitx^{
(n)}\bigr)-|\lambda_n|^{-1}+\frac{\zeta}{\expec[W]}|\lambda
_n|^{-2}\biggr|\geq\vep|\lambda_n|^{-2}\biggr)
&\leq&\frac{4|\lambda
_n|^{4}}{\vep^2} \Var\bigl(\sigma_2\bigl(\bfitx^{(n)}\bigr)\bigr) \nonumber\\[-4pt]\\[-12pt]
&=&o(1)\nonumber
\end{eqnarray}
by~(\ref{aim-a2}). Thus,~(\ref{aim-a}) follows from (\ref
{aim-a1}) and~(\ref{aim-a2}).\vadjust{\goodbreak}

To prove~(\ref{aim-a1}), we apply Lemma~\ref{lem-mean-sigma2},
in the setting that
%
%
\begin{equation} \nu_n(\lambda_n)=\nu_n(1+\lambda_n \ell_n
n^{-2\rho})=1+\lambda_n \ell_n n^{-2\rho}+ \zeta n^{-\eta
}+o(n^{-\eta}),
\end{equation}
so that, by Lemma~\ref{lem-mean-sigma2}(i),
%
%
\begin{eqnarray} \expec\biggl[\sum_{i\in[n]} w_i \weight(i)\biggr]&=&\frac{\sum
_{i\in[n]} w_i^2}{1-\nu_n(\lambda_n)} +o(n^{2\rho}|\lambda
_n|^{-2})\nonumber\\
&=&\nu_n(\lambda_n)\ell_n \bigl(|\lambda_n| \ell_n n^{-2\rho
}-\zeta n^{-\eta}+o(n^{-\eta})\bigr)^{-1}\nonumber\\
&&{} +o(n^{2\rho}|\lambda
_n|^{-2})\\
&=&|\lambda_n|^{-1} n^{2\rho}+\frac{\zeta}{\expec
[W]} n^{2\rho}|\lambda_n|^{-2} \nonumber\\
&&{}+o(|\lambda_n|^{-2} n^{2\rho}),\nonumber
\end{eqnarray}
which proves~(\ref{aim-a1}) with $\beta=-\zeta/\expec[W]$.

By Lemma~\ref{lem-mean-sigma2}(ii),
%
%
\begin{eqnarray}\qquad
\Var\biggl(\sum_{i\in[n]} w_i \weight(i)\biggr)&\leq& C
\biggl(\expec[w_T]^4\frac{1}{\ell_n} \sum_{i\in[n]} w_i^4+\expec
[w_T]^2\expec[w_T^2]\frac{1}{\ell_n} \sum_{i\in[n]}
w_i^3\biggr)\nonumber\\[-8pt]\\[-8pt]
&=&o(n^{4\rho}\lambda_n^{-4}),\nonumber
\end{eqnarray}
precisely when both terms in the middle inequality satisfy this bound.
We complete the proof by checking these estimates.
The first contribution is bounded by
%
%
\begin{equation} \frac{1}{\ell_n(1-\nu_n(\lambda_n))^4} \sum_{i\in
[n]} w_i^4 \leq\frac{C}{|\lambda_n|^4} n^{4\alpha+3\eta-1}
=o(n^{4\rho}|\lambda_n|^{-4}),
\end{equation}
since $4\alpha+3\eta-1=2(\tau-2)/(\tau-1)=2\rho<4\rho$.
The second contribution, instead, is bounded by
%
%
\begin{equation}
\frac{1}{\ell_n^2(1-\nu_n(\lambda_n))^5} \biggl(\sum_{i\in
[n]} w_i^3\biggr)^2 \leq\frac{C}{|\lambda_n|^5} n^{6\alpha+5\eta
-2} =o(n^{4\rho}|\lambda_n|^{-4}),
\end{equation}
since $6\alpha+5\eta-2=(3\tau-7)/(\tau-1)<4\rho=4(\tau-2)/(\tau-1)$.
This proves the required concentration for $\sigma_2(\bfitx^{
(n)})$ and hence
completes the proof of Theorem~\ref{thm-mult-coal}
for cluster weights and for any \textit{fixed} $\lambda$.

\subsection*{Convergence of the finite-dimensional distributions
random graph multiplicative coalescent}
So far, we have proved the convergence of
$\mathbf{X}^{(n)}(|\lambda_n|+\lambda)$
for a fixed time $\lambda$. By~\cite{AldLim98}, Lemma 26,
there exists an eternal multiplicative coalescent
with the same marginal for every $\lambda$.
By the strong Feller property of multiplicative coalescents
proved in~\cite{Aldo97}, as well as~\cite{AldLim98}, Lemma 27,
the convergence of $\mathbf{X}^{(n)}(|\lambda_n|+\lambda_1)$
implies that the future finite-dimensional distributions
$(\mathbf{X}^{(n)}(|\lambda_n|+\lambda_l))_{l=1}^k$
converge in distribution to the finite-dimensional distributions
of the eternal multiplicative coalescent. This completes the
proof of the convergence of the finite-dimensional distributions
in Theorem~\ref{thm-mult-coal} for cluster weights.

\subsection*{Convergence of cluster sizes from cluster weights}
By the adaptation of Theorem~\ref{thm-WC-3,4} to cluster weights
in Theorem~\ref{thm-scal-lim-cluster-weights}, we obtain
that $\weight_{\leq}(j)=|\Ccal_{\leq}(j)|(1+o_{\prob}(1))$,
so that the result immediately follows
for the cluster sizes.

\section*{Acknowledgments}

We thank Tom Kurtz for a discussion that helped us to simplify the
proof of Corollary~\ref{cor-conv-hit-times} substantially and to note
the extension to the convergence in the uniform topology in Remark
\ref{rem-conv-unif}. We thank David Aldous and Vlada Limic for help on
their results in~\cite{AldLim98}, which in particular clarified the
convergence of finite-dimensional distributions in Theorem~\ref{thm-mult-coal}.
We thank Sandra Kliem and two anonymous referees
for their valuable comments that helped us to substantially improve the
presentation.


%

%
\printaddresses

\end{document}